\documentclass[oneside,11pt]{book}
\usepackage{fullpage}
\usepackage{amsmath, amssymb,url}
\usepackage{amsthm}
\usepackage{textcomp}
\usepackage{graphicx}
\usepackage{latexsym, geometry}
\usepackage{color}
\usepackage[autostyle]{csquotes}
\usepackage[all]{xy}

\usepackage{enumitem}
\usepackage{mathtools}
\usepackage{ytableau}
\usepackage[enableskew,vcentermath]{youngtab}
\usepackage{tikz}

\newcommand{\unimodal}{\Delta}
\newcommand{\unimodless}{<_{\Delta}}
\newcommand{\Q}{\mathcal Q}
\newcommand{\J}{\mathcal J}

\newcommand{\TTT}{\mathcal T}
\newcommand{\symm}[1]{\mathfrak{S}_{#1}}
\newcommand{\uroot}[2]{\mathcal{R}_{#1}{\left(\symm{#2}\right)}}
\newcommand{\C}{\mathcal C}
\newcommand{\A}{\mathcal A}
\newcommand{\BBB}{{\mathcal{B}}}
\newcommand{\bbq}{\mathbb{Q}}
\newcommand{\F}{\mathcal F}
\newcommand\x{{\mathbf x}}
\newcommand\boxspan[1]{\boxplus\left\{#1\right\}}
\newcommand\spanbox[1]{{\left\langle #1 \right\rangle_\boxplus}}

\DeclareMathOperator\maj{\mathrm{maj}}
\DeclareMathOperator\imaj{\mathrm{imaj}}
\DeclareMathOperator\Des{Des}
\DeclareMathOperator\des{des}
\DeclareMathOperator\Inv{\mathrm{Inv}}
\DeclareMathOperator\inv{\mathrm{inv}}
\DeclareMathOperator\cDes{\mathrm{cDes} }
\DeclareMathOperator\CDes{\mathrm{CDES} }

\def\NN{{\mathbb N}}

\def\QSym{\mathcal{Q}}

\newcommand{\Pmod}[1]{\ (\mathrm{mod}\ #1)}

\newcommand{\ch}{{\operatorname{ch}}}
\newcommand{\Ct}{{Ct}}
\newcommand{\type}[1]{\mathbf{type}\,#1} 
\newcommand{\shape}[1]{\mathbf{sh}\,#1}
\newcommand{\RSK}[3]{#1\xrightarrow{RSK}\left({#2} ,{#3}\right)}

\def\SYT{{\rm SYT}}
\def\SSYT{{\rm SSYT}}

\theoremstyle{plain}
\newtheorem{theorem}{Theorem}[section]
\newtheorem{proposition}[theorem]{Proposition}
\newtheorem{lemma}[theorem]{Lemma}
\newtheorem{corollary}[theorem]{Corollary}
\newtheorem{conjecture}[theorem]{Conjecture}
\newtheorem{problem}[theorem]{Problem}
\newtheorem*{theorem*}{Theorem}
\newtheorem*{definition*}{Definition}
\newtheorem*{problem*}{Problem}
\newtheorem*{conjecture*}{Conjecture}
\theoremstyle{definition}
\newtheorem{definition}[theorem]{Definition}
\newtheorem{defn}[theorem]{Definition}
\newtheorem{example}[theorem]{Example}

\newtheorem{question}[theorem]{Question}
\newtheorem{remark}[theorem]{Remark}
\newtheorem{observation}[theorem]{Observation}

\numberwithin{section}{chapter}

\begin{document}
\begin{titlepage}
	\begin{center}
		
		\thispagestyle{empty}
		
		\vspace*{+7ex}
		{\fontsize{30}{38}\selectfont \bf Schur-positive Sets \\\vspace{+0.7ex}
		}
		\vspace{+16ex}
		
		{\LARGE Yuval Khachatryan-Raziel}\\
		\vspace{+20ex}
		{\Large Department of Mathematics}\\
		
		
		\vspace{+24ex}
		{\Large Ph.D.~Thesis\\Submitted to the Senate of Bar-Ilan University} \\
		
		\vspace{+1ex}
		\vspace{+10ex}
		Ramat-Gan, Israel \hspace{+40ex} {October, 2021}
	\end{center}
\end{titlepage}
\thispagestyle{empty}
\begin{center}
	This work was carried out under the guidance and supervision of \\
	Prof. Ron Adin and Prof. Yuval Roichman\\ 
	Department of Mathematics\\
	Bar-Ilan University
\end{center}
\pagebreak

\pagenumbering{gobble}
\chapter*{Acknowledgments}
I wish to express my deepest gratitude to several individuals whose contributions have been indispensable in the completion of this journey.

Foremost, my sincere thanks go to my academic supervisors, Ron Adin and Yuval Roichman, who guided me from the midst of my undergraduate studies up to the completion of this PhD. They have been extremely supportive and patient, fostering my growth both as a person and a mathematician. Their help, support and patience have been nothing short of unparalleled, and far exceeded my expectations.

I would also like to extend my appreciation to the professors and faculty members at Bar-Ilan University, whose contributions were pivotal to my development as a mathematician. Special thanks go to Uzi Vishne, Michael Schein, Nathen Keller, Reuven Cohen, and Gidi Amir. A note of profound gratitude must be made to the late Shmuel Dahari, who encouraged me to pursue an advanced degree and infused me with confidence. His readiness to answer every question will never be forgotten.

Special mention must also go to the participants of the Combinatorics Seminar at Bar Ilan University. I am particularly grateful to Eli Bagno, Yonah Cherniavsky, and Chaya Keller for their advice and help.

The journey would have been impossible without my friends from studies. Special gratitude is extended to Achya Baron, Tamar Baron, Ariel Weizman, Shira Gilat, and Adi Ben-Tsvi. A shout out to Tomer Bauer for his invaluable assistance with TeX and Sage.

To the secretaries - Malka Bachar, Yael Madar, and Michal Davidov - for their willingness to provide help whenever needed.

My deepest thanks go to my parents, Vardan and Varsenik Khachatryan, and my brother, Edward Khachatryan, for their unwavering support and belief in me.

I wish to express my deepest gratitude to my wife's family, especially my parents-in-law, Naomi and Amram Mashraky, for their help and support during my studies.

Lastly, but most importantly, I would like to express my deepest gratitude to my wife, Alona. Her steadfast support and patience throughout these long years of study have been a beacon of hope and endurance. To our children, Elisheva Miryam, David, Tamar, Ruth Tchiya, and Yochai Meir, your presence has filled my life with joy and purpose, motivating me to strive for success.

To anyone I may have inadvertently overlooked, please accept my sincere apologies.

Lastly, I dedicate my thesis to the memory of my grandparents Miryam and Edward Zilbershtein, and my uncle Andrei Zilberstein, who inspired me to foster a love for mathematics and science from a young age.

\thispagestyle{empty}

\tableofcontents
\mainmatter

\pagestyle{plain}
\pagenumbering{roman}
\chapter*{Abstract}
\addcontentsline{toc}{chapter}{Abstract}
The notion of descent set for permutations and standard Young tableaux (abbreviated to SYT) is classical. A significant problem in contemporary algebraic combinatorics is determining whether the associated quasisymmetric function is symmetric and Schur-positive.
Gessel proved that Knuth classes of permutations and inverse descent classes are Schur-positive.
Other examples of Schur-positive permutation sets follow from seminal works by Foata-Schützenberger and Gessel-Reutenauer: permutations of fixed Coxeter length, conjugacy classes, and $d$-roots of unity. 

Sometimes it is possible to define descents on maximal chains of certain posets using an $EL$-labeling
of the edges in the Hasse diagram. A classic example of a Schur-positive set of this type is the set of all maximal chains in the non-crossing partition lattice of type $A$, which was given by Stanley. In this case, the maximal chains can be realized by edge-labeled trees. A long-standing problem is to find other Schur-positive sets of maximal chains or trees.

In Chapter 3, we present an appealing set of maximal chains in the non-crossing
partition lattice (equivalently: a set of edge-labeled trees) that is Schur-positive:
we show that the subset of maximal chains in the non-crossing partition lattice, which may be realized 
by convex caterpillars, is Schur positive. We prove the following.

\begin{theorem*}[\ref{theorem-main}]
	The set of linearly ordered factorizations of the $n$-cycle $(1,2,\dots,n)$
	satisfies 
	$$\Q(\mathcal{L}_n)=\sum_{k=0}^{n-1} (k+1)s_{(n-k,1^k)},$$
	where the descent set of any $u\in \mathcal{L}_n$ is defined by the edge labeling of~\cite{AR_OMCNCPL}.
	In particular, $\mathcal{L}_n$ is Schur-positive.
\end{theorem*}

The proof uses the Goulden-Yong partial order on the edges together with 
a variant of an $EL$-labeling due to Björner. 

\medskip

The rest of our work studies {\em cyclic descent extensions} (CDE) of Schur-positive sets of permutations. The study of cyclic descent extensions has been an active research area in recent years.
Klyachko and Cellini introduced a natural notion of cyclic descent set for permutations. This notion was further studied by Dilks, Petersen, Stembridge, and many others.
An analogous notion for standard Young tableaux, using an axiomatic approach, was presented by Rhoades, Adin-Reiner-Roichman, and others. 

A significant problem in this field is whether a given set of permutations admits a cyclic descent extension.
The existence of cyclic extensions is related to significant problems in various fields, in particular, to the non-negativity of the coefficients of the Postnikov toric Schur functions and a long-standing open problem of Thrall on higher Lie characters. Adin, Hegedüs and Roichman, in \cite{AdinHegedusRoichman}, gave an algebraic characterization for Schur-positive sets that admit a CDE.

In Chapter 4, we study an important Schur-positive subset of the symmetric group, the set of all permutations of fixed Coxeter length.
We obtain sufficient conditions for the positive integer $m$, under which there is a CDE for permutations of length $m$ in the symmetric group $\symm{n}$, for almost all $n$.
The result involves surprising connections to generalized pentagonal numbers.
We prove the following.
\begin{theorem*}[\ref{theorem-mainCDesFixedInversions}]
	For $k < n$, the set of permutations in $\symm{n}$ with $k$ inversions has a cyclic descent extension if and only if $k$ is not a generalized pentagonal number.
\end{theorem*}

In Chapter 5, we state a conjectured necessary and sufficient condition for $d$-roots of unity in $\symm{n}$ to carry a cyclic descent extension. 
\begin{conjecture*}
	The set of roots of unity of order $d$ in $\symm{n}$ has a cyclic descent extension if and only if $d$ and $n$ are not coprime.
\end{conjecture*}

We prove the following results, towards this conjecture.

\begin{theorem*}[\ref{theorem-mainUrootPrimeVersion}]
	The set of roots of unity of order $d$ in $\symm{n}$:
	\begin{enumerate}
		\item Does not have a cyclic descent extension if $n$ and $d$ are coprime.
		\item Has a cyclic descent extension when $d$ is a prime power and $d$ and $n$
		are not coprime.
	\end{enumerate}
\end{theorem*}

\pagebreak
\pagenumbering{arabic}
\setcounter{page}{1}

\chapter{Introduction}
Given a set $A$ of combinatorial objects equipped with a descent function $Des:A\to 2^{[n-1]}$, define the quasi-symmetric function 
$$\Q(A) \coloneqq \sum_{a\in A}\F_{n,\Des(a)}$$
where $\F_{n,D}$ are Gessel's {\em fundamental quasi-symmetric functions}; see Subsection~\ref{subsection:QuasisymmetricFunctions} for more details. The following long-standing problem was posed by Gessel and Reutenauer in \cite{Gessel-Reutenauer}.
\begin{problem*}
	For which sets $A$ is $\Q(A)$ symmetric?
\end{problem*}

A symmetric function is called {\em Schur-positive} if all the coefficients in its expansion in the basis of Schur functions are non-negative; see Subsection~\ref{subsection: SymmetricFunctions} for definitions.
Determining whether a given symmetric function is Schur-positive is a significant problem in contemporary algebraic combinatorics; see Stanley's article~\cite{Stanley_problems}.
We say that the set $A$ of combinatorial objects is Schur-positive if $\Q(A)$ is Schur-positive; for details and examples, see Sections \ref{subsection: SymmetricFunctions} and~\ref{subsection: SchurPositiveExamples}.

Another approach to Schur-positivity is via a significant problem in combinatorial representation theory: Provide explicit formulas that express the character values of interesting representations as a weighted enumeration of combinatorial objects. The reader is referred to the survey article by Barcelo and Ram~\cite{BarceloRam}. Many such formulas have been found; examples include the Murnaghan-Nakayama rule \cite[Section 7.17]{EC2}, characters of Specht modules of ribbon shapes~\cite{Gessel}, Lie characters~\cite{Gessel-Reutenauer}, Gelfand models~\cite{AdinPostnikovRoichman2}, characters of homogeneous components of the coinvariant algebra~\cite{AdinPostnikovRoichman1}, characters induced from a lower ranked exterior algebra~\cite{ElizaldeRoichmanArc} and k-root enumerators~\cite{RoichmanKRoots}, among others. It turns out that the sets that appear in these formulas are Schur-positive sets.
An abstract framework for this phenomenon was proposed in~\cite{AR1}. In particular, it was shown in~\cite{AR1} that the equality of the characters of two Schur-positive sets is equivalent to the equidistribution of their descent sets, implying the equivalence of classical theorems in representation theory, for example, the equivalence of theorems of Lusztig-Stanley in invariant theory~\cite{Lusztig-Stanley} and Foata-Schützenberger in permutation statistics~\cite{FoataSchutz}. Their work is a development of the ideas of Gessel and Reutenauer in~\cite{Gessel-Reutenauer}. 

Note that equivalent - but different terminology is used in the literature in this field. In particular,  
Schur-positive sets are also called {\em fine sets}. In this work, we use only the term Schur-positive.

Schur-positive sets are often subsets of the symmetric group $\symm{n}$. However, there are important examples that are not of this form. In particular, a part of our work deals with Schur-positive sets of maximal chains in posets. Maximal chains in a labeled poset $P$ are equipped with a natural descent map. A well-known conjecture of Stanley~\cite[III, Ch. 21]{Stanley_thesis} 
implies that all examples of Schur-positive labeled posets, in this sense
correspond to intervals in the Young lattice.

Another way to equip the set of maximal chains with a descent map is by labeling the edges in the Hasse diagram.
A classic example of a Schur-positive set 
of this type, 
the set of all maximal chains in the non-crossing partition lattice of type $A$,
was given by Stanley~\cite{StanleyPFnNCPL}. 
An $EL$ edge-labeling of this poset was presented in an earlier work of Björner~\cite{Bjorner}; see also~\cite{BW, Mcnamara, AR_OMCNCPL}.

 Adin and Roichman presented in \cite{AR_OMCNCPL} a variant of the {\em EL}-labeling that was defined by Björner, using the fact that every maximal chain bijectively corresponds to a factorization of the cycle $(1\ \dots\  n)$ in $\symm{n}$ into a product of $n-1$ transpositions. This variant will be used in our work.

Goulden and Yong in \cite{Goulden-Yong} showed that every maximal chain in the non-crossing partition lattice can be interpreted as a geometric non-crossing tree with an ordering of the edges (see Chapter \ref{chapter: ConvexCaterpillars} for more details). In particular, their construction defines a poset structure on the edges of a non-crossing geometric tree. 
It is natural to ask when this order is a {\em linear order}.
In Chapter~\ref{chapter: ConvexCaterpillars}, we provide a detailed description of these trees, named {\em convex caterpillars}. We show that the set $\mathcal{L}_n$ of maximal chains that correspond to these trees is Schur-positive and calculate the corresponding quasisymmetric function.
We prove the following result:
\begin{theorem*}[\ref{theorem-main}]
	The set of linearly ordered factorizations of the $n$-cycle $(1,2,\dots,n)$
	satisfies 
	$$\Q(\mathcal{L}_n)=\sum_{k=0}^{n-1} (k+1)s_{(n-k,1^k)},$$
	where the descent set of any $u\in \mathcal{L}_n$ is defined by the edge labeling of~\cite{AR_OMCNCPL}.
	In particular, $\mathcal{L}_n$ is Schur-positive.
\end{theorem*}
Here and elsewhere, $s_\lambda$ denotes the Schur function indexed by a partition $\lambda$. See Subsection~\ref{subsection: SymmetricFunctions}.

The rest of our work is dedicated to studying {\em cyclic descent extensions} for several classes of Schur-positive sets in $\symm{n}$. The study of cyclic descent extensions has been an active research area in recent years.

One of the first works studying cyclic descents was a paper by Cellini~\cite{Cellini}, who introduced a natural notion of cyclic descent on permutations as follows:
$${\rm CDES}(\pi) \coloneqq \left\{1 \leq i \leq n: \pi_i > \pi_{i+1}\right\}\subseteq  \left[n\right]\ \text{for all}\ \pi\in\symm{n}.$$
with the convention $\pi_{n+1} \coloneqq \pi_1$.
A more restricted notion of a cyclic descent number was used by Klyachko in \cite{Klyachko-Cyclic}. Cyclic descent sets were further studied by Dilks, Petersen and Stembridge in \cite{Dilks_at_all_cylclic}.

A well-established notion of descent also exists for standard Young tableaux (SYT); however, there is no apparent cyclic analog. A significant milestone on the way to defining a cyclic descent set on SYT was the work of Rhoades~\cite{Rhoades_CDes}, who described a notion of cyclic descent set for SYT of rectangular shape. Some properties common to Cellini's definition for permutations and Rhoades' construction for SYT appeared in other combinatorial settings \cite{PetersenPyl_at_all_CDes, Pechnik, ElizaldeRoichmanCyclic, AdinElizaldeRoichmanCDes}. These properties led to an axiomatic definition of a cyclic descent extension in \cite{ARR}:
\begin{definition*}[\ref{definition-CyclicDescent}]
	Let $\TTT$ be a finite set, equipped with a {\em descent map} 
	$\Des: \TTT \longrightarrow 2^{[n-1]}$. 
	A {\em cyclic extension} of $\Des$ is
	a pair $(\cDes,p)$, where 
	$\cDes: \TTT \longrightarrow 2^{[n]}$ is a map 
	and $p: \TTT \longrightarrow \TTT$ is a bijection,
	satisfying the following axioms: for all $T$ in $\TTT$,
	\[
	\begin{array}{rl}
		\text{(extension)}   & \cDes(T) \cap [n-1] = \Des(T),\\
		\text{(equivariance)}& \cDes(p(T))  = {\rm sh}(\cDes(T)),\\
		\text{(non-Escher)}  & \varnothing \subsetneq \cDes(T) \subsetneq [n].\\
	\end{array}
	\]
\end{definition*}
Here $\mathrm{sh}$ denotes the cyclic shift of a subset of $[n]$, as in Definition~\ref{definition-CyclicShift}.

Cellini's and Rhoades' constructions are special cases of Definition~\ref{definition-CyclicDescent}. A naturally arising question is the existence of such an extension for specific subsets of the symmetric group $\symm{n}$(see Chapter~\ref{chapter: CyclicDescentInvDesClass} and the discussion at the beginning of \cite{AdinHegedusRoichman}).
It is also not apparent whether a cyclic descent set similar to Rhoades' can be defined for the set of SYT of every shape.
The question for SYT was solved by Adin, Reiner and Roichman in~\cite{ARR}. They proved that a cyclic descent extension exists for the set of SYT of skew shape $\lambda / \mu$ if and only if $\lambda / \mu$ is not a connected ribbon. They used Postnikov's toric symmetric functions. A constructive proof was found later by Huang in \cite{Huang}.

Adin, Hegedüs and Roichman provided in \cite{AdinHegedusRoichman} an algebraic necessary and sufficient condition for the existence of a cyclic descent extension on Schur-positive sets. In addition, they characterized the conjugacy classes in $\symm{n}$ that possess a cyclic descent extension. Their proof was not constructive - providing a combinatorial construction of a cyclic descent extension for conjugacy classes is a significant work in progress. An explicit construction for involutions with fixed points was described in \cite{AR_CdesInvolutions}.
Explicit constructions of a cyclic descent extension for Dyck paths appear in \cite{PetersenPyl_at_all_CDes}. This construction was extended to Motzkin paths by B.\ Han in \cite{BinHanCDes}.
Another recent development is a cyclic analog of quasi-symmetric functions by Adin, Gessel, Reiner and Roichman in~\cite{AGRR_CyclicQSym}.

In Chapter~\ref{chapter: CyclicDescentInvDesClass}, we study the existence of cyclic descent extensions on inverse descent classes and their generalizations.
In particular, we obtain the following result:
\begin{theorem*}[Theorem \ref{theorem-mainInvDevCDes}]
	The following sets in $\symm{n}$ have cyclic descent extensions:
	\begin{enumerate}
		\item $D_{2^J, n}^{-1}$ if and only if  $J\neq \emptyset$.
		\item $D_{[I,J],n}^{-1}$ for an interval 
		$$[I, J]=\left\{K: I\subseteq K\ and \ K\subseteq J\right\}$$ (in the poset $2^{[n-1]}$ ordered by inclusion),
		if and only if $I\subsetneq J$.
		\item For every maximal chain $C$ in the poset $[I,J]$, $D_{C, n}^{-1}$ has a cyclic descent extension if and only if  $I\subsetneq J$ and $2 \not| (|J| - |I|)$.
		In particular, $D_{C, n}^{-1}$ for a maximal chain $C$ in $2^{[n-1]}$ has a cyclic descent extension if and only if  $n-1$ is odd, or equivalently $n$ is even.		
	\end{enumerate}
\end{theorem*} 
 Here  $D_{S, n}^{-1}$ denotes, for a set $S$ of subsets of $[n]$, a {\em generalized inverse descent class}, as in Definition~\ref{definition:GeneralizedInvDesClass}. 

In the last part of Chapter~\ref{chapter: CyclicDescentInvDesClass}, we use this result to describe conditions for the existence of cyclic descent extensions for permutations with a given number of inversions.  We determine necessary and sufficient conditions for permutations with a given inversion number in $\symm{n}$ for sufficiently large $n$. Our result has a surprising connection to Euler's pentagonal number theorem.
We prove the following.
\begin{theorem*}[\ref{theorem-mainCDesFixedInversions}]
	For  $k < n$, the set of permutations in $\symm{n}$ with $k$ inversions has a cyclic descent extension if, and only if, $k$ is not a generalized pentagonal number.
\end{theorem*}

In the last chapter, we study roots of unity in $\symm{n}$ and conditions for the existence of a  cyclic descent extension, which can be seen as a generalization of the criterion given in \cite{AdinHegedusRoichman}.

In particular, we obtain the following result:
\begin{theorem*}[\ref{theorem-mainUrootPrimeVersion}]
The set of roots of unity of order $d$ in $\symm{n}$:
\begin{enumerate}
	\item Does not have a cyclic descent extension if $n$ and $d$ are co-prime.
	\item Has a cyclic descent extension when $d$ is a prime power and $d$ and $n$ are not co-prime.
\end{enumerate}
\end{theorem*}

Finally, we present the following conjecture, a subject for future work.
\begin{conjecture*}
The set of roots of unity of order $d$ in $\symm{n}$ has a cyclic descent extension if and only if $d$ and $n$ are not co-prime.
\end{conjecture*}

\chapter{Background}\label{section:background}

\section{Basic notation and definitions}
This section quickly reviews some notation and basic facts about permutations, partitions, and related structures.

We use $[n]$ to denote the set $\left\{1,\dots,n\right\}$. Unless explicitly stated, the arithmetic on $[n]$ is performed modulo $n$.
For a set $X$, the {\em power set of $X$}, consisting of all subsets of $X$, will be denoted  $2^X$. For sets $X$ and $Y$, the set of all functions from $X$ to $Y$ is denoted $Y^X$.

\subsection{Compositions and partitions}
\begin{definition}
	A {\em weak composition} of $n$ is a sequence $\alpha = \left(\alpha_1,\alpha_2,\dots\right)$ of non-negative integers such that $\sum_{k=1}^\infty\alpha_k = n$.
\end{definition}
\begin{definition}
	A {\em partition} of $n$ is a weakly decreasing sequence of non-negative integers $\lambda = (\lambda_1,\lambda_2,\dots)$ such that $\sum_{k=1}^\infty\lambda_k = n$. We denote $\lambda \vdash n$.
\end{definition}
\begin{definition}
	The {\em length} of a partition $\lambda = \left(\lambda_1,\lambda_2,\dots \right)$, denoted $l(\lambda)$, is the number of non-zero parts $\lambda_i$.
	The non-zero $\lambda_i$ are called the {\em parts} of $\lambda$.
\end{definition}
\begin{remark}
	Sometimes we use the notation $\lambda = (a_1^{i_1},\dots,a_k^{i_k})$
	where $i_k$ is the number of parts of $\lambda$ equal to $a_k$.
\end{remark}

\subsection{Permutations}
\begin{definition}
	The {\em symmetric group} or the group of {\em permutations} is the group of bijections from $[n]$ to itself. We denote the symmetric group $\symm{n}$.
\end{definition}

For a permutation $\pi\in\symm{n}$ we often write $\pi_i$ instead of $\pi(i)$.
In some cases, we will use cycle decomposition representation, which we now describe. 
\begin{definition} A cycle $\pi=\left(i_1\ \dots\ i_l\right)$ is a permutation in $\symm{n}$ such that $\pi(i_k) = i_{k+1}$ for all $1\leq k \leq l-1$, $\pi(i_l)=i_1$ and $\pi(m)=m$ if $m\notin\left\{i_1,\dots,i_l\right\}$. The number $l$ is called the length of the cycle. We say that cycles $(i_1,\dots,i_l)$ and $(i^\prime_1,\dots,i^\prime_{l^\prime})$ are \em{disjoint} if $i_{k}\neq{j_{k^\prime}}$ for all $1\leq k \leq l$ and $1 \leq k^\prime \leq l^\prime$. 
\end{definition}

Recall the basic fact from group theory.
\begin{theorem}\label{theorem:permutation_cylce_type}
	Every $\pi \in \symm{n}$ can be represented as a product of disjoint cycles $$\pi=\left(a_1\ \dots \ a_{l_1}\right)\left(a_{l_1+1}\dots a_{l_1+l_2} \right)\dots\left(a_{l_1 + \dots + l_{k-1}+1}\ \dots \ a_n\right).$$
	Furthermore:
	\begin{enumerate}
		\item The representation is unique, up to the order of the cycles.
		\item Two permutations $\pi, \sigma \in \symm{n}$ with $$\pi=\left(a_1 \dots \right)\dots\left(\dots a_n\right)$$ are conjugate in $\symm{n}$ if and only if there is $\tau\in\symm{n}$ such that 
		$$\sigma = \left(\tau(a_1)\dots\right)\dots\left(\dots\tau(a_n)\right).$$
		In such case we have $\sigma = \tau\pi\tau^{-1}$.
	\end{enumerate}
\end{theorem}

\begin{definition}
	The {\em cycle type} of a permutation
	$$\pi=\left(a_1\ \dots \ a_{l_1}\right)\left(a_{l_1+1}\dots a_{l_1+l_2} \right)\dots\left(a_{l_1 + \dots + l_{k-1}+1}\ \dots \ a_n\right)$$ is the partition of $n$ obtained by writing the lengths of cycles in the cycle decomposition of $\pi$ in weakly decreasing order. The cycle type of $\pi$ is denoted by $\type{\pi}$. 
\end{definition}

Part 2 of theorem \ref{theorem:permutation_cylce_type} is equivalent to the statement that two permutations in $\symm{n}$ are conjugate if and only if they have the same cycle type. The conjugacy class of all permutations $\pi$ with $\type{\pi}=\lambda$ will be denoted by $\C_\lambda$.

\section{Representations}
Representations and their characters play an essential role in understanding objects' symmetries and studying the structures of groups. We recall the basic definitions and theorems we need  from representation theory.
\begin{definition}
	An {\em action} of a group $G$ on a set $X$ is a function from $G\times X$ to $X$, written as $g\cdot x$ or simply $gx$ for all $g\in G$ and $x\in X$ with the following properties:
	\begin{enumerate}
		\item $g_1\cdot(g_2 \cdot x) = (g_1 g_2)\cdot x$ for all $g_1, g_2 \in G$ and $x \in X$.
		\item For the identity element $e$ of $G$ we have $e\cdot x = x$ for all $x\in X$.
	\end{enumerate}
\end{definition}

For a vector space $V$, we often identify a linear map $T: V \to V$ with its representation matrix relative to the standard basis. Let $GL(V)$ denote the general linear group of invertible linear transformations from $V$ to itself. For a field $\mathbb{F}$, let $GL_n(\mathbb{F})$ denote group of invertible matrices of size $n\times n$ with elements in $\mathbb{F}$. If $V$ is an $n$-dimensional vector space over $\mathbb{F}$, then $GL(V)$ and $GL_n(\mathbb{F})$ are isomorphic; therefore, unless stated explicitly, no distinction will be made between them.

\begin{definition}
	Let $G$ be a finite group, $\mathbb{F}$ a field, and $V$ a vector space over $\mathbb{F}$.
	A {\em  linear representation} of $G$ is a group homomorphism 
	$$\rho : G \to GL(V).$$ 
\end{definition}

\begin{observation}
	Every linear group representation $\rho:G\to GL(V)$ defines an action of $G$ on $V$
	by $g\cdot v = \rho(g)v$.
\end{observation}

\begin{definition}
For a set $X$ and field $\mathbb{F}$,  the {\em free vector space $\mathbb{F}^{(X)}$ generated by $X$ over $F$} is the direct sum of copies of $F$ indexed by $X$
$$\mathbb{F}^{(X)}=\bigoplus_{x\in X} \mathbb{F}$$ which is the subspace of the cartesian product 
$\prod\limits_{x\in X}F$ of elements with finitely many nonzero components. It can be equivalently viewed as the vector space of formal linear combinations $\sum\limits_{x\in X} a_x x$ with coefficients in $F$ with finitely many $x\in X$ with non-zero $a_x$. If we replace field $\mathbb{F}$ in the definition with a commutative ring $R$, then $R^{(X)}$ {\em is the a $R$-module over $X$.}
\end{definition}
If the set $X$ in the previous definition is the group $G$ itself, then $R^{(G)}$ can be given a ring structure as follows. 
\begin{definition}
	Let $G$ be a finite group and $R$ a commutative ring. Then $R\left[G\right]$, the {\em group ring} or the {\em group algebra} of $G$ over $R$ is the free $R$-module over $G$, together with multiplication between two elements defined as follows:
	$$\sum_{g\in G}a_g g \sum_{h\in G} b_h h = \sum_{k\in G}\left(\sum_{gh=k} a_g b_h\right)k.$$
\end{definition}

Note that for every finite group $G$ and finite set $X$, if $G$ acts on $X$, then this action can be extended to an invertible linear transformation from $\mathbb{F}^{(X)}$ to itself by linear extension of the action, namely:
$$g\cdot\left(\sum_{x\in X} a_x x\right)=\sum_{x\in X} a_x (g\cdot x).$$ 
In turn, if $V$ is a vector space over the field $\mathbb{F}$, then every group representation $\rho: G \to GL(V)$ turns the vector space $V$ into ab $\mathbb{F}[G]$ module with elements of $\mathbb{F}\left[G\right]$ acting on an element of $V$ by
$$\left(\sum_{g\in G}a_g g\right)v=\sum_{g\in G} a_g\rho(g)v.$$ 
On the other hand, every module $V$ over $\mathbb{F}[G]$ is a vector space, and the action of $\mathbb{F}[G]$ is a linear transformation. Hence, when a module over $\mathbb{F}[G]$ is a finite-dimensional vector space over $\mathbb{F}$, it defines a matrix representation of $G$.

\section{Shapes and Standard Young Tableaux}\label{section:ShapesAndSYT}
Standard Young tableaux are a classical object widely studied in many branches of mathematics. Often, they are used to describe the irreducible representations of $\symm{n}$. This section provides the basic definitions that will be used later. Our exposition is mainly based on \cite{Sagan} and \cite{EC2}, which one can refer to for a more detailed treatment of the subject.

\begin{definition}
	Let $\lambda \vdash n$. The {\em Ferrer's diagram} or {\em shape} of $\lambda$ is  an array $n$ of cells, having $l(\lambda)$ left justified rows, such that row $i$ contains $\lambda_i$ cells for $1 \leq i \leq l(\lambda)$. 
\end{definition}
\begin{remark}
	In this text, we use the {\em English notation}, where row $i$ appears above the row $i+1$. Sometimes Ferrer's diagrams are called {\em Young diagrams}.
\end{remark}

The cell in row $i$ and column $j$ has coordinates $(i,j)$.

If $m \leq n$ and $\mu \vdash m$ and $\lambda \vdash n$, we denote $\mu \subseteq \lambda$ if $l(\mu) \leq l(\lambda)$ and $\mu_i \leq \lambda_i$ for every $1 \leq i \leq l(\mu)$.

\begin{definition}
	Let $m \leq n$, $\mu \vdash m$, $\lambda \vdash n$. The {\em skew shape} of $\lambda / \mu$ is an array of $n-m$ cells obtained by removing the first $\mu_i$ cells from row $i$ in the shape of $\lambda$. If $m=0$, the skew shape $\lambda / \mu$ is called normal or just a shape, and we write $\lambda$ instead of $\lambda/\mu$.         
\end{definition}

\begin{definition}                 
	A {\em Young tableau} of shape $\lambda$ is an array $T$ obtained by filling the cells of the Young diagram of $\lambda$ by the numbers $1,\dots,n$ bijectively.
	A young tableau T of shape $\lambda$ is also called a {\em $\lambda$-tableau} and is denoted $T^\lambda$. We also write $\shape{T} = \lambda$. 
\end{definition}

\begin{definition}
	The tableau $T$ is a {\em standard Young tableau} if it is a Young tableau and every row read left to right and column read from top to bottom of $T$ is an increasing sequence. The set of all standard Young tableaux of shape $\lambda$ will be denoted $\SYT(\lambda)$. The size of $\SYT(\lambda)$ will be denoted  $f^\lambda$.
\end{definition}

\begin{definition}
	A semistandard Young tableau of shape $\lambda$ is an array $T$ obtained by replacing the cells of shape $\lambda$ with positive integers, such that 
	the rows of $T$ are non-decreasing sequences and the columns of $T$ are increasing sequences. The set of all semistandard Young tableaux of shape $\lambda$ will be denoted $\SSYT(\lambda)$.
\end{definition}

\begin{definition}\label{definition:TypeSSYT}
	For a standard or semistandard Young tableaux $T$, the {\em content} or {\em type} of $T$ is the sequence $\type{T} = \mu=(\mu_1,\dots,)$ where $\mu_i$ is equal to the number of times that $i$ appears in $T$. The set of semistandard Young tableaux of shape $\lambda$ and content $\mu$ is denoted with $\mathcal{T}_{\lambda\mu}^0$. The size of $\mathcal{T}_{\lambda\mu}^{0}$ is called the Kostka number $K_{\lambda\mu}$.
\end{definition}

\begin{remark}
	All the definitions involving shapes can be generalized to skew shapes by replacing $\lambda$ with $\lambda/\mu$ in the definition to obtain
	skew Young tableaux, standard skew Young tableaux and semistandard skew Young tableaux, $f^{\lambda / \mu}$, etc.  
\end{remark}

\section{Statistics on permutations and Young tableaux}
\begin{definition}\label{definition-CombinatrialStatistic}
Let $X$ be a set. A {\em statistic} on $X$ is a function $f: X \to A$ where $A \subseteq \NN\cup\left\{ 0 \right\}$ or $A \subseteq 2^{\NN\cup\left\{0\right\} }$.
\end{definition} 

The study of statistics on finite objects has a crucial role in combinatorics. The notion of equidistribution is crucial as it often allows us to unify questions that deal with seemingly different statistics. 
\begin{definition}\label{definition-equidistribution}
Let $f: X \to A$ be a statistic on the set $X$ and $g: Y \to A$ a statistic on the set $Y$. We say that $f$ and $g$ are {\em equidistributed} if for every $a\in A$ we have 
$$\left|\{x\in X:f(x)=a\}\right|=\left|\{y\in Y: g(x)=a\}\right|.$$
\end{definition}

\begin{remark}\label{remark: multisetStatistic}
	The notions of a statistic and a distribution of a statistic are naturally extended to multisets.
\end{remark}

 In the rest of this section, we describe some of the classical statistics on $\symm{n}$ and their generalizations on $\SYT(n)$ that will be needed later.

\begin{definition}
	Let $\pi\in\symm{n}$. We say that $i$ is a descent of $\pi$ if $\pi_i > \pi_{i+1}$. We denote the descent set of $\pi$ by
	$$\Des(\pi) \coloneqq \left\{i:\pi_i > \pi_{i+1}\right\}.$$
\end{definition}
The notion of descent can be extended to other objects. For example, it is easily extended to finite sequences of integers. In Chapter~\ref{chapter: ConvexCaterpillars} of this work, we study descents defined on non-crossing partitions of $[n]$.

Extension of descents on the standard Young tableaux plays a significant role in algebraic combinatorics and is defined as follows.

\begin{definition}
	The descent set of a  standard (skew) Young tableau $T$ is 
	$$\Des(T) \coloneqq \left\{i:\text{$i+1$ appears in a lower row of $T$ than $i$}\right\}.$$ 
\end{definition}

The Robinson-Schensted-Knuth correspondence, or RSK in short, is a bijection between permutations on $n$ and ordered pairs of standard Young tableaux of the same shape. It also connects the descents of permutations and standard Young tableaux.

\begin{theorem}\label{theorem-RSKProperties}\cite[Chapter 3]{Sagan}
	There exists a bijection 
	$$RSK:\symm{n}\to\biguplus_{\lambda\vdash n}\SYT(n)\times\SYT(n)$$
	such that if $\RSK{\pi}{P_\pi}{Q_\pi}$ then:
	\begin{enumerate}
		\item $\shape{P_\pi} = \shape{Q_\pi}.$
		\item $\RSK{\pi^{-1}}{Q_\pi}{P_\pi}.$
		\item $\Des(\pi)=\Des(Q_\pi)$ and $\Des(\pi^{-1})=\Des(P_\pi)$.		
	\end{enumerate}
\end{theorem}

\begin{definition}
	Let $X$ be a set with descent statistic $Des:X \to\left[n\right]$. For $x\in X$ the {\em descent number} of $x$ is the size of its descent set:
	$$\des(x) \coloneqq |Des(x)|.$$
\end{definition}

Inversions are similar to descents and can be seen as their generalization.
\begin{definition}
	Let $\pi\in\symm{n}$. We say that $(i,j)$ is an inversion of $\pi$ if $i<j$ and $\pi_i > \pi_j$. We denote the inversion set of $\pi$ with
	$$\Inv(\pi) \coloneqq \left\{(i,j):\text{$i<j$ and $\pi_i > \pi_j$}\right\}.$$
	The number of inversions of $\pi$ is called the {\em inversion number} of $\pi$ and is denoted with $\inv(\pi)$.
\end{definition}

A statistic closely related to descents and inversions is the major index.
\begin{definition}\label{definition:MajorAndInverseMajorIndex}
	Let $\pi\in\symm{n}$. The {\em major index} of $\pi$ is defined as the sum of its descents:
	$$\maj(\pi) \coloneqq \sum_{i \in \Des{\pi} } i.$$
	The {\em inverse major index} of $\pi$ is defined as the major index of its inverse.
	$$
	\imaj(\pi)=\maj(\pi^{-1}).
	$$	
\end{definition}

\begin{definition}\label{definition:CoxeterLength}
The standard {\em Coxeter generators} of $\symm{n}$ are the transpositions $s_i=(i\ i+1)$
for $1\leq i \leq n-1$. The minimal number of Coxeter generators required to express the permutation $\pi\in\symm{n}$ as a product of Coxeter generators is called the {\em Coxeter length} of  $\pi$.
\end{definition}
The connection between the length of a permutation and its inversion number is stated below.
\begin{theorem}\cite[Chapter 1]{BjornerCombinatoricsCoxeter}\label{thm:CoxeterLengthInversionNumber}
	For every $\pi\in\symm{n}$, $\inv(\pi)$ is equal to its Coxeter length.
\end{theorem}

\section{Symmetric and quasi-symmetric functions}\label{section:prel_quasi}

\subsection{Symmetric functions}
\label{subsection: SymmetricFunctions}
Let $\x = \left\{x_n\right\}_{n\in\mathbb{N}}$ be a set of commuting indeterminates. 

Symmetric and quasi-symmetric functions in ${\bf x}$ can be defined over various (commutative) rings of coefficients, including the ring of integers $\mathbb{Z}$; for simplicity we define it over the field $\bbq$ of rational numbers. 
For a sequence $\mu=\left\{\mu_n\right\}_{n\in\mathbb{N}}$, let 
$$\x^\mu = \prod_{n\in\mathbb{N}}x_n^{\mu_n}.$$
For formal power series $f(\x)=\sum\limits_{a\in\mathcal{A}}c_a x^a$, $[x^a]f$ will denote $c_a$,  the coefficient of $x^\alpha$ in $f$. 
\begin{definition}                                 
	A {\em symmetric function} in the variables $x_1,x_2,\ldots$ is a formal power series $f({\bf x})\in \bbq[[{\bf x}]]$, of bounded degree, 
	such that for any three sequences (of the same length $k$) of positive integers, $(a_1,\ldots,a_k)$, $(i_1,\dots,i_k)$ and $(j_1,\dots, j_k)$, the coefficients of $x_{i_1}^{a_1}\cdots x_{i_k}^{a_k}$ and of $x_{j_1}^{a_1}\cdots x_{j_k}^{a_k}$ in $f$ are the same:
	\[
	[x_{i_1}^{a_1}\cdots x_{i_k}^{a_k}]f = 
	[x_{j_1}^{a_1}\cdots x_{j_k}^{a_k}]f.
	\]
	A symmetric function $f$ is called {\em homogeneous} of degree $n$ if every for every $\alpha$ such that $\left[x^\alpha\right] \neq 0$, we have $\sum\limits_{i=1}^\infty\alpha_i = n$. The vector space of homogeneous functions of degree $n$ is denoted  $\Lambda_n$.
\end{definition}

It is easy to see, that if $f\in\Lambda_n$ and $\left[\x^\mu\right]f\neq 0$,
then by definition of symmetric function, $[\x^\lambda]f=[\x^\mu]f$ where $\lambda$ is the partition obtained by writing the parts of $\mu$ in non-increasing order.
Hence, every $f\in\Lambda_n$ can be written as a linear combination of 
$$m_\lambda = 
	\sum_{\substack{\left(x_{i_1},\dots,x_l\right) \\ 
	\text{$i_j \neq i_k$ for $j\neq k$} } }x_{i_1}^{\lambda_1}\dots x_{i_l}^{\lambda_l}$$
where $\lambda =\left(\lambda_1,\dots,\lambda_l\right)$  ranges over all partitions of $n$. 
It is also clear that $m_\lambda$ are independent, therefore 
$${m_\lambda:\lambda \vdash n}$$ is a basis o $\Lambda_n$. The function $m_\lambda$ is called {\em monomial symmetric function corresponding to $\lambda$.}
We will work with another basis for $\Lambda_n$. 

\begin{definition}
	For a semistandard Young tableau $T$, define $x^T$ to be $x^\mu$ where $\mu=\type{T}$ as in Definition~\ref{definition:TypeSSYT}.
\end{definition}
\begin{definition}
	{\em Schur functions}, indexed by partitions, are defined as follows:
	$$s_\lambda=\sum_{T\in \SSYT(\lambda)}x^T,$$
	where $\SSYT(\lambda)$ is the set of semistandard Young tableaux of shape $\lambda$, defined in Section~\ref{section:ShapesAndSYT}.
\end{definition}
\begin{remark}
	{\em Skew Schur functions} are defined analogously by replacing partitions with skew shapes in the definition.
\end{remark}

There is a dictionary relating symmetric functions to characters of the symmetric group $\symm{n}$. The irreducible characters of $\symm{n}$ are indexed by partitions $\lambda \vdash n$ and denoted $\chi^\lambda$. 
The {\em Frobenius characteristic map} $\ch$ from class functions on $\symm{n}$ to symmetric functions is defined by $\ch(\chi^{\lambda}) = s_{\lambda}$, and extended by linearity.
The Frobenius characteristic map has two properties that allow reducing the study of characters to the study of symmetric function.
First, it is an isometry between the space of class functions on $\symm{n}$ and $\Lambda_n$, where the inner product on characters is defined by $\left<\chi^\mu,\chi^\nu\right> = \delta_{\mu \nu}$ and the inner product on $\Lambda_n$ is defined by $\left<s_\mu,s_\nu\right> = \delta_{\mu\nu}.$
The other property is that if $\chi$ and $\mu$ are irreducible characters for $\symm{n}$ and $\symm{m}$ respectively, then $\ch(\chi\cdot\mu)=s_\mu s_\nu$ 
where $\chi\cdot\mu$ denotes the induced character $(\chi \otimes \mu)\uparrow^{\symm{n+m} }$.
A symmetric function $f$ is {\em Schur-positive} if in the expansion of $f$ in the basis of symmetric functions $$f=\sum_{\lambda\vdash n}a_\lambda s_\lambda$$
we have $a_\lambda \geq 0$ and $a_\lambda$ is an integer for every $\lambda \vdash n$.

\subsection{Quasisymmetric functions}\label{subsection:QuasisymmetricFunctions}
Quasisymmetric functions allow us to reduce questions about representations and symmetric functions to the study of descents.
The following definition of a quasi-symmetric function can be found in~\cite[7.19]{EC2}.
\begin{definition}                                 
	A {\em quasi-symmetric function} in the variables $x_1,x_2,\ldots$  is a formal power series $f({\bf x})\in \bbq[[{\bf x}]]$, of bounded degree, 
	such that for any three sequences (of the same length $k$) of positive integers, $(a_1,\ldots,a_k)$, $(i_1,\dots,i_k)$ and $(j_1,\dots, j_k)$, where the last two are {\em increasing}, the coefficients of $x_{i_1}^{a_1}\cdots x_{i_k}^{a_k}$ and of $x_{j_1}^{a_1}\cdots x_{j_k}^{a_k}$ in $f$ are the same:
	\[
	[x_{i_1}^{a_1}\cdots x_{i_k}^{a_k}]f = 
	[x_{j_1}^{a_1}\cdots x_{j_k}^{a_k}]f
	\]
	whenever $i_1 < \ldots < i_k$ and $j_1 < \ldots < j_k$. 
	The space of quasisymmetric functions is denoted  $\QSym$. A quasisymmetric function $f$ is homogeneous of degree $n$ if $[x^\alpha]f$ implies that $\sum\limits_{i=1}^\infty \alpha_i = n$. The space of homogeneous quasisymmetric functions of degree $n$ is denoted  $\QSym_n$.
\end{definition}
Every symmetric function is quasi-symmetric but not conversely. For example, the function $\sum_{i<j}{x_i^2 x_j}$ is quasi-symmetric but not symmetric.

The {\em monomial quasisymmetric functions} $M_\alpha$ are defined by 
$$M_\alpha=\sum_{i_1<\dots<i_k}x_{i_1}^{\alpha_1}\dots x_{i_k}^{\alpha_k}$$
where $\alpha$ runs over compositions of $n$ that do not have zero parts, viewed as finite sequences. 

It is easy to see that $M_\alpha$ is a basis of $\QSym_n$. However, we will use another basis that provides a connection from descents to characters and symmetric functions.
\begin{definition}\label{Definition-GesselFundamental}
For each subset $D \subseteq [n-1]$ define the {\em
	fundamental quasi-symmetric function}
$$\F_{n,D}(x) := \sum_{\substack{i_1\le i_2 \le \ldots \le i_n \\ {i_j < i_{j+1} \text{ if } j \in D}}} 
x_{i_1} x_{i_2} \cdots x_{i_n}.$$
\end{definition}

\begin{definition}\label{definition-AssociatedQuasisymmetrical}
	Let $\BBB$ be a set of combinatorial objects, equipped with a {\em descent map} $\Des: \BBB \to 2^{[n-1]}$ which associates to each
	element $b\in \BBB$ a subset $\Des(b) \subseteq [n-1]$. Define the
	{\em quasi-symmetric function associated with $\BBB$} by
	\[
	\Q(\BBB) := \sum\limits_{b\in \BBB} 
	\F_{n,\Des(b)}.
	\]
\end{definition}

\begin{definition}\label{SchurPositiveSet}
	We say that $\BBB$ is a {\em fine set} or {\em Schur-positive } if $\Q(\BBB)$ is symmetric and the coefficients $a_\lambda$ are non-negative integers in the expansion 
	$$\Q(\BBB) = \sum\limits_{\lambda \vdash n} a_\lambda s_\lambda.$$ 
	where $s_\lambda$ are Schur functions.
\end{definition}

The following fundamental theorem is due to Gessel.

\begin{theorem}{\rm \cite[Theorem 7.19.7]{EC2}}\label{theorem-GesselSYTSchur} 
	For every shape $\lambda \vdash n$,
	\[
	\Q({\SYT(\lambda)})=s_{\lambda}.
	\]
\end{theorem}

\begin{remark}\label{Remark:Skew-Schur}
	Skew Schur functions are defined similarly to Schur functions by replacing shape with skew shape. The result of Theorem~\ref{theorem-GesselSYTSchur} also holds for skew shapes. Namely, we have:
	$$\Q(\SYT(\lambda / \mu)) = s_{\lambda / \mu}.$$
\end{remark}

Every subset of $J=\NN$ can be associated with the sequence $\left\{j_n\right\}_{n\in\NN}$ where $$j_{n}=\begin{cases}
	1 & n\in J\\
	0 & n\notin J
\end{cases}.$$
\begin{definition}
	For a subset $J\subseteq\NN$, define $X^J=x^{\{j_n\}_{n\in\NN}}$ where $\{j_n\}_{n\in\NN}$ is the sequence 
	associated with $J$. 
\end{definition}

\begin{theorem}\label{theorem-Sp-combin}\cite{AdinHegedusRoichman}
	For every set $\A$, equipped with a descent map $\Des: \A \to 2^{[n-1]}$, the following statements are equivalent:
	$$\Q(\A)=\sum_{\lambda \vdash n}c_\lambda s_\lambda$$ 
	and
	\begin{equation*}\label{eq:equidist}
		\sum_{a \in \A } x^{\Des(a)}=\sum_{\lambda\vdash n}c_\lambda\sum_{T\in SYT(\lambda)}x^{\Des(T)}. 
	\end{equation*}
	In particular, if $\A$ is Schur-positive, then $c_\lambda \geq 0$. 
\end{theorem}
\begin{remark}
	The theorem in the paper is stated for subsets of $\symm{n}$. However, the proof remains correct if we replace $\A\subseteq\symm{n}$ with \enquote{$\A$ equipped with a descent map}.
\end{remark}
The Frobenius characteristic map can be redefined in terms of quasisymmetric functions as follows:
\[
\ch(\chi^\lambda) = \sum_{T \in SYT(\lambda)} \F_{n,\Des(T)}.
\]

We often rely on the following observation. Hence, even though it is evident, we state it explicitly.
\begin{lemma}\label{lemma:QSymDisjointUnion}
	Let $A$ and $B$ disjoint sets equipped with a descent map 
	$$\Des:A\uplus B \to 2^{n-1}$$
	for some $n\in\NN$. Then
	$$\Q(A\uplus B) = \Q(A)+\Q(B).$$
\end{lemma}
 
\subsection{Examples of Schur-positive sets}\label{subsection: SchurPositiveExamples}
Recall from the Subsection \ref{subsection: SymmetricFunctions} that a symmetric function is {\em Schur-positive} if all the coefficients in its expansion in the basis of Schur functions are non-negative. Determining whether a given symmetric function is Schur-positive is a significant problem in contemporary algebraic combinatorics~\cite{Stanley_problems}.

The following problem is also long-standing; see, for example, \cite {Gessel-Reutenauer}. 

\begin{problem}\label{prob:symmetric}
	Given a set $A$, equipped with a descent map, 
	is $\Q(A)$ symmetric?
	In case of an affirmative answer, is it Schur-positive?
\end{problem}

In the remainder of the section, we describe a few classic examples of Schur-positive sets in $\symm{n}$.
These sets and their generalizations are studied in more depth in Chapters~\ref{chapter: CyclicDescentInvDesClass} and~\ref{chapter: CyclicDescentUnityRoot} in perspective of cyclic descent extensions, which are defined in Section~\ref{section: CyclicDescentExtensions}.

Before proceeding to the mentioned above examples, we note that 
Schur-positivity questions are also crucial for sets that are not subsets of $\symm{n}$, as it implies the existence of a representation of $\symm{n}$. We have already seen an example in the previous subsection, the set of (skew) Young tableaux. Another set - of maximal chains in the non-crossing partition lattice and its subset of convex caterpillars will be studied in Chapter~\ref{chapter: ConvexCaterpillars}.

\begin{definition}
	The permutations $\pi,\sigma\in\symm{n}$ are said to be {\em Knuth equivalent} if $P_\sigma=P_\sigma$ where $P_\sigma$ is the first SYT in the RSK bijection.
	A {\em Knuth class} of $T$ is the set 
	$$K_T=\left\{\pi\in\symm{n}:P_\pi=T\right\}.$$
\end{definition}

Schur-positivity of a Knuth class is an immediate result of theorems \ref{theorem-GesselSYTSchur} and \ref{theorem-RSKProperties}. To see that, observe that if $\RSK{\pi}{P_\pi}{Q_\pi}$ then $\Des(\pi)=\Des(Q_\pi)$ and that by theorem \ref{theorem-RSKProperties} we have bijection from  $K_T$ to $SYT_{\lambda}$ where $\lambda = \shape{T}$ and therefore the descents of $K_T$ and $SYT(\lambda)$ are equally distributed. Combining with theorem~\ref{theorem-GesselSYTSchur} we obtain
$$\Q(K_T)=\Q(SYT(\lambda)=s_\lambda.$$

Inverse descent classes defined below are another classic example.
\begin{definition}\label{definition-InvDesClass}
	Let $J\subseteq [n-1]$. The {\em inverse descent class} of $\symm{n}$ corresponding to $J$ is defined by 
	$$D_{J,n}^{-1} = \left\{\pi\in\symm{n}:\Des(\pi^{-1})=J\right\}.$$
 \end{definition}

Schur-positivity of inverse descent classes follows from theorem \ref{theorem-RSKProperties} and the following equality:
\begin{equation*}
\begin{aligned}
	D_{J,n}^{-1} &= \left\{\pi\in\symm{n}:\Des(\pi^{-1})=J\right\}  \\
	 &= \biguplus_{\Des{T}=J}\left\{\pi:\RSK{\pi}{P_\pi}{T}\right\} \\
	 &= \biguplus_{\Des{T}=J}K_T.
\end{aligned}
\end{equation*}

Combining with lemma~\ref{lemma:QSymDisjointUnion} we obtain
$$\Q(D_{J,n}^{-1})=\Q(\biguplus_{\Des{T}=J}K_T)=\sum_{\Des{T}=J}s_{\shape{T}}.$$

A classic result by Foata and Schützenberger proved in \cite{FoataSchutz} (also see Section 1.4 in \cite{EC1} for more accessible proof) states that there exists a bijection $\phi:\symm{n}\to\symm{n}$ such that 
\begin{enumerate}
	\item $\maj(\pi)= \inv(\phi(\pi))$.
	\item $\Des(\pi^{-1})=\Des(\phi(\pi)^{-1})$.
\end{enumerate}

It immediately followsfrom Defiition~\ref{definition:CoxeterLength} that Coxeter length of $\pi$ and $\pi^{-1}$ are equal. In addition, by  Theorem~\ref{thm:CoxeterLengthInversionNumber}, the Coxeter length of a permutation $\pi$ in $\symm{n}$ is equal to $\inv(\pi)$. Therefore, we $\inv(\pi)=\inv(\pi^{-1})$. Hence, the mapping $\phi$ satisfies the following equations: 
\begin{enumerate}
	\item $\imaj(\pi)=\maj(\pi^{-1})=\inv(\phi(\pi^{-1}))=\inv(\phi(\pi^{-1})^{-1})$.
	\item $\Des(\pi)=\Des(\phi(\pi^{-1})^{-1})$	
\end{enumerate}
If we denote $\psi(\pi)=(\phi^{-1}(\pi^{-1}))^{-1}$ we obtain a bijective function from $\symm{n}$ to $\symm{n}$ such that
\begin{enumerate}
	\item $\inv(\pi)=\imaj(\psi(\pi))$.
	\item $\Des(\pi)=\Des(\psi(\pi))$.
\end{enumerate}
Hence the distributions of $\Des$ on the sets 
$$A=\left\{\pi\in\symm{n}:\inv(\pi)=m\right\}$$ and $$B=\left\{\pi\in\symm{n}:\imaj(\pi)=m\right\}$$ are equal. By Theorem~\ref{theorem-Sp-combin}, the Schur-positivity of a set $X$ depends only on the distribution of $\Des$ on $X$. Therefore, $A$ is Schur-positive only and only if $B$ is. It is obvious that $B$, the set of permutations $\pi$ in $\symm{n}$ with $\\imaj(\pi)=m$, can be expressed as a disjoint union of inverse descent classes and is Schur-positive. Therefore, the set of permutations with fixed length $m$ is also Schur-positive.

The Schur-positive sets mentioned above will be explored in more depth in Chapter~\ref{chapter: CyclicDescentInvDesClass}. 

Another important class is the set of permutations with a fixed cycle type. In \cite{Gessel-Reutenauer}, it was shown that they are Schur-positive. Their generalization - the roots of unity are studied in Chapter~\ref{chapter: CyclicDescentUnityRoot}. 

\chapter{Schur-positivity of convex caterpillars}
\label{chapter: ConvexCaterpillars}
	\section{Introduction}

A symmetric function is called {\em Schur-positive} if all the
coefficients in its expansion in terms of Schur functions are non-negative.
 	
With a set $A$ of combinatorial objects, equipped with a {\em descent map} $\Des: A \to 2^{[n-1]}$, one associates the
quasi-symmetric function
\[
\Q(A) := \sum\limits_{\pi\in A} \F_{n,\Des(\pi)}
\]
where $\F_{n,D}$ (for $D \subseteq [n-1]$) are Gessel's {\em fundamental quasi-symmetric functions;} see Subsection~\ref{section:prel_quasi} for more details.

\medskip

This chapter focuses on Schur-positive sets of maximal chains in partially ordered sets. Maximal chains in a labeled poset $P$ are equipped with a natural descent map.
A well-known conjecture of Stanley~\cite[III, Ch. 21]{Stanley_thesis} 
implies that all examples of Schur-positive labeled posets, in this sense,
correspond to intervals in the Young lattice.


Another way to equip the set of maximal chains with a descent map is by the edges in the Hasse diagram.
A classic example of a Schur-positive set 
of this type, the set of all maximal chains in the non-crossing partition lattice of type $A$,
was given by Stanley~\cite{StanleyPFnNCPL}. 
An $EL$ edge-labeling of this poset was presented in an earlier work of Björner~\cite{Bjorner}; see also~\cite{BW, Mcnamara, AR_OMCNCPL}.

\medskip

The goal of this chapter is to present an interesting set of maximal chains in the non-crossing partition lattice $NC_n$ 
(equivalently: a set of edge-labeled trees) 
which is Schur-positive. 
We will use a variant of Björner's $EL$-labeling, presented in~\cite{AR_OMCNCPL}.

\medskip

It is well known that 
maximal chains in the non-crossing partition lattice may be interpreted as factorizations of the $n$-cycle $(1,2,\dots,n)$
into a product of $n-1$ transpositions.

\begin{definition}
	A factorization $t_1\cdots t_{n-1}$ of the $n$-cycle $(1,2,\dots,n)$ as a product of transpositions is called {\em linearly ordered} if, 
	for every $1 \le i \le n-2$, $t_i$ and $t_{i+1}$ have a common letter.
\end{definition}

\begin{example}
	The factorization $(1\ 4)(1\ 3)(1\ 2)$ of $(1\ 2\ 3\ 4)$ is linear, while $(1\ 3)(3\ 4)(1\ 2)$ is not. 
\end{example}

This definition is motivated by Theorem~\ref{theorem-linear} below.
Denote the set of linearly ordered factorizations of  $(1,2,\dots,n)$ by $\mathcal{L}_n$.

\begin{proposition}\label{proposition:LnSize}
	For every $n\ge 1$, the number of linearly ordered factorizations of the $n$-cycle $(1,2,\dots,n)$ is 
	\[
	|\mathcal{L}_n|=n2^{n-3}.
	\]	\end{proposition}
Our main result is

\begin{theorem}\label{theorem-main}
	The set of linearly ordered factorizations of the $n$-cycle $(1,2,\dots,n)$

	satisfies 
	$$\Q(\mathcal{L}_n)=\sum_{k=0}^{n-1} (k+1)s_{(n-k,1^k)},$$
	where the descent set of any $u\in \mathcal{L}_n$ is defined by the edge labeling of~\cite{AR_OMCNCPL}.
	In particular, $\mathcal{L}_n$ is Schur-positive.
\end{theorem}


It should be noted that Theorem~\ref{theorem-main} does not follow from Stanley's proof of the Schur-positivity of the set of all maximal chains in $NC_n$; Stanley's action on maximal chains does not preserve linearly ordered chains. 

We prove Theorem~\ref{theorem-main} by 
translating it into the language of geometric trees called convex caterpillars.

\begin{definition}
	A tree is called a {\em caterpillar} if the subgraph obtained by removing all its leaves is a path. This path is called the \emph{spine} of the caterpillar.
\end{definition}

\begin{definition}
	A {\em convex caterpillar} of order $n$ is a caterpillar drawn in the plane such that
	\begin{enumerate}
		\item[(a)]
		the vertices are in convex position (say, the vertices of a regular polygon) and are labeled ${1,\ldots,n}$ clockwise;
		\item[(b)]
		the edges are drawn as non-crossing straight line segments; and
		\item[(c)]
		the spine forms a cyclic interval $(a,a+1),(a+1,a+2),\ldots,(b-1,b)$ in $[n]$.
	\end{enumerate}
	
	Denote by $\Ct_n$ the set of convex caterpillars of order $n$.
\end{definition}

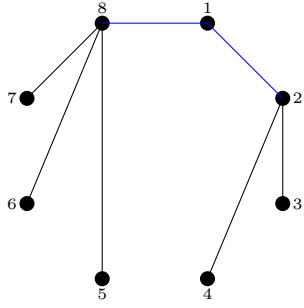
\begin{figure}[hbt]
	\begin{center}
		\begin{tikzpicture}[scale=1]
			\fill (2.4,3.4) circle (0.1) node[above]{\tiny 1};
			\fill (3.4,2.4)	circle (0.1) node[right]{\tiny 2};
			\fill (3.4,1) circle (0.1) node[right]{\tiny 3}; 
			\fill (2.4,0) circle (0.1) node[below]{\tiny 4};
			\fill (1,0) circle (0.1) node[below]{\tiny 5};
			\fill (0,1)	circle (0.1) node[left]{\tiny 6};
			\fill (0,2.4) circle (0.1) node[left]{\tiny 7};
			\fill (1,3.4) circle (0.1) node[above]{\tiny 8};
			
			\draw (3.4,2.4)--(3.4,1);
			\draw (0,1)--(1,3.4)--(1,0);
			\draw (3.4,2.4)--(2.4,0); 
			\draw (0,2.4)--(1,3.4);
			
			\draw[blue] (1,3.4)--(2.4,3.4); 
			\draw[blue] (3.4,2.4)--(2.4,3.4);
		\end{tikzpicture}
	\end{center}
	\caption{A convex caterpillar and its spine}\label{fig:1}
\end{figure}

\begin{example}
	Figure~\ref{fig:1} shows a convex caterpillar $c \in \Ct_8$, with spine consisting of the edges $(8,1)$ and $(1,2)$, forming a cyclic interval.
\end{example}

Goulden and Yong~\cite{Goulden-Yong} introduced a mapping from factorizations of $(1,2,\ldots,n)$ to non-crossing geometric trees.
This mapping is not injective: to recover the factorization from the tree, one has to choose a linear extension of a specific partial order on the edges, which we call the {\em Goulden-Yong partial order}; see Definition~\ref{definition-GLOrder} below.

In a previous work~\cite{YK-CenterPaper}, we proved that the Goulden-Yong order is linear if and only if the geometric tree is a convex caterpillar; see Theorem~\ref{theorem-linear} below.
It follows that the Goulden-Yong map, restricted to the set $\mathcal{L}_n$ of linearly ordered factorizations, is a bijection between $\mathcal {L}_n$ and the set $\Ct_n$ of convex caterpillars of order $n$.

\begin{definition}\label{definition-Des-c}
	The {\em descent set} of a linearly ordered factorization $u = (t_1, \ldots, t_{n-1}) \in \mathcal{L}_n$ is
	\[
	\Des(u) := \{i\in [n-2] \,:\, t_i=(b,c) \text{ and } t_{i+1}=(b,a) \text{ with }\ c > a\}.
	\]
\end{definition}

\begin{example} 
	The convex caterpillar $c \in \Ct_8$, drawn in Figure~\ref{fig:1},
	corresponds to the linearly ordered word
	\[
	u = \left( (7,8),(6,8),(5,8),(1,8),(1,2),(2,4),(2,3) \right) \in U_8,
	\]
	for which $\Des(u) = \{1,2,3,4,6\}$.
\end{example}

\begin{remark} In \cite{AR_OMCNCPL}, the authors define a map $\phi$ from the set denoted here $\mathcal{L}_n$ to the symmetric group $\symm{n-1}$; for a detailed description see Subsection~\ref{section:phi} below. 
The map $\phi$ is an $EL$-labeling of the non-crossing partition lattice. This property, relations to Björner's $EL$-labeling 
and other positivity phenomena will be discussed in another paper.
\end{remark}

It turns out that our Definition~\ref{definition-Des-c} above fits nicely with this map.

\begin{lemma}
	For any $u \in \mathcal{L}_n$, 
	\[
	\Des(\phi(u)) = \Des(u).
	\]
\end{lemma}

See Proposition~\ref{proposition:phi_caterpillar_des} below.
We further show that the number of caterpillars with a given descent set depends only on the cardinality of the descent set.	

\begin{lemma}
	For every subset $J\subseteq [n-2]$,
	\[
	|\{c \in \Ct_n:\ \Des(c)=J\}|= |J| + 1.
	\]
\end{lemma}

These two key lemmas are used to prove Theorem~\ref{theorem-main}.


\section{Background}
This section provides the necessary definitions and historical background to explain the main results. More information can be found in the references.

\subsection{Maximal chains in the non-crossing partition lattice}

The systematic study of non-crossing partitions began with Kreweras \cite{Kreweras-NonCrossing} and Poupard \cite{Poupard-NonCrossing}.
Surveys of results and connections with various areas of mathematics can be found in 
\cite{Simion-NonCrossing} and \cite{Armstrong}.

A {\em non-crossing partition} of the set $[n]$ is a partition $\pi$ of $[n]$ into nonempty blocks with the following property: 
for every $a<b<c<d$ in $[n]$, if some block $B$ of $\pi$ contains $a$ and $c$ and some block $B^\prime$ of $\pi$ contains $b$ and $d$, then $B=B^\prime$. 
Let $NC_n$ be the set of all non-crossing partitions of $[n]$.
Define a partial order on $NC_n$, by refinement: 
$\pi \leq \sigma$ if every block of $\pi$ is contained in a block of $\sigma$.
This partial order turns $NC_n$ into a graded lattice.

An \emph{edge labeling} of a poset $P$ is a function from the set of edges of the Hasse diagram of $P$ to the set of integers. Several different edge labelings of $NC_n$ were defined and studied 
by Björner \cite{Bjorner}, Stanley \cite{StanleyPFnNCPL}, and Adin and Roichman \cite{AR_OMCNCPL}.
Let $\Lambda$ be an edge labeling of $NC_{n+1}$,
and let $F_{n+1}$ be the set of maximal chains in $NC_{n+1}$.
For each maximal chain $\mathfrak{m} : \pi_0 < \pi_1 < \dots < \pi_n$ define
\[
\Lambda^*(\mathfrak{m}) := \left( \Lambda(\pi_0, \pi_1), \ldots, \Lambda(\pi_{n-1}, \pi_n) \right) \in \NN^n,
\]
with a corresponding {\em descent set}
\[
\Des(\mathfrak{m}) := \left\{ i \in [n-1] \,:\, \Lambda(\pi_{i-1}, \pi_i) >
\Lambda(\pi_i, \pi_{i+1}) \right\}.
\]

\medskip

The non-crossing partition lattice is intimately related to cycle factorizations.
The $n$-cycle $(1,2,\ldots,n)$ can be written as a product of $n-1$ transpositions. There is a well-known bijection between such factorizations and the maximal chains in $NC_{n+1}$; see, for example, \cite[Lemma 4.3]{NC_SurprisingLocations}. 
A classical result of Hurwitz states that the number of such factorizations is $n^{n-2}$ \cite{HurwitzClassical, HurwitzReproof},
thus equal to the number of labeled trees of order $n$.
In the next section, we will describe a connection between maximal chains and 	
geometric trees.



\subsection{The Goulden-Yong partial order}

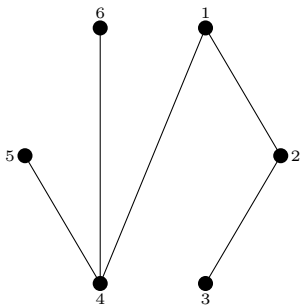
\begin{figure}[hbt]
	\begin{center}
		\begin{tikzpicture}[scale=1]
			\fill (2.4,3.4) circle (0.1) node[above]{\tiny 1};
			\fill (3.4,1.7)	circle (0.1) node[right]{\tiny 2};
			\fill (2.4,0) circle (0.1) node[below]{\tiny 3};
			\fill (1,0) circle (0.1) node[below]{\tiny 4};
			\fill (0,1.7)	circle (0.1) node[left]{\tiny 5};
			\fill (1,3.4) circle (0.1) node[above]{\tiny 6};
			
			\draw (2.4, 3.4) -- (3.4, 1.7);
			\draw (3.4, 1.7) -- (2.4, 0);
			\draw (2.4, 3.4) -- (1, 0);
			\draw (1, 0) -- (0, 1.7);
			\draw (1, 0) -- (1, 3.4);
		\end{tikzpicture}
	\end{center}
	\caption{$G(w)$ for $w = ((1,4),(4,6),(4,5),(1,2),(2,3))$}\label{fig:2}
\end{figure}

With each sequence of $n-1$ different transpositions $w = (t_1,\dots,t_{n-1})$, associate a geometric graph $G(w)$ as follows.
The vertex set is the  set of vertices of a regular $n$-gon, labeled clockwise $1,2,\dots,n$.
The edges correspond to the given transpositions $t_1,\dots,t_{n-1}$, where the edge corresponding to a transposition $t_k = (i,j)$ is the line segment connecting vertices $i$ and $j$.
See Figure~\ref{fig:2} for the geometric graph $G(w)$ corresponding to $w = ((1,4),(4,6),(4,5),(1,2),(2,3))$.

Let $F_n$ be the set of all factorizations of the $n$-cycle $(1,2,\ldots,n)$ into a product of $n-1$ transpositions.
Write each element of $F_n$ as a sequence $(t_1, \ldots,t_{n-1})$, where $t_1 \cdots t_{n-1} = (1,2,\ldots,n)$.
The following theorem of Goulden and Yong gives the necessary and sufficient conditions for a sequence of $n-1$ transpositions to belong to $F_n$.

\begin{theorem}\cite[Theorem 2.2]{Goulden-Yong}\label{theorem-GY}
	A sequence of transpositions $w = (t_1, \ldots, t_{n-1})$ belongs to $F_n$ if and only if the following three conditions hold:
	\begin{enumerate}
		\item $G(w)$ is a tree.
		\item $G(w)$ is non-crossing; namely: two edges may intersect only in a common vertex.
		\item Cyclically decreasing neighbors: For every 
		$1 \leq i < j \leq n-1$, if $t_i = (a, c)$ and $t_j = (a, b)$, then $c >_a b$. Here, $<_a$ is the linear order $a <_a a + 1 <_a \dots <_a a-1$.
	\end{enumerate}
\end{theorem}

For example, the graph in Figure~\ref{fig:2} corresponds to a sequence $w \in F_6$ and therefore satisfies the conditions of Theorem~\ref{theorem-GY}.

Note that a sequence $w = (t_1, \ldots, t_{n-1}) \in F_n$ carries more information than its Goulden-Yong tree $G(w)$: It defines a {\em linear order} on the edges, with the edge corresponding to $t_i$ preceding the edge corresponding to $t_j$ whenever $i < j$. 
How much of this information can be recovered from the tree?

\begin{defn}\label{definition-GLOrder}
	Let $T$ be a non-crossing geometric tree 
	(namely, satisfying conditions 1 and 2 of Theorem~\ref{theorem-GY}) 
	on the set of vertices of a regular $n$-gon, labeled clockwise $1,2,\dots,n$.
	Define relations $\prec_T$ and $\le_T$ on the set of edges of $T$ as follows:
	\begin{enumerate}
		\item For edges $s$ and $t$ of $T$, we have $s \prec_T t$, if $s=(x,z)$ and $t=(x,y)$ have a common vertex $x$ and $z >_x y$ as in condition 3 of Theorem~\ref{theorem-GY}. 
		\item We have $(a,b) \le_T (c,d)$ if there exists a sequence of edges $(a,b) = t_0, \ldots, t_k =(c,d)$ $(k \ge 0)$ such that for every $0 \leq i \leq k-1$, the edges $t_i $ and $t_{i+1}$ satisfy $t_i \prec_T t_{i+1}$.
		\end{enumerate}
\end{defn}

\begin{lemma}\label{lemma-GouldenYongOrderOnT} Let $T$ be a non-crossing geometric tree. Then $\le_T$ is a partial order on the set of edges of $T$.
\end{lemma}

\begin{proof}
	One needs to show that $\le_T$ is reflexive, transitive and antisymmetric relation on the edges of $T$. Observe that $\le_T$ is reflexive, since for every edge $t$ of $T$, the sequence $)$ trivially satisfies the conditions of $\le_T$ as in Definition~\ref{definition-GLOrder}, hence $t \le_T t$ holds. Next, observe that if $s \le_T t$ and $t \le_T u$, then by definition of $\le_T$ there exist sequences $t = t_0, t_1, \dots, t_k = s$, $s = s_0 < \dots s_l = u$ such that $t_i \prec_T t_{i+1}$ for every $0\leq i \leq k-1$ and $s_j \prec_T s_{j+1}$ for every $0 \leq j \leq l-1$.
	Then the sequence $t = t_0,\dots,t_k=s=s_0,s_1,\dots,s_l=u$
	 is the sequence that satisfies the definition of $\le_T$ for $t$ and $u$. 
	 
	Before showing that $\le_T$ is antisymmetric, note that $\prec_T$ is a strict order on any subset of edges of $T$ that share a common vertex. In addition, $\prec_T$ is an amtisymmetric relation as well. Let $e$ and $f$ edges in $T$ such that $e \le_T f$. 
	By definition of $\le_T$ there exists a sequence $e = t_0,\dots,t_k=f$ such that $t_i \prec_T t_{i+1}$ for every $0 \leq i \leq k-1$ such that $t_i \prec_T t_{i+1}$. Therefore, if we replace every subsequence  $t_1, \dots, t_k$ of form $(v, u_m),\dots,(v, u_{m+l})$ with $(v, u_m), (v, u_{m+l})$, we obtain a subsequence $s_1,\dots,s_{k^\prime}$ of $t_0,\dots,t_k$ with the same property $s_i \prec_T s_{i+1}$ for every $0 \leq i \leq k^\prime - 1$, every two consecutive edges share a common vertex, and for every vertex there are at most two edges  in the sequence $s_0,\dots,s_{k^\prime}$ that have it as an endpoint. If we denote $e = (a, u_1)$, $s_i = (u_i, u_{i+1})$ for $1 \leq i \leq k^\prime - 1$, $f= (u_{k^\prime-1}, b)$ the sequence of vertices $a, u_0, \dots, u_{k^\prime-1}, b$ is a path from $a$ to $b$ in the tree $T$.  One the other hand, if $f \le_T e$, by the same reasoning, one shows an existence of a sequence $f = s^\prime_1,\dots,s^\prime_{k^{\prime\prime}}$ with the same properties as $s_0,\dots,s_{k^\prime}$. Since there is exactly one path between two vertices in a tree, we have $k^\prime = k^{\prime\prime}$ and that $(s_0,\dots,s_{k^\prime})=(s^\prime_{k^\prime\prime},\dots,s^\prime_9).$ Hence, we have both $s_i \prec_T s_{i+1}$ and $s_{i+1} \prec_T s_i$ for every $0\leq i \leq k^\prime-1$, which is a contradiction, unless $k^\prime=0$. Therefore $\le_T$ prime is antisymmetric.
	 We have shown that $\le_T$ is reflexive, transitive and antisymmetric. Therefore, we are done. 
	\end{proof}

We call $\le_T$ the {\em Goulden-Yong partial order} corresponding to $T$.

\begin{observation}\label{conclusion:linear_extension}
	For every factorization $w = (t_1,\ldots,t_n) \in F_n$, the order

	$t_1 < t_2 < \ldots < t_n$ is a linear extension of the Goulden-Yong order $<_{G(w)}$.
\end{observation}

\begin{example}
	In Figure~\ref{fig:2}, the tree $T = G(w)$ gives the partial order satisfying $(1,4) <_T (4,6) <_T (4,5)$ and $(1,4) <_T (1,2) <_T (2,3)$.
	It is not a linear order. The order $(1,4) < (4,6) < (4,5) < (1, 2), (2,3)$
	is a linear extension of it.
\end{example}

\section{Convex caterpillars}

In this section, we prove Theorem \ref{theorem-main} using the properties of convex caterpillars.
\subsection{Basic properties of convex caterpillars}
Let us use the following conventions.
All arithmetical operations on the elements of $[n]$ will be performed modulo $n$. We write $[a, b]$ to denote the cyclic interval $\{a, a+1, \dots, b\}$. Using this notation, for the edges $(a\ b), (a\ c)$ of a geometric non-crossing tree $T$, we have $(a\ c) <_T (a\ b)$ if and only if $b \in [a, c]$.

The following result was proved in \cite{YK-CenterPaper}. Later in this section, we provide a somewhat different proof, details of which will be used later.

\begin{theorem}\cite[Theorem 3.2]{YK-CenterPaper}\label{theorem-linear}
	The Goulden-Yong order on the edge set of a non-crossing geometric tree $T$ is linear (total) if and only if $T$ is a convex caterpillar.
\end{theorem}

The following observation follows from the fact that a linear extension of a Goulden-Yong order $<_T$ on the edges of 
the geometric non-crossing tree $T$ corresponds to a factorization of the cycle 
$(1\dots n)$ into $n-1$ transpositions.
\begin{observation}\label{observation:consectuve_edges}
	If $T$ is a geometric non-crossing tree and $<_T$ is linear, then any two consecutive edges in $<_T$, viewed as transpositions in $\symm{n}$, do not commute and therefore have a common vertex.
\end{observation}

The following lemma gives sufficient conditions for $<_T$ to be non-linear.
\begin{lemma}\label{lemma:non_linear}
	Let $T$ be a non-crossing geometric tree. The order $<_T$ is not linear in the following cases.
	\begin{enumerate}
		\item  There are edges $(a\ b), (c\ d), (e\ f)$ of $T$ such that $(a\ b) <_T (c\ d), (e\ f)$ and 
		$c,d \in [a,b-1]$ and $e, f \in [b, a-1]$.
		\item $T$ has edges $(a\ b), (c\ d), (e\ f)$ such that
		$(c\ d), (e\ f) < (a\ b)$ and $c,d \in [b+1, a]$ and $e, f \in [a+1, b]$.
	\end{enumerate}
\end{lemma}
\begin{proof}
	We prove the first case, the second one being similar by reversing directions. 
	Assume that $<_T$ is linear and that the first case holds. For every $v \in [a+1, b-1]$, the edge $(v\ b)$ is smaller than $(a\ b)$ in $<_T$. The same is true for any edge $(a\ v)$ with $v \in [b+1, a-1]$.
	Combining this fact with the non-crossing property of $T$, we deduce that the endpoints of any edge that is larger than $(a\ b)$ in $<_T$ are either both in $[a, b-1]$ or both in $[b, a-1]$. On the other hand, by Observation~\ref{observation:consectuve_edges}, every pair of consecutive edges shares a common vertex. Therefore, the endpoints of any two consecutive edges greater than $(a\ b)$ in $<_T$ must all be in either $[a, b-1]$ or $[b, a-1]$. Using induction, we conclude that the endpoints of the edges greater than $(a\ b)$ in $<_T$ are all in $[a, b-1]$ or all in $[b, a-1]$, contradicting the assumption.
\end{proof}

We are ready to prove Theorem \ref{theorem-linear}.

\begin{proof}[Proof of Theorem \ref{theorem-linear}]
	Assume that $T$ is a convex caterpillar. If the spine of $T$ does not contain any edges, then $T$ is a star; hence, every pair of edges is comparable because they have a common vertex, and the Goulden-Yong order is linear.
	
	Otherwise, let $(a\ a+1), (a+1\ a+2)\dots,(b-2\ b-1), (b-1\ b)$ be the spine of $T$. For every pair of edges $(k\ l)$ and $(k\ m)$ that share a common vertex $k$, $(k\ l) <_T (k\ m)$ holds if and only if $m \in [k+1, l-1]$ where $[k+1, l-1]$ denotes the cyclic interval $\{k+1,\dots,l-1 \}$. Note that it is simply a restatement of the fact that the neighbors of $k$ are ordered counterclockwise.
	Therefore, all the edges of the spine are comparable and $(a\ a+1) <_T (a+1\ a+2)\dots <_T (b-1\ b)$. It also implies that if $(k\ k+1)$ is an edge in the spine, then for every edge $(k\ l)$ that has $k$ as an endpoint, $(k\ l) \leq_T (k\ k+1)$ and for every edge $(k+1\ m)$ that has $k+1$ as an endpoint, we have $(k\ k+1) \leq_T (k+1\ m)$. It follows that if the edges $(k\ k+1),(k+1\ k+2)\dots,(m-1\ m)$ are in the spine of $T$, then for every edge $(k\ j)$ connected to $k$ and every edge $(m\ l)$ connected to $m$ we have $(k\ j) \leq_T (k\ k+1) <_T \dots <_T(m-1\ m) \dots \leq_T (m\ l)$. Any edge in the caterpillar is either in the spine or has a common vertex with an edge in the spine. Therefore, any two edges are comparable.
	
	To prove the converse statement, assume that $<_T$ is linear. Enumerate the edges such that $(a_1\ b_1)<_T\dots<_T(a_{n-1}\ b_{n-1})$. Let $t_l$ denote the $l$-th edge in $<_T$, namely $(a_l\ b_l)$.
	
	By Observation~\ref{observation:consectuve_edges}, since $<_T$ is linear, every pair of adjacent edges shares a common vertex. 
	Now we claim that the first edge must be of the form $(i\ i+1)$ for some $i$. To see that, denote $t_1 = (i\ j)$. By the proof of Lemma~\ref{lemma:non_linear} the endpoints of every other edge of $T$, which is necessarily larger that $t_1$ in the linear order, must be either both in $[i,j-1]$ or both in $[j,i-1]$. By Lemma~\ref{lemma:non_linear}, either there are no edges with endpoints in $[i, j-1]$ or no edges with endpoints in $[j, i-1]$, which implies that either $[i,j]$ or $[j, i]$ consists of precisely two points, which means that either $i=j+1$ or $j=i+1$ holds.

	Next, observe that $i$ must be a leaf since the existence of an edge $(i\ j) <_T (i\ i+1)$, contradicts the fact that $(i\ i+1)$ is the first edge in $<_T$.
	Next, we use induction to show that for every $1\leq m \leq n-1$, the following statements hold:
	\begin{enumerate}
		\item The endpoints of the first $m$ edges in $<_T$ form a cyclic interval $[j, k]$ in $[n]$.
		\item Vertices $j, j+1, \dots, i-1, i$ are leaves in $T$.
		\item Edges $(i\ i+1), (i+1\ i+2), \dots (k-1\ k)$ are edges in $T$ and occur among the first $m$ edges.
		\item Every edge that has at least one of $j, j+1 \dots, k-1$ as an endpoint occurs among the first $m$ edges.
		\item Either $t_m$ = $(k-1\ k)$ or $t_m = (k\ j)$ is true. 
	\end{enumerate}
	
	As we have shown, $t_1 = (i\ i+1)$ for some $1\leq i \leq n$ and $i$ is a leaf. Therefore, the statements hold for $m=1$. Assume that the statement holds for $m$. According to the induction hypothesis, we have either $t_m = (k\ j)$ or $t_m = (k-1\ k)$. From the linearity of $<_T$, $t_{m+1}$ and $t_m$ have a common vertex. Again, by the induction hypothesis $j$ and $k-1$ cannot be endpoints of $t_{m+1}$, therefore, $t_{m+1}$ must have $k$ as an endpoint.
	Next, we show that $t_{m+1}$ is either $(k\ k+1)$ or $(k\ j-1)$. Assume $t_{m+1} = (k\ l)$ for $l\neq j-1, k+1$. There must be an edge (necessarily not among the first $m$ edges) with endpoint $j-1$, and similarly an edge with endpoint $k+1$. Therefore we must have edges $(k\ l) <_T (s\ t), (u\ v)$ with $u,v \in [k\ l-1]$ and $s, t \in [l, k-1]$ contradicting Lemma~\ref{lemma:non_linear}.
	Now if $t_{m+1}=(k\ k+1)$, we are done, since statements $1$ and $2$ hold by the induction hypothesis for $m$, $3$ and $5$ hold for $m+1$, and $4$ holds because $(k\ k+1)$ must be the greatest edge in $<_T$ that has $k$ as an endpoint according to Part 3 of Theorem~\ref{theorem-GY}.
	
	If $t_{m+1} = (k\ j-1)$, then for every $v \in [k+1, j-2]$ we have $(j-1\ v) <_T (j-1\ k)$, which implies that if $(j-1\ v) = t_l$ then $l\leq m$, which is impossible since $v\notin [j,k]$ according to the induction assumption.
	On the other hand, for every $v \in [j, k-1]$, $(j-1\ v)$ cannot be an edge since, by assumption, every edge with endpoints $j, \dots k-1$ occurs among the first $m$ edges. Therefore, $j-1$ must be a leaf.
	Again, it is easy to check that assumptions $1, 2, 3, 4, 5$ still hold for $m+1$.
	Now, for $m = n-1$, we see that $T$ must be a geometric caterpillar because, by construction, the vertices that are not leaves are $i+1, i+2, \dots, k$ for some $k$, with edges $(i+1\ i+2), \dots, (k-1, k)$ connecting them.
\end{proof}

For example, the tree in Figure~\ref{fig:2} is a caterpillar but not a convex one. The corresponding Goulden-Yong order is not linear.

\begin{corollary}\label{cor:linear}
	A non-crossing geometric tree $T$ on $n$ vertices is a convex caterpillar if, and only if, there is a unique $w \in F_n$ such that $G(w) = T$. 
\end{corollary}

We can view each edge of $T$ as a transposition. Therefore, in the rest of the paper, we shall identify a convex caterpillar $c \in \Ct_n$ with the corresponding sequence of transpositions $(t_1,\ldots,t_{n-1}) \in F_n$.

\begin{proposition}\label{proposition:common_vertex}
	In a convex caterpillar $c=(t_1,\dots ,t_{n-1})$:
	\begin{enumerate}
		\item 
		Any two consecutive edges $t_i$ and $t_{i+1}$ share a common vertex.
		\item 
		The first edge $t_1$ is of the form $(a, a+1)$ for some $a$.
		The same holds for the last edge $t_{n-1}$.
	\end{enumerate}
\end{proposition}
\begin{proof}
	The first part is simply a restatement of Observation~\ref{observation:consectuve_edges}. The second part follows from the proof of Theorem~\ref{theorem-linear}.
\end{proof}

\begin{definition} Let $e$ be an edge of a caterpillar $c$.
	\begin{enumerate}
		\item 
		We say that $e$ is a \emph{branch} if (at least) one of its endpoints is a leaf.
		\item 
		We say that $e$ is a \emph{link} if its endpoints have cyclically consecutive labels.
	\end{enumerate}
\end{definition}

We get the following observation by carefully reading the proof of Theorem~\ref{theorem-linear}.
\begin{observation}\label{observation:first_last_edge}
	An edge of a convex caterpillar $c$ is both a link and a branch if and only if it is either the first or the last edge of $c$.
\end{observation}

\begin{lemma}\label{CaterpillarStructure}
	Let $c = (t_1,\dots, t_{n-1}) \in \Ct_n$. The following statements hold. 
	\begin{enumerate}
		\item 
		The endpoints of the first $k$ edges form a cyclic interval in $[n],$ for every $1\leq k \leq n-1$.
		\item 
		If the first edge is $(i,i+1)$, then the endpoints of the first $k$ branches which are leaves are $i, i-1, \dots, i-k+1$, in that order.
		\item 
		If the first edge is $(i, i+1)$ then the first $k$ links are $(i,i+1),(i+1,i+2),\ \dots, \ (i+k-1, i+k)$.
		\item For $1\leq k < n-1$,
		the product of the first $k$ edges, viewed as transpositions, is equal to the cycle $(l\ l+1\ \dots\ m)$ where $l$ is the leaf endpoint of the last branch among the first $k$ edges and $(m-1, m)$ is the last link among the first $k$ edges.
	\end{enumerate}
\end{lemma}
\begin{proof}
	Parts $1$, $2$ and $3$ follow from the proof of Theorem~\ref{theorem-linear}. Part $4$ follows by induction and using the fact that if the product of the first $k$ edges (viewed as transpositions) is the cycle $(l\ l+1 \dots\ l+k)$, where the cyclic interval $[l, l+k]$ is formed by the endpoints of the first $k$ edges, then $t_{k+1}$, viewed as transposition, is either $(l+k\ l+k+1)$ or $(l+k\ l-1)$. Multiplying the cycle $(l\ l+1\dots\ l+k)$ by $t_{k+1}$, it is easily verified that the statement is true for $k+1$.
\end{proof}

\begin{corollary}\label{conc:caterpillarstructure2}
	Every $c = (t_1,\dots, t_{n-1}) \in \Ct_n$ is completely determined by its first edge $t_1$ and the set of indices $i$ for which $t_i$ is a branch, , which can be any subset of [n-1] containing $1$ and $n-1$.
\end{corollary}

This proves, in particular, Proposition~\ref{proposition:LnSize}.

\subsection{A labeling of maximal chains}\label{section:phi}
The following labeling of maximal chains in the non-crossing partition lattice was introduced by Adin and Roichman in \cite{AR_OMCNCPL} and is closely related to the $EL$-labeling introduced by Björner in \cite{Bjorner}.
This subsection describes this labeling, denoted by $\phi$.

\medskip

Recall, from Definition \ref{definition-Des-c}, the notion of a descent set of a convex caterpillar. 

Next, we show the connection to the descents defined by the map $\phi$ in \cite{AR_OMCNCPL}. First, let us describe $\phi$.
For $w=(t_1,\dots,t_{n-1}) \in F_n$ define the partial products $\sigma_j = t_j\dots t_{n-1}$ with $\sigma_n = id$. By definition, $\sigma_j = t_j\sigma_{j+1}$. For $1 \leq j \leq n-1$ define 
$$A_j = \left\{1 \leq i \ n-1 : \sigma_j(i) > \sigma_{j+1}(i) \right\}.$$

We get the following statement from the discussion preceding Definition 3.2 in \cite{AR_OMCNCPL}.
\begin{proposition}\label{proposition:phi_justification}
	The following statements hold.
	\begin{enumerate}
		\item For each $1 \leq j \leq n-1$, $|A_j| = 1$.
		\item The map $\pi_w$ defined by $$\pi_w(j) = i\ if\ A_j = \left\{ i \right\}$$ is a permutation in $\symm{n-1}$.
	\end{enumerate}
\end{proposition}
\begin{defn}\label{definition:PhiAR}\cite[Definition 3.2]{AR_OMCNCPL} Define $\phi:F_n\to\symm{n-1}$ by $$\phi(w) = \pi_w.$$
\end{defn}
Define for each $w\in F_n:$ $$\Des(w)=\Des(\phi(w)).$$

\subsection{Descents of convex caterpillars}

We now calculate the restriction of $\phi$ to $Ct_n$. Before we proceed, note that the arguments used to prove Theorem~\ref{theorem-linear} hold for the last edges. Therefore, we can formulate statements analogous to the ones used in the proof for the first edges of $c$ for the last edges of $c$ instead. We sum them up in the following lemma. 
\begin{lemma}\label{lemma:ReverseCaterpillarStructure}
Let $c=(t_1,\dots,t_{n-1})\in\Ct_n$. The following statements hold for every $1\leq m \leq n-1$:
\begin{enumerate}
    \item The last edge of $c$, $t_{n-1}$ is of the form $(i-1\ i)$ for some $i$ and $i$ is a leaf of $c$.
    \item The endpoints of the last $m$ edges in $T$ for a cyclic interval in $[j, k]$ in $[n]$.
    \item Vertices $k, k-1, \dots,i$ are leaves in $T$ and are endpoints of the last $m$ edges.
    \item Edges $(j, j+1),\dots,(i-1, i)$ are edges in $T$ and occur among the last $m$ edges.
    \item Every edge that has $j+1,\dots,k$ as an endpoint occurs among the last $m$ edges.
    \item Either $t_{n-m} = (k, k+1)$ or $t_{n-m} = (k, j)$ is true.
\end{enumerate}
\end{lemma}

Note that the product of the last $m$ edges in a caterpillar is just $\sigma_{n-m}$. Parts 5 and 6 of the preceding lemma can be reformulated in terms of $\sigma_j$ and $t_j$, as follows.
\begin{corollary}\label{corollary:suffixProduct}
Let $c = (t_1,\dots,t_n)\in\Ct_n$. Then, for every $1\leq j \leq n-1$, $\sigma_j$ is equal to the cycle $(l\ l+1\ \dots\ m)$ for some $l\in[n]$ and $t_j = (l, l+1)$ or $t_j = (l, m)$. 
\end{corollary}

The following proposition gives us a simple method to calculate the values of $\phi$ restricted to $\Ct_n$.
\begin{proposition}\label{proposition:alt_description_phi}
	Let $c\in Ct_n$, $1 \leq j \leq n-1$ and $\sigma_{j+1} = (l\ l+1\ \dots m)$. Then 
	$$\phi(c)(j) = 
	\begin{cases*}
		l & if $\sigma_{j+1}(l) \neq l$ and $l \neq n$ \\
		m & if $\sigma_{j+1}(l) = l$ and $l > m$ or $l=n$ \\
		m-1 & if $\sigma_{j+1}(l) = l$ and $l < m$ 
	\end{cases*}
	$$
\end{proposition}
\begin{proof}
	By Corollary~\ref{corollary:suffixProduct}, we have $\sigma_j=(l\ l+1\ \dots m)$, and either $t_j = (l, l+1)$ and $\sigma_{j+1}=(l+1\ \dots\ m)$ hold, or $t_j = (l, m)$ and $\sigma_{j+1} = (l\ l+1\ \dots \ m-1)$. By Proposition~\ref{proposition:phi_justification} there is exactly one $1\leq i \leq n-1$ such that $\sigma_j(i) > \sigma_{j+1}(i)$, which is equal to $\phi(c)(j)$ by Definition~\ref{definition:PhiAR}. We calculate it for each case.
	Before proceeding, note that if $t_j = (l, l+1)$, then $\sigma_{j+1}(m) = l+1$ and if $t_j=(l, m)$ then $\sigma_{j+1}(m) = m$. In both cases, we have $\sigma_{j+1}(m) \neq l = \sigma_{j+1}(m)$.
	
	First, observe that if $l=n$, then whether $t_j=(l, m)$ or $t_j=(l-1\ l)$ we have
	$$\sigma_j(m) = n > \sigma_{j+1}(n),$$ which implies $\phi(c)(j)=m$.
	
	If $l\neq n$ and $t_j=(l, l+1)$ then
	$$\sigma_{j+1}(l) = l < l+1 = \sigma_j(l),$$ 
	hence $\phi(c)(j) = l$.
	
	If $l\neq n$ and $t_j = (m, l)$ and $l > m$, then 
	$$\sigma_{j}(m) = l > m = \sigma_{j+1}(m),$$
	and $\\phi(c)(j) = m$.
	Otherwise, 
	$$\sigma_{j+1}(m-1) = l < m = \sigma_j(m),$$
	hence, $\phi(c)(j) = m-1$.
\end{proof}

\begin{proposition}\label{proposition:phi_caterpillar_des}
	The descent set of a convex caterpillar, defined in Definition~\ref{definition-Des-c}, 
	coincides with the descent set defined using the map $\phi$.
\end{proposition}
\begin{proof}
Let $c=(t_1,\dots,t_{n-1})\in\Ct_n$. We show that for every $1\leq i \leq n-2$, $j$ is a descent of $\phi(c)$ if and only $j$ is a descent of $c$ as in Definition~\ref{definition-Des-c}, namely $t_j = (k, l)$, $t_{j+1} = (l, m)$ for $1\leq k,l,m\leq n$, such that $k > m$. We do this by verifying the statement by calculating the $\phi(c)(j), \phi(c)(j+1)$ using Proposition~\ref{proposition:alt_description_phi} for every $1\leq j \leq n-2$ and every possible configuration of $t_j, t_{j-1}$.

First, we check that the statement is true in the special case $j=n-2$.
If $j = n-2$, then by Lemma~\ref{lemma:ReverseCaterpillarStructure}, $\sigma_{n-1} = (i,i+1)$ for some $1\leq i \leq n$ and two possibilities for $t_{n-2}$.
\begin{enumerate}
\item $t_{n-2} = (i-1, i)$. If $i-1, i \neq n$, then $t_{n-2} = (i-1, i), t_{n-1}=(i, i+1)$ and $i-1 < i < i+1$, thus $n-2$ is not a decent of $c$ and is not a descent of $\phi(c)$, since $\phi(c)(n-2) = i-1 < i = \phi(c)(n-1)$. 

If we have $i=n$, then $t_{n-2} = (n-1, n), t_{n-1} = (n, 1)$, which means that $n-2$ is a descent of $c$. But we also have $\phi(c)(n-1) = 1$ and $\phi(c)(n-2) = n-1$, which means that $n-1$ is also a descent of $\phi(c)$. 

On the other hand, if $i-1 = n$, we have $t_{n-2}=(n, 1)$ and $t_{n-1}=(1, 2)$, which means that $n-2$ is a descent of $c$ and is a descent of $\phi(c)$, since $\phi(c)(n-2) = 2 >1= \phi(c)(n-1)$. 

\item $t_{n-2} = (i, i+2)$. If none of $i, i+1$ equals $n$, then $t_{n-2} = (i, i+2), t_{n-1}=(i, i+1)$, where $i+2 > i+1$, which implies that $n-2$ is a descent of $c$ and also a descent of $\phi(c)$ since $\phi(n-2) = i+1 > i = \phi(n-1)$. 

If $i+1=n$, then $t_{n-2} = (n-1, 1), t_{n-1} = (n-1, n)$, which implies that $n-2$ is not a descent of $c$ and also not a descent of $\phi(c)$, since $\phi(c)(n-2) = 1 < n-1 = \phi(c)(n-1)$. 

If $i=n$, then $t_{n-2} = (n, 2), t_{n-1} = (n, 1)$, which implies that $n-2$ is a descent of $c$ and also a descent of $\phi(c)$, as we have $\phi(c)(n-2) = 2 > 1 = \phi(c)(n-1)$.
\end{enumerate}

We show that the statement is true for $1\leq j<n-2$. By Lemma~\ref{lemma:ReverseCaterpillarStructure}, we have $\sigma_{j+2} = (l, l+1\  \dots\ m)$ for some $1\leq l \leq n$ and four possibilities $t_j, t_{j+1}$.
\begin{enumerate}

\item We have $t_j = (l-2, l-1)$ and $t_{j+1}=(l-1, l)$. Then $j$ is not a descent of $c$, if $l-2 < l$, which implies $l-2 < l-1 < l$. Applying Proposition~\ref{proposition:alt_description_phi}, we get $\phi(j)(l-2)$ and $\phi(j)(j+1)=l-1$, which means $j$ is also not descent of $\phi(c)$.

If $j$ is a descent of $c$, then either $l-1=n$ or $l-2=n$. If we have$l-2=n$, then $t_j = (n,1), t_{j+1} = (1,2)$ and $\sigma_{j+2} = (2\ \dots\ m)$. From Proposition~\ref{proposition:alt_description_phi}, it follows that $\phi(c)(j) = m$ and $\phi(c)(j+1) = 1$, which means that $j$ is  also a descent of $\phi(c)$.
In case $l-1 = n$, we have $t_j = (n-1, n), t_j=(n,1)$ and $\sigma_{j+2}=(1,\dots,m)$. Applying Proposition~\ref{proposition:alt_description_phi}, we get $\phi(c)(j)=n-1$, $\phi(c)(j+1)=m$, which implies that $j$ is also a descent of $\phi(c)$.

\item We have $t_j = (l-1, l)$ and $t_{j+1} = (l, m+1)$. In this case, $j$ is not a descent of $c$ if $l-1 < m+1$ holds, which implies that $l-1 < l < m < m+1$. Then we have $\phi(c)(j) = l-1$ and $\phi(c)(j+1) = m$, again indicating that $j$ is not a descent of $\phi(c)$. 

If $j$ is a descent of $c$, we have $l-1 > m+1$, together with the following possibilities:
\begin{enumerate}
\item We have $l-1 = n$, $t_j=(n,1)$, $t_{j+1} = (1, m+1)$ and $\sigma_{j+1}(c)=(1\ \dots\ m)$. Then $\phi(c)(j) = m+1$ and $\phi(c)(j+1) = m$ and $j$ is a descent of $\phi(c)$.
\item We have $m=n$, $t_j=(l-1, l)$, $t_{j+1} = (l, 1)$ and $\sigma_{j+2}(l \dots n)$. In this case 
$$\phi(c)(j) = l-1 > 1 = \phi(c)(j+1)$$
and $j$ is a descent of $\phi(c)$.
\item We have $m < m+1 < l-1 < l$ and $\sigma_{j+2} = (l \dots m)$. Then $\phi(c)(j) = l-1 > m = \phi(c)(j+1)$ and $j$ is a descent of $\phi(c)$.
\end{enumerate}
In each case, $j$ is a descent of $c$ if and only if $j$ is a descent of $\phi(c)$.

\item We have $t_j = (l-1, m+1)$ and $t_{j+1} = (l-1, l)$. The condition for $j$ to not be a descent of $c$ is that $m+1 < l$, which implies $m+1 < l-1 < l$. Then 
$\phi(c)(j) = m+1$ and $\phi(c)(j+1) = l-1$ and $j$ is also not a descent of $\phi(c)$.
If $j$ is a descent of $c$, then $m+1 > l-1$, which implies $m+1 > m > l > l-1$. In this case, we have $\phi(c)(j) = m$ and $\phi(c)(j+1) = l-1$, which implies that $j$ is also a descent of $\phi(c)$.

\item We have $t_j = (l, m+2)$ and $t_{j+1}=(l, m+1)$. Then $j$ is not a descent of $c$ if
$m+2 < m+1$, which holds only if $m+1=n$. Then $t_{j} = (l,1)$ and $t_{j+1}=(l, n)$ with $\sigma_{j+2} = (l\ \dots\ m)$. Then $\phi(c)(j)=1$ and $\phi(c)(j+1)=n-1$ and $j$ is also not a descent of $\phi(c)$. 

When $j$ is a descent of $c$, then $m+2 > m+1$. In this case, we have $l > m+1$ if and only if $l>m+2$. If $l>m+1, m+2$, then we have $$\phi(c)(j) = m+2 > \phi(c)(j+1) = m+1$$ and if $l<m+1$ we have $$\phi(c)(j) = m+1 > \phi(c)(j+1) = m.$$ 
In both cases, $j$ is a descent of $\phi(c)$.
\end{enumerate}

We have exhausted all possible configurations for $j$, and $t_j, t_{j+1}$, and our proof is complete.

\end{proof}

\subsection{Schur-positivity of convex caterpillars}

\begin{definition}
	Let $c=(t_1,\dots,t_{n-1})$ be a convex caterpillar, and let $i$ be the index of the first edge with $1$ as its endpoint. 
	The edge $t_i$ is called \emph{main edge} of $c$ and the index $i$ is called \emph{main index} of $c$, denoted $\mathtt{I}(c)$.
\end{definition}

For example, for
$c = \left( (4,5),(5,6),(3,6),(1,6),(1,2) \right)$ we have $\mathtt{I}(c) = 4$.

\medskip

Using Lemma~\ref{CaterpillarStructure}, we prove the following explicit description of the descents of a convex caterpillar $c$, based on $\mathtt{I}(c)$ and on the geometry of $c$.

\begin{lemma}\label{DescentBranchLemma1}
	Let $c \in \Ct_n$ and $i \in [n-2]$. Then:
	\begin{enumerate}
		\item 
		For $1 \le i < \mathtt{I}(c)-1$, $i$ is a descent of $c$ if and only if $t_{i+1}$ is a branch of $c$. 
		\item 
		For $i = \mathtt{I}(c)-1$,  $i$ is always a descent of $c$.
		\item 
		For $i=\mathtt{I}(c)$, $i$ is a descent of $c$ if and only if $1$ is not a leaf of $c$.
		\item 
		For $\mathtt{I}(c) < i \le n-2$, $i$ is a descent of $c$ if and only if $t_i$ is a branch of $c$. 
	\end{enumerate}
\end{lemma}
\begin{proof} We prove each case separately
	\begin{enumerate} 
		\item First, suppose that $i = \mathtt{I}(c)$. Then $t_i = (a\ b)$ and $t_{i+1} = (b\ 1)$ for some $2 \leq a, b \leq n$.
		Obviously, $a > 1$, therefore, $t_i$ is a descent. 
		
		\item If $t_{i+1}$ is a branch, then $t_i = (a\ b),\ t_{i+1} = (a\ c)$ for some $a, b, c > 1$. If $(a\ b)$. By Lemma \ref{CaterpillarStructure} $b = a-1$ if $(a\ b)$ is a link or $b = c+1$ if $(a\ b)$ is a branch, the endpoints of the first $i+1$ edges form the cyclic interval $[c, a]$. Since $i < i(c)$ we $1 < c < b < a$, so $i$ is a descent. On the other hand, if $t_{i+1}$ is a link, then $t_i = (a, b)$, $t_{i+1} = (b\ b+1)$, and also $a$ is between $b+1$ and $b$ in $<_b$. Since $1 < a$, we have $a < b$, so $i$ is not a descent.
		\item Now, if $i = \mathtt{I}(c)$ and $1$ is a leaf. Then we have 
		$ t_i = (a\ 1),\ t_{i+1} = (a\ b).$
		Since $1 < b$, $t_i$ is not a descent. In contrast, if $t_i$ is a link, then $t_i = (1\ a),\ t_{i+1}=(1\ b)$ and since both $a$ and $b$ are sorted counterclockwise and are greater than $1$ in cyclic order $<_1$, we have $b < a$, so $i$ is a descent.  
		\item \label{partC} Now suppose that $i > \mathtt{I}(c)$. Then if $t_i$ is a branch, we have 
		$t_i = (a\ b)\ t_{i+1} = (a\ c)$ where $b$ and $c$ are ordered counterclockwise and $a < b, c < n$, which implies that $c < b$ and that $t_i$ is a descent. On the other hand, if $t_i$ is a link, we have $t_i = (a\ a+1),\ t_{i+1} = (a+1\ k)$ where $k > a$, and $i$ is not a descent.
	\end{enumerate}
\end{proof}

Combining  Lemmas~\ref{CaterpillarStructure}  and~\ref{DescentBranchLemma1} with Corollary~\ref{conc:caterpillarstructure2},
we deduce the following key proposition.

\begin{proposition}\label{UniqueLemma}
	A convex caterpillar $c$ is uniquely determined by $\mathtt{I}(c)$ and $\Des(c)$.
\end{proposition}
\begin{proof}
	By Lemma \ref{conc:caterpillarstructure2}, it suffices to show that the pair $(\mathtt{I}(c),\Des(c))$   
	determines 
	the first edge and branches. Note that by Observation~\ref{observation:first_last_edge}, the first and last edges are always branches.
	
	Denote $k:= \mathtt{I}(c)$.
	When $i < k$, combining Lemma~\ref{DescentBranchLemma1} and Lemma~\ref{CaterpillarStructure} allows us to determine whether the $i$-th edge is a branch or a link. By Part 2 of Lemma~\ref{DescentBranchLemma1}, we know whether $1$ is a leaf and whether $e_k$ is a branch. The first edge is determined in both cases by applying parts $2$ and $3$ of Lemma~\ref{CaterpillarStructure}. Branches with indices greater than $k$ are determined by Part $3$ of Lemma~\ref{DescentBranchLemma1}. Therefore, the first edge and the branches are entirely determined by the descent set and $k=\mathtt{I}(c)$ as desired.
\end{proof}

The following lemma describes the possible values of $\mathtt{I}(c)$, given the descent set of $c$. 

\begin{lemma}
	Let $c \in \Ct_n$. Then $\mathtt{I}(c) = 1$ or $\mathtt{I}(c) - 1 \in \Des(c)$.
\end{lemma}
\begin{proof}
	Let $X\subseteq\left[n-2\right]$ and suppose that $\Des(c) = X$. We show that $\mathtt{I}(c) = 1$ or $\mathtt{I}(c)\in X + 1$. 
	It is clear that if $\mathtt{I}(c)\neq 1$ then there exists a $i \in \mathtt{des}_C(1)+1$ such that $\mathtt{I}(c)=i$ according to part $1$ of the lemma~\ref{DescentBranchLemma1}.
\end{proof}

\begin{lemma}
	For every subset $J\subseteq\left[n-2\right]$ and every $i\in\left(1+J\right)\cup\left\{1\right\}$, there exists a unique $c \in \Ct_n$ such that $\Des(c)=J$ and $\mathtt{I}(c)=i$. 
\end{lemma}
\begin{proof}
	Recall that every caterpillar is determined by its first edge and true branches, where every $J\subseteq\left\{2,\dots,n-2\right\}$ can appear as the set of true branches of a caterpillar.
	Placing $\mathtt{I}(c)$ after $x\in J$ results in a proper set of true branches, which defines a caterpillar. Now, suppose that $i\notin J$. Then $(1, 2)$ can be the first edge of the leaf since $1$ is a leaf; hence $1$ is not a descent, and the branches correspond to the members of $X$. If $1\in X$, then $(n\ 1)$ can be the first edge, with the rest of the branches defined by the descents.
\end{proof}

\begin{corollary}\label{conclusion-caterpillarsnumber}
	For every subset $J\subseteq\left[n-2\right]$, the number of convex caterpillars with a descent set $J$ equals $\left| J \right| + 1$.
\end{corollary}

The following observation is well known.


\begin{observation}\label{lemma:hook}
	For every $0 \le k \le n-1$ we have
	$$\{\Des(T) \,:\, T\in \SYT(n-k,1^k)\} = \{J \subseteq [n-1] \,:\, |J|=k\},$$
	each set being obtained exactly once.
\end{observation}

\begin{proof}[Proof of Theorem~\ref{theorem-main}]
	Combine Corollary~\ref{conclusion-caterpillarsnumber} with Observation~\ref{lemma:hook} and Theorem~\ref{theorem-GesselSYTSchur} to deduce
	$$ \Q(Ct_n)=\sum_{k=0}^{n+1}(k+1)\sum\limits_{\substack{J\subseteq [n-1] \\ |J|=k}} \F_{n,J} =
	\sum_{k=0}^{n+1}(k+1)s_{n-k,1^k}.
	$$
\end{proof}

\chapter{Cyclic descent extensions}
\label{chapter: CyclicDescentInvDesClass}
\section{Introduction}\label{section: CyclicDescentExtensions}
This chapter and the next study {\em cyclic descent extensions}.
We will build upon recent work by Adin, Hegedüs and Roichman \cite{AdinHegedusRoichman}.

Recall that the {\em descent set} of a permutation $\pi=\left[ \pi_1,\dots,\pi_n \right]\in\symm{n}$ is
$$\Des(\pi) \coloneqq \left\{ 1 \leq i \leq n-1: \pi_i > \pi_{i+1} \right\} \subseteq \left[n-1\right].$$
Its {\em cyclic descent set} was defined by Cellini in \cite{Cellini} as
$$\CDes(\pi) \coloneqq \left\{1 \leq i \leq n: \pi_i > \pi_{i+1}\right\}\subseteq  \left[n\right],$$
with the convention that $\pi_{n+1} \coloneqq \pi_1$.
\begin{definition}\label{definition-CyclicShift}
The {\em cyclic shift} of a set is the bijection  ${\rm sh} : 2^{[n]} \longrightarrow 2^{[n]}$ defined by
$$\mathrm{sh}(P) \coloneqq \left\{i+1\pmod n: i\in P\right\}.$$
Define also the {\em cyclic shift} of a permutation $\pi \in\symm{n}$ by  $$\rho(\pi)=\left[\pi_n,\pi_1,\dots,\pi_{n-1}\right].$$ 
\end{definition}

\begin{observation} 
$\CDes$ together with $\rho$ satisfy the following 3 conditions.
\begin{enumerate}
	\item $\CDes(\pi)\cap\left[n\right]=\Des(\pi$).
	\item $\CDes\left(\rho\left(\pi\right)\right) = \rm{sh}\left(\CDes\left(\pi\right)\right)$.
	\item For all $\pi\in\symm{n}$, $\emptyset \subsetneq \CDes(\pi) \subsetneq \left[n\right]$.
\end{enumerate}
\end{observation}

We are interested in generalizing this situation as follows.
\begin{definition}\label{definition-CyclicDescent}\cite{ARR}
	Let $\TTT$ be a finite set, equipped with a {\em descent map} 
	$\Des: \TTT \longrightarrow 2^{[n-1]}$. 
	A {\em cyclic extension} of $\Des$ is
	a pair $(\cDes,p)$, where 
	$\cDes: \TTT \longrightarrow 2^{[n]}$ is a map 
	and $p: \TTT \longrightarrow \TTT$ is a bijection,
	satisfying the following axioms:  for all $T$ in  $\TTT$,
	\[
	\begin{array}{rl}
		\text{(extension)}   & \cDes(T) \cap [n-1] = \Des(T),\\
		\text{(equivariance)}& \cDes(p(T))  = {\rm sh}(\cDes(T)),\\
		\text{(non-Escher)}  & \varnothing \subsetneq \cDes(T) \subsetneq [n].
	\end{array}
	\]
\end{definition}
The term non-Escher refers to M.\ C.\ Escher's drawing ``Ascending and Descending"~\cite{WikipediaEscher}, which illustrates the impossibility
of the cases $\cDes(\cdot) = \varnothing$ and $\cDes(\cdot) = [n]$ for permutations in $\symm{n}$.

 As stated earlier, Cellini's $\CDes$(with $p = \rho$) is a cyclic extension of $\Des$ on $\TTT=\symm{n}$. 
 
 Before we proceed, let us stress several properties of cyclic descent extensions. 
\begin{example}
Consider the following subset of $\symm{4}$.
$$T=\left\{ \left[1,2,3,4\right],\left[2,1,3,4\right],\left[2,3,1,4\right],\left[2,3,4,1\right],\left[3,4,1,2\right],[3,1,4,2]\right\}$$
Define $p$ by
$$
p:\left[\begin{array}{cc}
	[1,2,3,4] & \mapsto[2,1,3,4]\mapsto[2,3,1,4]\mapsto[2,3,4,1]\mapsto\left[1,2,3,4\right]\\
	& [3,4,1,2]\mapsto[3,1,4,2]\mapsto\left[3,4,1,2\right]
\end{array}\right]
$$ 
with $cDes$ defined by
\begin{equation*}
\begin{split}
\cDes\left[1,2,3,4\right]&=\left\{ 4\right\} \\
\cDes\left[2,1,3,4\right]&=\left\{ 1\right\} \\
\cDes\left[2,3,1,4\right]&=\left\{ 2\right\} \\
\cDes\left[2,3,4,1\right]&=\left\{ 3\right\} \\
\cDes\left[3,4,1,2\right]&=\left\{ 2,4\right\} \\
\cDes\left[3,1,4,2\right]&=\left\{ 1,3\right\}. 
\end{split}
\end{equation*}
We see that  $(cDes, p)$ is a cyclic descent extension.
We can also define 
$$
\overline{p}:\left[\begin{array}{cc}
	[1,2,3,4] & \mapsto[2,1,3,4]\mapsto[3,4,1,2]\mapsto[2,3,4,1]\mapsto\left[1,2,3,4\right]\\
	& [2,3,1,4]\mapsto[3,1,4,2]\mapsto\left[3,4,1,2\right]
\end{array}\right]
$$
with 
$\overline{\cDes}$ defined by 

\begin{equation*}
	\begin{split}
		\overline{\cDes}\left[1,2,3,4\right]&=\left\{ 4\right\} \\
		\overline{\cDes}\left[2,1,3,4\right]&=\left\{ 1\right\} \\
		\overline{\cDes}\left[3,4,1,2\right]&=\left\{ 2\right\} \\
		\overline{\cDes}\left[2,3,4,1\right]&=\left\{ 3\right\} \\
		\overline{\cDes}\left[2,3,1,4\right]&=\left\{ 2,4\right\} \\
		\overline{\cDes}\left[3,1,4,2\right]&=\left\{ 1,3\right\}. 
	\end{split}
\end{equation*}
Then $(\overline{cDes},\overline{p})$ is also a cyclic descent extension.
\end{example}
However, while cyclic descent extension is not unique, the following holds.
\begin{theorem}\label{theorem:cDistributionUnique}\cite{ARR}
	For a finite set $\TTT$ equipped with descent map $\Des:\TTT\longrightarrow2^[n-1]$, If $(\cDes,p)$ and $(\overline{\cDes}, \overline{p})$ are two cyclic descent extensions of $\Des$ then $\cDes$ and $\overline{\cDes}$ are equidistributed.
\end{theorem}

We restrict our attention to cyclic descent extensions on Schur-positive sets. Recall that a set is Schur-positive if, and only if, it can be represented as a disjoint union of sets, such $\Des$ on every set in the disjoint union is equidistributed with $\Des$ on $SYT(\lambda)$ for some shape $\lambda$. 
Therefore, the study of cyclic descent extensions on Schur-positive sets reduces to the study of cyclic descent extensions on Standard Young Tableaux, which can be reduced to the study of the corresponding quasi-symmetric function $\Q(A)$.

A necessary and sufficient criterion for a skew shape $\lambda / \mu$ that determines whether $\SYT(\lambda/\mu)$ has a cyclic descent extension was found in~\cite{ARR}.

Recall that a {\em ribbon} is a skew shape that doesn't contain a $2 \times 2$ square. 
\begin{example}
\[
\young(::\hfil\hfil,::\hfil,\hfil\hfil\hfil,\hfil)
.
\]  
\end{example}

\begin{theorem}[{\cite[Theorem 1.1]{ARR}}]\label{conj1}
	Let $\lambda/\mu$ be a skew shape with $n$ cells.
	The descent map $\Des$ on $\SYT(\lambda/\mu)$ has a cyclic extension $(\cDes,p)$ if and only if $\lambda/\mu$ is not a connected ribbon.
\end{theorem}
The proof was not constructive.
A constructive proof was later given by 
Huang~\cite{Huang}.

While this solves the problem of the existence of a cyclic descent extension for sets that have descents equidistributed with the $SYT$ of some shape, it does not solve the problem in general. The disjoint union $C=A\uplus B$ of Schur-positive sets $A$ and $B$ can have a cyclic descent extension, while neither $A$ nor $B$ have it.
\begin{example}\label{example:cDesHookExample}
	Let $\lambda = (n-k,1^k)$ and $\mu = (n-k-1, 1^{k+1})$.
	 Then $\Des(\SYT(\lambda))$ is just the set of all subsets of $[n-1]$ of size $k$ and $\Des(\SYT(\mu))$ is the set of all subsets of $[n-1]$ of size $k+1$. It is clear that both $\SYT(\lambda)$ and $\SYT(\mu)$ do not have cyclic descent extension since all their descents have the same size, and shift and restriction to $[n-1]$ must change the size. However their disjoint union $\SYT(\lambda) \uplus \SYT(\mu)$ has a cyclic descent extension, since one can define $cDes(T):=Des(T) \cup \{n\}$ for $T \in SYT(\lambda)$ and $cDes(T):=Des(T)$ for $T \in SYT(\mu)$, with a properly defined bijection $p$. 
\end{example}

On the other hand, the following holds.
\begin{lemma}\label{lemma-DisjointUnionCDes}
	The class of sets with a cyclic descent extension is closed under disjoint union.
\end{lemma}
\begin{proof} if $A$ and $B$ are disjoint sets with $\Des_A : A \to 2^{[n-1]}$ and $\Des_B : B \to 2^{[n-1]}$ such that both $A$ and $B$ have a cyclic descent extension $(\cDes_A, p_A)$ and $(\cDes_B, p_B)$ then $A \uplus B$ has a cyclic descent extension $(\cDes,p)$ defined by 
$$
(\cDes(T),p(T))=\begin{cases}
	\left(cDes_{A}\left(T\right),p_{A}(T)\right) & T\in A\\
	\left(cDes_{B}\left(T\right),p_{B}(T)\right) & T\in B
\end{cases}.
$$
\end{proof}

By Example~\ref{example:cDesHookExample}, if we have $\Q(A)=s_{(n-k,1^k)}+s_{(n-k-1,1^{k+1})}$, or equivalently, the descents of A and ${\SYT(n-k,1^k)\uplus \SYT(n-k-1,1^{k+1})}$ are equidistributed, then $A$ has a cyclic descent extension.

A necessary and sufficient condition on $Q(A)$, under which the Schur-positive set $A$ has a cyclic descent extension and was proved in \cite{AdinHegedusRoichman}. We proceed to describe this condition.
\begin{definition}
The shape $(n-k,1\dots,1)$ is called a {\em hook}. The numbers $k$  and $n-k-1$ are {\em leg length} and {\em arm length} of the hook, respectively.
\end{definition}
\begin{example}
	The hook $h=(5,1,1,1)$ 
	\[
	\young(\hfill\hfill\hfill\hfill\hfill,\hfill,\hfill,\hfill)
	\]
	has leg length $3$ and arm length $4$.
\end{example}

We rely on the following lemma from \cite{AdinHegedusRoichman}.
\begin{lemma}\label{lemma-HookMulPrefixCondition}
	Let $A$ be a Schur-positive set with $\Q(A) = \ch(\phi)$. Denote
	$$
	m_{\lambda}\coloneqq\left\langle \phi,\chi^{\lambda}\right\rangle =\left\langle \Q\left(A\right),s_{\lambda}\right\rangle.
	$$
	Let $h_k$ = $m_{\left(n-k,1^k\right)}$. Then 
	$$
	h_{k}=\left|\left\{ a\in A:\Des\left(a\right)=\left[k\right]\right\}\right| 
	$$ 
\end{lemma}

\begin{definition}
	Let $A$ be a Schur-positive set. Define 
	$$h_k(A) = \left<\Q(A), s_{(n-k,1^k)}\right>$$ and the {\em hook polynomial} of $A$ by
	$$ H_A(x) = \sum\limits_{k=0}^{n-1} h_k(A) x^k. $$
\end{definition}

The following sufficient and necessary condition for cyclic descent extension on a Schur-positive set to exist was proved in \cite{AdinHegedusRoichman}.

\begin{lemma}\label{lemma-CDesExistenceAlgebraicCondition}
	Let $A$ be a Schur-positive set. Then a cyclic descent extension of $Des$ on $A$ exists if and only if $1+x$ divides the hook polynomial $H_A(x)$ and the coefficients of $\frac{H_A(X)}{1+x}$ are non-negative.
\end{lemma}

Before we move on, we state a simple property of hook polynomials that follows from a similar property of quasi-symmetric functions.
\begin{lemma}\label{lemma:HookPolynomialDisjointUnion}
	Let $A$ and $B$ be disjoint Schur-positive sets.
	Then $$H_{A\uplus B}(x) = H_A(X) + H_B(X).$$
\end{lemma}

The rest of this work deals with proof of the existence of cyclic descent extension for particular Schur-positive sets in $\symm{n}$.

Adin, Hegedüs and Roichman in~\cite{AdinHegedusRoichman} solved the problem for permutations with a given cycle type. 
This chapter considers permutations with a given inversion number and permutations with a given inverse descent set. The next chapter deals with proof of the existence of cyclic descent extensions for roots of unity in $\symm{n}$.

\section{Inverse descent classes and their cyclic descent extensions}

We restate Definition~\ref{definition-InvDesClass}
\begin{definition}
	Let $J \subseteq [n-1]$. Define the {\em inverse descent class} of $J$ in $\symm{n}$ by 
	$$D_{J, n}^{-1} = \left\{\pi\in\symm{n}|Des(\pi^{-1}) = J\right\}.$$   
\end{definition}

We are also interested in a generalization of inverse descent classes.
\begin{definition}\label{definition:GeneralizedInvDesClass}
	Let $\J \subseteq 2^{[n]}$. 
	Define {\em generalized inverse descent class}
	$$D_{\J, n}^{-1} = \left\{\pi\in\symm{n}|\Des(\pi^{-1})\in\J\right\}.$$
	or equivalently 
	$$D_{\J, n}^{-1}=\biguplus_{J \in \J}D_{J,n}^{-1}.$$
\end{definition}

It is easy to see that generalized inverse decent classes are Schur-positive since they are disjoint unions of Schur-positive sets.
We are interested in determining the conditions for generalized inverse descent classes to have cyclic descent extensions.

To answer this question, we first calculate 
$H_{D^{-1}_{J, n}}(x)$ for $J\subseteq [n]$.
By discussion after Definition \ref{definition-InvDesClass},
the inverse descent class $D_{J,n}^{-1}$ is the disjoint union of the Knuth classes for standard Young tableaux  $T$ that have $\Des{T}=J$. Namely, we have
$$D_{J,n}^{-1} = \biguplus_{\Des{T}=J}K_T.$$
Now, observe that $Des$ is a bijection between the hooks with leg length $k$ and the subsets of $[n-1]$ of size $k$.
Hence, there is exactly one hook $T$ with $\Des(T)=J$, and it has of shape $(n-|J|, 1^{|J|})$.
Therefore we have $h_{|J|} = 1$ and $h_k = 0$ for $k\neq |J|$.
This result is summarized in the following lemma.
\begin{lemma}\label{Lemma-InvDesHookPolynomial}
	Let $J \subseteq [n]$.
		$$ h_k(D^{-1}_{J, n}) = \delta_{k, |J|}.$$ where $\delta$ is the Kronecker delta.
	Equivalently,
	$H_{D_{J,n}^{-1}}(x)=x^{|J|}$.	
\end{lemma}

Since $1+x$ does not divide $x^k$, we immediately obtain the following result by lemma~\ref{lemma-CDesExistenceAlgebraicCondition}.
\begin{theorem}
	For every $J\subseteq [n-1]$, no cyclic descent extension exists for $D_{J,n}^{-1}$.
\end{theorem}

For the generalized inverse descent class, we obtain the following result:
\begin{lemma}\label{lemma:InverseDecscentCDesLemma}
	Let $\mathcal{J}\in 2^{n-1}$. Let $j_k = \big|\left\{J\in\mathcal{J}\mid |J|=k\right\}\big|$. Then a cyclic descent extension for $D^{-1}_{\mathcal{J}, n}$ exists if and only if
	$$\frac{\sum\limits_{k=0}^{n-1}j_k x^k}{1+x}$$ is a polynomial with non-negative integer coefficients.
\end{lemma}
\begin{proof}
	Recall that $D_{\J, n}^{-1}=\biguplus_{J\in\J}D_{J,n}^{-1}$.
	By lemma~\ref{lemma:HookPolynomialDisjointUnion} we have 
	\begin{equation*}
	\begin{split}
	H_{D_{\mathcal{J}, n}^{-1}}(x)=
	H_{\bigcup\limits_{J \in \mathcal{J} } D_{J, n}^{-1}}(x) 
	= \sum\limits_{J\in\mathcal{J} }H_J(x) 
	= \sum\limits_{J\in\mathcal{J} }x^{|J|}
	\\ = \sum_{k=0}^{n-1}\sum\limits_{\substack{ J\in \mathcal{J} \\ |J| = k} } x^k = \sum\limits_{k=0}^{n-1}j_k x^k
	\end{split}.
	\end{equation*}
	and by lemma \ref{lemma-CDesExistenceAlgebraicCondition} has a cyclic descent extension only if $\sum_{k=0}^{n-1}j_k x^k$ is divisible by $1+x$ and the quotient polynomial has non-negative integer coefficients, as desired.
\end{proof}

Applying the theorem, we obtain the following result.
\begin{theorem}\label{theorem-mainInvDevCDes}
	The following sets in $\symm{n}$ have a cyclic descent extension under the following conditions
	\begin{enumerate}
		\item $D_{2^J, n}^{-1}$ if and only if  $J\neq \emptyset$.
		\item $D_{[I,J],n}^{-1}$ for an interval 
		$$[I, J]=\left\{K: I\subseteq K\ and \ K\subseteq J\right\}$$ (in the poset $2^{[n-1]}$ ordered by inclusion),
		 if and only if  $I\subsetneq J$.
		\item For a maximal chain $C$ in the poset $[I,J]$, $D_{C, n}^{-1}$ has a cyclic descent extension if and only $2 \not| (|J| - |I|)$.
		In particular, $D_{C, n}^{-1}$ for a maximal chain $C$ in $2^{[n-1]}$ has a cyclic descent extension if and only if  $n-1$ is odd, or equivalently $n$ is even.		
	\end{enumerate}
\end{theorem}
\begin{proof}
	\begin{enumerate}
		\item The number of subsets of $J$ of size $k$ is $j_k=\binom{|J|}{k}$. Hence, we have 
		$$\sum_{k=0}^{|J|}j_kx^k=\sum_{k=0}^{|J|}\binom{|J|}{k}x^k=(1+x)^{|J|}$$ which is divisible by $1+x$ unless $|J| = 0$. The quotient 
		$$\frac{(1+x)^{|J|}}{1+x}=(1+x)^{|J-1|}$$ is clearly a polynomial with non-negative integer coefficients. Hence $D_{2^J, n}^{-1}$
		has a cyclic descent extension unless $J=\emptyset$, as desired.
		\item The number of subsets $D^{-1}_{[I,J],n}$ of size $k$ is $\binom{|J|-|I|}{k-|I|}.$ Thus
		$$\sum_kj_kx^k=x^{|I|}(1+x)^{|J|-|I|}$$.
		\item For a maximal chain $C$ in $[I,J]$, we have $$\sum_k j_k x^k = \sum_{k=|I|}^{k=|J|}x^k$$ which is divisible by $1+x$ with quotient having non-negative integer coefficients when $|J|-|I|$ is odd, and is not divisible by $1+x$ otherwise.
	\end{enumerate}
	
\end{proof}

\section{Permutations with given inversion number}
It is known that the set of permutations in $\symm{n}$ with a given inversion number is Schur-positive. See for example \cite{AR1}, Section 9.5.
Hence, it is natural to ask the following question.
\begin{question}\label{question: FixedInvNumberCDes}
	Is there a cyclic descent extension for the set of permutations in $\symm{n}$ with a given inversion number $k$?
\end{question} 
In this section, we give an answer to this question for case $k<n$, proving the following theorem.
\begin{theorem*}[\ref{theorem-mainCDesFixedInversions}]
	For $k < n$, the set of permutations in $\symm{n}$ with $k$ inversions has a cyclic descent extension if and only if $k$ is not a generalized pentagonal number.
\end{theorem*}

Recall the definition of generalized pentagonal numbers.
\begin{definition}\label{definition-generalizedPentagonalNumbers}
	A number $n\in\NN$ is a \emph{generalized pentagonal numbers} if there is $m\in\NN$ such that $n=\frac{(3m\pm1)m}{2}$.
\end{definition}

The rest of this section is devoted to proving this theorem. First, we introduce some notation and definitions. Let $\inv(\pi)$ denote the number of inversions $\pi \in \symm{n}$.Recall that by Definition~\ref{definition:MajorAndInverseMajorIndex} the {\em major index} of a permutation $\pi$ denoted $\maj(\pi)$ is the sum of descents, namely 
${\maj(\pi) = \sum\limits_{i \in Des(\pi)}i}$
and the {\em inverse major index} of $\pi$, denoted $\imaj(\pi)$ is the major index of $\pi^{-1}$, specifically, ${\imaj(\pi) = \sum\limits_{i \in Des(\pi^{-1})} i}$. We define the following sets.
\begin{definition} 
	Let $k\in\NN\cup\{0\}.$ Denote $\inv_k$ the set of permutations in $\symm{n}$ with inversion number $k$, that is
	$$\inv_k(\symm{n})=\{\pi\in\symm{n}:inv(\pi)-k\}.$$
	The set of permutations in $\symm{n}$ with inverse major index $k$ is denoted $\imaj_k(\symm{n})$, namely
	$$\imaj_k(\symm{n})=\{\pi\in\symm{n}:\imaj(\pi)=k\}.$$
\end{definition}

Before going into detail, let us describe the structure of the proof.
First, using a classical result by Foata and Schützenberger from~\cite{FoataSchutz}, we show that the descents of $\inv_k(\symm{n})$ and $\imaj_k(\symm{n})$ are equidistributed, which implies that that $\inv_k(\symm{n})$ admits a cyclic descent extension if and only if $\imaj_k(\symm{n})$ does.

However, $\imaj_k(\symm{n})$ can be expressed as a disjoint union of inverse descent classes, which allows us easily to calculate $H_{\inv_k(\symm{n}}(x)$, the hook polynomial of $\inv_k(\symm{n})$. 

Then we check when the conditions of Lemma~\ref{lemma-CDesExistenceAlgebraicCondition} hold for $H_{\imaj_k(\symm{n})}(x)$. We show then, that divisibility by $(1+x)$ for $n>k$, is equivalent to coefficients of $x^k$ being $0$ in the expansion
$$\prod_{0<n}(1-x^n)=\sum_{0\leq k}a_k x^k,$$
which, by the famous Euler's pentagonal number theorem, implies that $k$ is not a generalized pentagonal number. 

Using the same tools used in the proof of Euler's pentagonal number theorem, we show the second condition of Lemma~\ref{lemma-CDesExistenceAlgebraicCondition} for $H_{inv_k(\symm{n})}$ $(1+x)$ divides $H_{\imaj_k(\symm{n})}$ still holds. That is, the coefficients $\frac{ H_{\inv_k(\symm{n})(x)}}{1+x}$ non-negative which finishes our proof of Theorem~\ref{theorem-mainCDesFixedInversions}.

We proceed to carry out the details of the proof.

We begin citing the classical result of Foata and Schützenberger from \cite{FoataSchutz}. 
\begin{theorem}[\cite{FoataSchutz}, Theorem 1]\label{MajorIndexInvThm}
	There exists a bijection $\phi:\symm{n} \to \symm{n}$ such that $$\Des(\phi(\pi)^{-1}) = \Des(\pi^{-1})$$ and
	$$\inv(\phi(\pi) ) = \maj(\pi).$$
\end{theorem}
Recall that by Theorem~\ref{thm:CoxeterLengthInversionNumber}, the inversion number of $\pi \in \symm{n}$ is equal to its length in Coxeter generators $s_1,\dots,s_{n-1}$ and therefore $\inv(\pi) = \inv(\pi^{-1})$. Using this fact, we show that a bijection exists with similar relations for $\imaj$ and $\inv$.
\begin{corollary}
	There exists a bijection $\psi:\symm{n} \to \symm{n}$ such that 
	$$\Des(\psi(\tau)) = \Des(\tau)$$ and $$\inv(\psi(\tau)) = \imaj(\tau).$$
\end{corollary}
\begin{proof}
Let us define $\psi : \symm{n} \to \symm{n}$ by 
$$ \psi(\tau) = \phi(\tau^{-1})^{-1}.$$ 
By the preceding theorem, we have
$$\Des\left(\psi(\tau)\right) = \Des\left(\phi(\tau^{-1})^{-1}\right)
= \Des(\left(\tau^{-1}\right)^{-1}) = \Des(\tau)$$
and
$$\inv(\psi(\tau)) = \inv\left(\phi(\tau^{-1})^{-1}\right)
	= \inv\left( \phi(\tau^{-1}) \right) = \maj\left(\tau^{-1}\right) = \imaj(\tau),$$
which proves our claim.
\end{proof}

For our purposes, it is more useful to restate it in the following form.
\begin{corollary}\label{corollary:invKimajKEquiDistrib}
	The sets
	$\inv_k(\symm{n})$
	and 
	$\imaj_k(\symm{n})$
	have equidistributed descents.
\end{corollary}

From the preceding corollary, by Definition~\ref{definition-AssociatedQuasisymmetrical}, $\Q\left(\inv_k(\symm{n})\right) = \Q\left( \imaj_k(\symm{n}) \right)$. Therefore, the existence of cyclic descent extension for $\inv_k(\symm{n})$ is equivalent to the existence of cyclic descent extension for $\imaj_k(\symm{n})$.

It is easy to see that $$\imaj_k(\symm{n}) = \left\{\pi \in \symm{n} \big | \maj(\pi^{-1}) = k \right\}$$ is a disjoint union of inverse descent classes, namely
$$\imaj_k(\symm{n}) = \left\{\pi \in \symm{n} \big | \maj(\pi^{-1}) = k \right\} = 
\biguplus_{\substack {J\subseteq [n-1] \\ \sum_{j \in J}j=k} }D_{J,n} ^{-1},$$

which allows us to apply Lemmas~\ref{lemma:HookPolynomialDisjointUnion} and~\ref{Lemma-InvDesHookPolynomial} to obtain
\begin{equation*}
\begin{split}
	H_{\inv_k(\symm{n})}(x) = H_{\imaj_k(\symm{n})}(x) 
	= \sum\limits_{\substack {J\subseteq [n-1] \\ \sum_{j \in J}j=k}  }H_{D_{J,n}^{-1}}(x)
	= \sum\limits_{\substack {J\subseteq [n-1] \\ \sum_{j \in J}j=k} }x^{|J|}.
\end{split}
\end{equation*}

Note that the sets satisfying $J\subseteq [n-1]$ and $\sum\limits_{j\in J} = k$ are the partitions of $k$ into $|J|$ distinct parts of size equal or smaller than $n-1$.

Denote the number of partitions of $k$ into $m$ distinct parts $pd(k, m)$. Let $Pd_k(x)$ be the generating function for the number of partitions of $k$ into distinct parts, namely, the polynomial 
$$Pd_k(x) = \sum pd(k,m)x^m.$$

It is easy to see the following.
\begin{observation}\label{observation-HookPolynomialandPdnkEquivalence}
For $n>k$, $h_m\bigl(\imaj_k(\symm{n})\bigr)=pd(k,m)$ and $$ H_{\imaj_k(\symm{n})}(x)=\sum\limits_{{\substack {J\subseteq [n-1] \\ \sum_{j \in J}j=k} } }x^{|J|} = \sum pd(k,m)x^m = Pd_k(x).$$
\end{observation}

Using the preceding observation, we obtain a necessary condition for the existence of a cyclic descent extension on $\imaj_k(\symm{n})$ and $\inv_k(\symm{n})$ when $k < n$.
\begin{lemma}\label{lemma-imajCycDesExistanceCondition1}
	Let $k,n\in\NN$ such that $k<n$. If a cyclic descent extension on $\inv_k(\symm{n})$ exists, then $a_k=0$ in the expansion
	$$\prod_{1 \leq n}(1-x^n) = \sum_{0 \leq k}a_k x^k.$$
\end{lemma}
\begin{proof}

Let us examine the series $\sum\limits_{k} Pd_k(t)x^{k}$. We have
\begin{equation*}
	\begin{split}
	\sum_{0\leq k} Pd_k(t)x^k &= \sum_{0\leq k} x^k \sum_{m=0}^{k}pd(k,m)t^m = \sum_{0\leq m}\sum_{m \leq k}pd(k,m)t^mx^k \\
	&= \sum_{m\leq 0}t^m\sum_{m \leq k}\left(\sum\limits_{\substack {1 \leq i_1 < \dots < i_m \\ i_1 + \dots + i_m =k }}\prod_{j=1}^mx^{i_j}\right) \\
	& =\sum_{0 \leq m}\left(\sum_{1 \leq i_1 < \dots < i_m}\prod_{j=1}^m tx^{i_j}\right)=\prod_{1\leq n}(1+tx^n).
	\end{split}
\end{equation*}

Substituting $t=-1$ in $\prod\limits_{1 \leq n}(1+tx^n)$, we obtain the equality
$$\prod_{1 \leq n}(1-x^n) = \sum_{0 \leq k}Pd_k(-1)x^k=\sum_{0\leq k}a_k x^k.$$
Since $k<n$, by Observation~\ref{observation-HookPolynomialandPdnkEquivalence}, 
we have $H\bigl(\imaj_k\symm{n}\bigr)(x)=Pd_k(x)$.
By Lemma~\ref{lemma-CDesExistenceAlgebraicCondition}, if a cyclic descent extension for $\imaj_k(\symm{n})$ then $(1+x) | Pd_k(x)$.
But $Pd_k(t)$ if is divisible by $(1+t)$, then $a_k = Pd_k(-1)=0$ proving our claim.
\end{proof}

The following result, known as the pentagonal number theorem, is due to Euler. 
\begin{theorem}[Pentagonal number theorem] We have 
	$$\prod\limits_{1 \leq n} (1-x^n) = 1+ \sum\limits_{1 \leq m}(-1)^m\left(x^{\frac{(3m-1)m}{2}} + x^{\frac{(3m+1)m}{2} }\right)$$
\end{theorem} 
Some notation and arguments from the proof of this theorem, as in \cite[Proposition 1.8.7]{EC1}, are used later. Therefore, we provide a slightly modified version here to ease comprehension. 
\begin{proof}
Let $Q(k)$ be the set of partitions of $k$ into distinct parts. Let $q_e(k)$  be the number of partitions of $k$ into an even number of distinct parts, and let $q_o(k)$ be the number of partitions of $k$ into an odd number of distinct parts. From the proof of Lemma~\ref{lemma-imajCycDesExistanceCondition1}, we have $a_k=Pd_k(-1)$ in the expansion 
$$\prod_{1 \leq n}(1-x^n) = \sum_{0 \leq k}a_k x^k.$$
	
On the other hand, 
$$Pd_k(-1)=\sum_{m=1}^k pd(k,m)(-1)^m =  = q_e(k)-q_o(k) = \sum_{\lambda\in Q(k)}(-1)^{l(\lambda)}.$$
Thus, to prove the theorem, it is enough to show that  
\begin{equation*}
	a_k = q_e(k)-q_o(k)=\left\{ 
		\begin{array}{rl}
		(-1)^m & \text{if}\ k=\frac{(3m\pm 1)m}{2} \text\ {for some}\ m\in \NN\\
		0   & \text{otherwise}.
	\end{array}\right.
\end{equation*}
It is easy to see that $a_0=1$. Therefore the following arguments assume that $k>1$.
To prove the statement $k\neq \frac{(3m\pm 1)m}{2}$, we show that $a_k=0$ by defining an involution $\phi:Q(k) \to Q(k)$ such that $l(\lambda) \not\equiv l(\phi(\lambda))(\bmod2)$ for all $\lambda \in Q(k)$. When $k=\frac{(3m\pm 1)m}{2}$, we will define a partition $\mu\in Q(k)$ with $l(\mu)=m$ and an involution $\phi:Q(k)\setminus\{\mu\}\to Q(k)\setminus\{\mu\}$ such that $l(\lambda)\not\equiv l(\phi(\lambda))(\bmod2)$ for all $\lambda \in Q(k)\setminus\{\mu\}$.
	
For a partition $\lambda \in Q(k)$, let $L_\lambda$ denote the length of the last row of the Ferrers diagram of $\lambda$ and let $D_\lambda$ the number of rows $i$ in $\lambda$ for which $\lambda_i = \lambda_1-i+1$ and $\lambda_i>0$. 
For example for $\lambda = (7,6,5,3,2)$ we have $D_\lambda=3$ and $L_\lambda=2$. 
We define $\phi$ on $Q(k)$ as follows.

If $D_\lambda < L_\lambda$ define $\phi(\lambda) = (\phi(\lambda)_1,\dots,\phi(\lambda)_{l(\lambda)+1})$ by
	
\begin{equation*}
	\phi(\lambda)_i=\left\{ 
	\begin{array}{rl}
	\lambda_i-1 & \text{if}\ 1\leq i \leq D_\lambda\\
	\lambda_i   & \text{if}\ D_\lambda+1 \leq i \leq l(\lambda)\\
	D_\lambda   & \text{if}\  i = l(\lambda)+1.
	\end{array}
	\right.
\end{equation*}
	
If $D_\lambda \geq L_\lambda$. Define $\phi(\lambda) = (\phi(\lambda)_1,\dots,\phi(\lambda)_{l(\lambda)-1})$ by 
\begin{equation*}
	\phi(\lambda)_i=\left\{ 
	\begin{array}{rl}
	\lambda_i+1 & \text{if}\ 1\leq i \leq L_\lambda\\
	\lambda_i   & \text{if}\ L_\lambda+1 \leq i \leq l(\lambda)-1.
	\end{array}
	\right.
\end{equation*}
When $D_\lambda < L_\lambda$ then $\phi(\lambda)\in Q(k)$ unless $l(\lambda)=D_\lambda=L_\lambda-1$. In this case, we have $\lambda=(\lambda_1,\lambda_1-1,\dots,\lambda_1-l(\lambda)+1=l(\lambda)+1)$ and $\phi(\lambda)=(\lambda_1-1,\dots,l(\lambda), l(\lambda))$, and since two smallest parts of $\phi(\lambda)$ are equal, $\phi(\lambda)\notin Q(k)$. Then $k=\frac{(3m+1)m}{2}$ for $m=D_\lambda$. 

When $D_\lambda\geq L_\lambda$ then $\phi(\lambda)\in Q(k)$, unless ${L_\lambda=D_\lambda=l_\lambda}$. In this case, we have ${\lambda=(\lambda_1,\lambda_1-1,\dots,\lambda_1-l(\lambda)+1=l(\lambda))}$ and $\phi(\lambda) = (\lambda_1+1,\lambda_1,\dots,\lambda_1-l(\lambda)+2)$. But then $\phi(\lambda)\notin Q(k)$, since the sum of the parts of $\phi(\lambda)$ is equal to $k-1$. In this case $k=\frac{(3m-1)m}{2}$ for $m=l(\lambda)$.

By the definition of $\phi$, we have ${l(\phi(\lambda))=l(\lambda)\pm1}$. In addition, when both ${\lambda\in Q(k)}$ and ${\phi(\lambda)\in Q(k)}$, then $\phi(\phi(\lambda))=\lambda$. To see that, consider the following geometric interpretation. First, identify $\lambda$ with its Ferrers diagram. Then, if $D_\lambda < L_\lambda$, remove the rightmost cell in each of the first $D_\lambda$ rows of $\lambda$ and add a new row of length $D_\lambda$ below to obtain $\phi(\lambda)$. Otherwise, if $L_\lambda\leq D_\lambda$, remove the last row of $\lambda$, and add a new cell to the right in each one of the first $L_\lambda$ rows of $\lambda$.  It is easy to see that when both $\lambda$ and $\phi(lambda)$ are both in $Q(k)$, applying $\phi$ twice puts the moved cells of the shape $\lambda$ back in their original position without moving other cells and therefore $\phi(\phi(\lambda))=\lambda$. 
As we have shown, $\lambda\in Q(k)$ such that $\phi(\lambda)
\notin Q(k)$ exists if and only if $k=\frac{(3m\pm1)m}{2}$ some 
$m\in\NN$, therefore if $k\neq \frac{(3m\pm1)m}{2}$ for every 
$m\in\NN$, $\phi$ is an involution on $Q(k)$ with the property 
$l(\phi(\lambda))=l(\lambda)\pm1$.

There is no $k\in\NN$ such that $k=\frac{(3m_1+1)m_1}{2}=\frac{(3m_2-1)m_2}{2}$ such that $m_1\neq m_2$ and $m_1,m_2
in \NN$. Therefore, there is at most one $\mu\in Q(k)$ for which one $\phi(\mu)\notin Q(k)$. Therefore, $\phi$ is an involution on $Q(k)$ when $k\neq \frac{(3m\pm1)m}{2}$ and is an involution on $Q(k)\setminus\{\mu\}$ where ${\mu = (2m,\dots,m+1)}$ if $k=\frac{(3m+1)m}{2}$, and ${\mu=(2m-1,\dots,m)}$ if $k=\frac{(3m-1)m}{2}$. It follows that $\phi$ is an involution on $Q(k)\setminus{\mu}$ with the property the property $l(\phi(\lambda))=l(\lambda)\pm 1$ for $k=\frac{(3m\pm1)m}{2}$. Observe that $l(\mu)=m$.

Therefore, for $k\in\NN$ we have
\begin{equation*}a_k=q_e(k)-q_o(k)=\sum_{\lambda\in Q(k)}(-1)^{l(\lambda)}=
	\left\{\begin{array}{rl}
		(-1)^m & \text{if}\ k=\frac{(3m\pm1)m}{2}\ \text{for some}\ m\in\NN\\
		0 & \text{otherwise}
	\end{array}\right.
\end{equation*}
and the proof is finished.
\end{proof} 

Now we proceed to show that if $Pd_k(-1)=0$, then the coefficients of $\frac{Pd_k(t)}{(1+x)}$ are nonnegative integers. We make the following observation by closely examining the proof of the pentagonal number theorem.
\begin{observation}
	For $k\neq \frac{(3m\pm1)m}{2}$, there exists an involution $$\phi:Q(k)\to Q(k)$$ such that $l(\phi(\lambda))=l(\lambda)\pm 1$.
\end{observation}

Using the involution $\phi$, we finish the proof of Theorem~\ref{theorem-mainCDesFixedInversions}.

\begin{theorem}\label{theorem-mainCDesFixedInversions}
	For $k < n$, the set of permutations in $\symm{n}$ with $k$ inversions has a cyclic descent extension if and only if $k$ is not a generalized pentagonal number.
\end{theorem}
\begin{proof}
	First, observe that if $k$ is a generalized pentagonal number, then $Pd_k(-1)\neq 0$ by Euler's Pentagonal Number Theorem, which implies that $(1+t)$ does not divide $Pd_k(t)$.
	
	For the opposite direction, we show that if $k$ is not a generalized pentagonal number, then $\frac{Pd_k(t)}{1+t}$ is a polynomial with integer non-negative coefficients. 
	Let $d$ be the degree of $Pd_k(t)$, which is equal to
	$\max\{l(\lambda)|\lambda\in Q(k)\}$.
	Denote $pd(k,m)=a_m$ and 
	$$b_m=\big|\{\lambda\in Q(k)|l(\phi(\lambda))=l(\lambda)+1\}\big|.$$
	Observe that $a_1 = b_1$, as there is exactly one partition $\lambda=(k)$ of $k$ with $1$ distinct part and we obviously have $l(\phi((k)))=l(k-1,1)=2$. 
	
	Also, we must have $b_d=0$ as $l(\phi(\lambda))=l(\lambda)+1=d+1$ contradicts the fact that $d$ is the maximal possible length of an element $\lambda$ of $Q(k)$.
	
	Since we have $l(\phi(\lambda))=l(\lambda)\pm1$, and  $\phi(\phi(\lambda))=\lambda$, the number of elements $\lambda$ of $Q(k)$ with $l(\lambda)=i$ for which $l(\phi(\lambda))=l(\lambda)-1=i-1$ is equal to the number of elements of $Q(k)$ with length $i-1$ for which $l(\phi(\lambda))=l(\lambda)=(i-1)+1$ which is exactly $b_{i-1}$. Similarily, the number of elements of $Q(\lambda)$ with length $i$ for which $l(\phi(\lambda))=l(\lambda)+1=i+1$ is $b_i$ by definition. Hence, $a_i=b_{i-1}+b_i$.
	 
	Combining everything we have $a_1 = b_1 = 1$, $b_d = 0$, $a_i = b_{i-1}+b_i$.
	Multiplying $1+t$ and $b_1 t+\dots b_{d-1}t^{d-1}$ we have
	\begin{equation*}
		\begin{split}
		(1+t)\sum_{i=1}^{d-1}b_it^i 
		& = b_1 t + \sum_{i=2}^{d-1}(b_{i-1}+b_i)t^i + b_{d-1}t^d\\
		& = a_1 t + \sum_{i=2}^{d-1}a_i t^i + a_d t^d = \sum_{i=1}^d a_i t^i
		  = Pd_k(t).
		\end{split}
	\end{equation*}
	It is easy to see that the coefficients $b_i$ are not negative. By Corollary~\ref{corollary:invKimajKEquiDistrib} and Observation~\ref{observation-HookPolynomialandPdnkEquivalence} we have $H_{\inv_k}(t)=H_{\imaj_k}(t)=Pd_k(t)$. The result follows.
\end{proof}

\chapter{On cyclic descents of roots of unity}
\label{chapter: CyclicDescentUnityRoot}
\section{Background and main results}
Adin, Hegedüs and Roichman characterized the conjugacy classes in the symmetric group, which carry a cyclic descent extension~\cite{AdinHegedusRoichman}.

\begin{theorem}\label{theorem:CDesHegedusCycleStructure}\cite[Theorem 1.4]{AdinHegedusRoichman}
	Let  $\lambda$ be a partition of $n$, and let $C_\lambda \subseteq \symm{n}$ be the corresponding conjugacy class. The descent map $\Des$ on $C_\lambda$ has a cyclic extension $(\cDes, p)$ if and only if $\lambda$ is not of
	the form $(r^s)$ for some square-free $r$.
\end{theorem}

Motivated by this result, we explore the existence of a cyclic descent extension for roots of unity of order $d$ in $\symm{n}$. 
\label{definition-CyclicDescentUnityRoot}
\begin{definition}
	A permutation $\pi \in \symm{n}$ is a \em{root of unity of order $d$}, 
	or \em{$d$-root of unity}, if $\pi^d=e$ where $e$ is the identity permutation in $\symm{n}$. 
	The set of roots of unity order $d$ in $\symm{n}$ will be denoted by 
	$\uroot{d}{n}$.
\end{definition}

The following observation describes the roots of unity of order $d$ in $\symm{n}$.
\begin{observation}\label{observation:PermutationOrder}
	Let $\pi$ be a permutation with cycle decomposition $\pi = \sigma_1\dots \sigma_k$ with corresponding cycle lengths $c_1,\dots,c_k$. Then the order of $\pi$ in $\symm{n}$ is $$ord(\pi)=lcm(c_1,\dots,c_k).$$ 
	Hence, a permutation $\pi$ is a root of unity of order $d$ if and only if $\pi$ is a product of disjoint cycles with lengths that divide $d$.
\end{observation}
Therefore, the set of roots of unity of order $d$ in $\symm{n}$ is a disjoint union of conjugacy classes.

\begin{remark}\label{remark:UrootStructure}
	For convenience, for partition $\lambda=(\lambda_1,\dots,\lambda_k)$, we write $lcm(\lambda)=\lambda(\lambda_1,\dots,\lambda_k$). Using this notation, we have 
	$$\uroot{d}{n}=\biguplus_{lcm(\lambda) | d}C_\lambda.$$
\end{remark}

Calculations suggest the following conjecture.
\begin{conjecture}\label{conjecture-MainConjectureUnityRootCdes}
	A cyclic descent extension for the set of roots of unity of order $d$ in $\symm{n}$ exists if, and only if, $(d,n)\neq 1$.
\end{conjecture}

In this chapter, we show that the $(d,n)=1$ condition is necessary and prove the conjecture for the case when $d=p^k$ for some prime number $p\in\NN$ and $k\in\NN$.
In particular, we prove the following:

\begin{theorem}\label{theorem-mainUrootPrimeVersion}
	The set of roots of unity of order $d$ in $\symm{n}$:
	\begin{enumerate}
		\item Does not have a cyclic descent extension if $n$ and $d$ are coprime.
		\item Has a cyclic descent extension when $d$ is a prime power and $d$ and $n$ are not co-prime.
	\end{enumerate}
\end{theorem}

This chapter is structured as follows. In Section~\ref{section-Ch5BasicPropositions}, we show that the condition $(d,n)\neq 1$ is necessary and show that the set of cyclic roots of unity of order $p^k$ for a prime number $p$ and $k\in\NN$ in $\symm{n}$ admits a cyclic descent extension if and only if the set of roots of order $p$ does. In addition, prove the theorem for the case when $d$ is a power of $2$.

In Section~\ref{section-Ch5UnimodalReduction}, we show that for a Schur-positive subset $A$ of $\symm{n}$, the existence of cyclic descent extension can be reformulated in terms of multiplicities unimodal permutations in $A$ instead of hook-lengths in $\Q(A)$. In particular, we introduce the \emph{unimodal maxima polynomial} $U_A(x)$, which counts the unimodal permutations in $A$ by their maxima. We show that $A$ has a cyclic descent extension if and only if $U_A(x)$ satisfies the same conditions of Lemma~\ref{lemma-CDesExistenceAlgebraicCondition} that $H_A(x)$ does. This equivalence allows us to attack the problem of the existence of cyclic descent extension for $\symm{n}$ using the existing tools and results for unimodal permutations in $\symm{n}$.

The last sections are devoted to proving the Conjecture~\ref{conjecture-MainConjectureUnityRootCdes} for the case when $d$ is a power of a prime number $p$, which has three stages. First, in Section~\ref{section-ch5PrimeCase}, we prove the conjecture for the case when $d=n$ for a prime number $d$. After that, in Section~\ref{section-Ch5Gluing}, we explore the structure of unimodal permutations and their cycle structure.
In Section~\ref{section-Ch5FinalProof}, we explore the implications of the cycle structure on the distribution of the maxima. Then we derive the divisibility conditions of Lemma~\ref{lemma-CDesExistenceAlgebraicCondition} for $U_{\uroot{d}{n}}(x)$, finishing the proof of Theorem~\ref{theorem-mainUrootPrimeVersion}.

\section{Basic propositions}\label{section-Ch5BasicPropositions}
In this section, we show the necessity of $d$ and $p$ being coprime for the existence of a cyclic descent extension for $\uroot{d}{n}$. We also show that the existence of a cyclic descent extension for $\uroot{p}{n}$ for a prime $p$ implies the existence of a cyclic descent extension for $\uroot{p^k}{n}$ for every $k\geq 1$. We conclude this section with proof of the Theorem~\ref{theorem-mainUrootPrimeVersion} for the case of $d=2^k$, which doesn't require the machinery developed in the later sections of this chapter.

First, we prove that the condition that $d$ and $n$ are not coprime is necessary.

\begin{lemma}\label{lemma:easyMainConjecturePart}
	If $d$ and $n$ are coprime, then $\uroot{d}{n}$  has no cyclic descent extensions.
\end{lemma}
\begin{proof}
	By remark~\ref{remark:UrootStructure} we have 
	$$\uroot{d}{n}=\biguplus_{lcm(\lambda)|d}C_\lambda.$$
	But $lcm(\lambda) | d$ for $\lambda=(r^s)$ for square-free $r$ means that $r|d$, and since necessarily $r|n$ and  $(n,d)=1$ we must have $r=1$ and $\lambda=(1^n)$,
	which means $C_\lambda= C_{(1^n)}=\{e\}$. Otherwise, $\lambda\neq(r^s)$ for every square-free $r$ and $s\in\NN$. Therefore, by Theorem~\ref{theorem:CDesHegedusCycleStructure}, for every $\lambda$, when $lcm(\lambda)|d$,
	$C_\lambda$ has a cyclic descent extension, unless $C_\lambda=\{e\}$.	
	Denote $\biguplus_{lcm(\lambda)|d,C_\lambda\neq\{e\}}C_\lambda=C.$ 
	By Lemmas~\ref{lemma-DisjointUnionCDes} and~\ref{lemma:HookPolynomialDisjointUnion} $C$ has a cyclic descent extension and $(1+x)|H_C(x)$. On the other hand $H_{\{e\}}(x)=1$. Obviously
	$\uroot{d}{n}=C\uplus\{e\}$, and 
	$$H_{\uroot{d}{n}}(x)=H_C(x)+H_{\{e\}}(x)=H_C(x)+1.$$
	Since $(1+x)|H_C(x)$, $(1+x)\not| (H_C(x)+1$, which implies that $\uroot{d}{n}$ has no cyclic descent extension.
\end{proof}

Now we show that when $d=p^k$ for some prime number $p$, to prove the cyclic descent extensibility of $\uroot{d}{n}$, it is enough to prove the cyclic descent extensibility of $\uroot{p}{n}$.  
\begin{lemma}
	Assume that $p|n$, $p$ is prime and a cyclic descent extension exists for $\uroot{p}{n}$. Then for every $d=p^k$ with $k>1$, $\uroot{d}{n}$ has a cyclic descent extension.
\end{lemma}
\begin{proof}
	Let $C$ denote $\uroot{d}{n} \setminus \uroot{p}{n}$. If $C=\emptyset$ then we are done, since $\uroot{d}{n}=\uroot{p}{n}$ in this case.
	Otherwise, we have
	$$C=\biguplus_{\substack{lcm(\lambda) | p^k \\ lcm(\lambda) > p} }C_\lambda.$$
	But every conjugacy class $C_\lambda$ in the union has a cycle of length $p^t$ for $t > 1$, which is not square-free, namely $\lambda$ is not of the form $(r^s)$ with square-free $r$. It follows that $C_\lambda$ has a cyclic descent extension, which implies the same for $C$ by Lemma~\ref{lemma-DisjointUnionCDes}. By  assumption, $\uroot{p}{n}$ has a cyclic descent extension, and so does $\uroot{d}{n}$, as a disjoint union of $C$ and $\uroot{p}{n}$.
\end{proof}

We proceed to prove Theorem~\ref{theorem-mainUrootPrimeVersion} when $d=2$. We do it now, as the proof does not need the tools we develop in later sections. In addition, we can later assume that $d\neq 2^k$ when proving the general case, thus significantly simplifying the proof.
\begin{definition}
	A permutation $\pi$ in $\symm{n}$ is called an {\em involution} if it is a root of unity of order $2$. 
\end{definition}

\begin{remark}In our notation, the set of involutions is just $\uroot{2}{n}$.
\end{remark}
\begin{theorem}\label{theorem-CdesInvolutions}
	The set $\uroot{2}{n}$ has a cyclic descent extension if and only if $n$ is even.
\end{theorem}

\begin{proof}
	Recall that $RSK$ is a bijective function between permutations in $\symm{n}$ and pairs of Standard Young Tableaux $(P,Q)$ of the same shape, such that if $RSK(\pi)=(P,Q)$ then $RSK(\pi^{-1}) = (Q, P)$ and $\Des(\pi)=\Des(Q)$. Since, $\pi$ is an involution if and only if  $\pi^{-1}=\pi$, we have from the properties of $RSK$ that $\pi$ is an involution if and only if  $RSK(\pi)=(P,P)$ for some $P$. Hence, we have $Des(\uroot{2}{n})=Des(SYT(n))$ as multisets. Therefore,
	$$\Q(\uroot{2}{n})=\sum_{\lambda \vdash n}s_\lambda,$$
	which obviously implies that $s_{(n-k,1^k)}$ appears exactly once in $\Q(\uroot{2}{d})$. Hence we have, 
	$$H_{\uroot{2}{n}}(x)=\sum_{i=0}^{n-1} x^i$$ which is divisible by $(1+x)$ if and only if  $n$ is even, as claimed.
\end{proof}

\section{Reduction to unimodal permutations}\label{section-Ch5UnimodalReduction}
In this section, we show that for a Schur-positive subset $A$ of $\symm{n}$, the study of hook multiplicities in $\Q(A)$ can be reduced to the study of multiplicities of unimodal permutations in $A$. 
First, we prove a slightly different form of lemma~\ref{lemma-HookMulPrefixCondition}.
\begin{lemma}\label{lemma:HookMulSuffixCondition}
	Let $A$ be a Schur-positive set with $\Q(A) = \ch(\phi)$ and $h_k$ defined as in lemma~\ref{lemma-HookMulPrefixCondition}. Then 
	$$h_k = \big|\bigl\{ a\in A:\Des(a)=\left\{ n-k,\dots,n-1\right\} \bigr\}\big|.$$ 
\end{lemma}
\begin{proof}
	Let $T\in \SYT(n-k,1^{k})$ be the hook with $1,\dots,n-k$ in the first row, and ${1,n+k+1,\dots,n}$ in the first column.
	$$T = \ytableausetup{mathmode,boxframe=normal,boxsize = 2em}
	\begin{ytableau}
		1& \none[\dots]& {\scriptstyle n-k}\\
		{\scriptscriptstyle n-k+1}\\
		\none[\vdots]\\
		n
	\end{ytableau}.$$
Clearly, we have $\Des(T)=\{n-k,\dots,n-1\}.$.
Hence, $$h_k \leq \big|\bigl\{ a\in A:\Des(a)=\left\{ n-k,\dots,n-1\right\} \bigr\}\big|.$$
Now, if 
$\Des(S)=\left\{n-k,\dots,n-1\right\}$ for some standard Young tableau $S$, then the numbers $1,\dots,n-k$ must appear in the first row, since $i+1$ appearing south of $i$ would imply that $i\in\Des(S)$. 
Now, $n-k\in\Des(S)$ implies that $n-k+1$ appears south of $n-k$ and can not appear in the first row. On the other hand, it must appear in position $(2,1)$; otherwise, there must be a number larger than $n-k+1$, which appears west or north of $n-k$, contradicting $S$ being a standard Young tableau. By induction, $n-k+j$ appears in position $(j+1,1)$ for every $1\leq j \leq k$. Hence, we must have $S=T$. Hence, there exists a unique standard Young tableau $T$ with 
$Des(T)=\left\{n-k,\dots,n-1\right\}$ . 
Since $A$ is Schur-positive we have by Theorem~\ref{theorem-Sp-combin} 
$$\sum_{a\in A} x^{\Des(A)}=\sum_{\lambda \vdash n}c_\lambda\sum_{T\in\SYT(\lambda)}x^{\Des(T)},$$
 where $c_\lambda$ is the multiplicity of $s_\lambda$ in $\Q(A)$.
The desired result follows by comparing coefficients of $x^{\left\{n-k,\dots,n-1\right\} }$ on both sides of the equation.
\end{proof}
For a permutation $\pi\in\symm{n}$, the condition
$$\Des(\pi)=\left\{ n-k,\dots,n-1\right\}$$
implies that $\pi_i < \pi_{i+1}$ for $1<i<n-k$ and $\pi_i>\pi_{i+1}$ for $i \geq n-k$. Namely, if we view $\pi$ as a finite sequence, it is an ascending sequence up to the position $n-k$ and a descending sequence after. Recall that such sequences are called {\em unimodal.} In particular, the last lemma can be restated for permutations in the following form.
\begin{lemma}\label{lemma:UnimodalHookEquality}
	Let $A$ be a Schur-positive subset of $\symm{n}$. Then, $h_k$ equals the number of unimodal permutations in $A$ that have the maximum at $n-k$.
\end{lemma} 
Therefore, the calculation of $h_k(\uroot{d}{n})$ reduces to counting the number of unimodal permutations that have the maximum at $n-k$.
\begin{definition}\label{definion:UnimodalMaximaPolynomial}
	For a subset $A$ of $\symm{n}$, let $u_k(A)$ denote the number of unimodal permutations in $A$ that have the maximum at $k$. Define the {\em unimodal maxima polynomial} of $A$ by
	$$U_A(x)=\sum_{k=1}^n u_k(A) x^k.$$
\end{definition}
The following lemma allows us to work with $U_A(x)$ instead of $H_A(x)$ to check an existence of a cyclic descent extension for a Schur-positive set $A$, as in Lemma~\ref{lemma-CDesExistenceAlgebraicCondition} 
\begin{lemma}
	Let $A$ be a Schur-positive subset of $\symm{n}$ with unimodal maxima and hook polynomials $U_{A}(x)$ and $H_A(x)$ respectively.  Then $U_A(x)$ is divisible by $1+x$ with the quotient having non-negative coefficients if and only if $H_A(x)$ is.
\end{lemma}
\begin{proof}
	By Lemma~\ref{lemma:UnimodalHookEquality}, we have $h_k(A)=u_{n-k}(A)$,
	which implies the equality $U_A(x)=x^n H_A(\frac{1}{x})$. If $$H_A(x)=(1+x)(a_0+\dots a_{n-1}x^{n-1}),$$ then 
	\begin{equation*}
		\begin{split}
			U_A(x) &= x^n(1+\frac{1}{x})(a_0+\dots a_{n-1}\frac{1}{x^{n-1} }) \\
			&= (1+x)(a_0 x^{n-1}+\dots a_{n-1})
		\end{split}
	\end{equation*}
	Hence, divisibility by $1+x$ of $H_A(x)$ with quotient coefficients being non-negative implies the same for $U_A(x)$. The converse follows similarly.
\end{proof}
\begin{corollary}\label{corolloary: UnimodalPolynomialCdesCondition}
	A Schur-positive subset $A$ of $\symm{n}$ has a cyclic descent extension if and only if $U_A(x)$ is divisible by $1+x$ with coefficients of the quotient being non-negative integers.
\end{corollary}

\section{The case of prime $n$}\label{section-ch5PrimeCase}
In this section, we prove the existence of a cyclic descent extension for the set $\uroot{p}{p}$ when $p$ is a prime number.

Before we proceed to calculations, we cite a result by Gannon from~\cite{Gannon-Unimodal}. Let $\unimodal_{n}^{m}$ denote the set of all unimodal $n$-cycles with the maximum at $m$.
\begin{lemma}\label{lemma:unimodalNCycleWithMaxMultiplicity}
	The number of unimodal cycles of length $n$, with maximum at $m$ is
	\begin{align*}
		|\unimodal_{n}^{m}| &=\delta_{n,1}+
		\frac{(-1)^{n-m}}{n}\sum_{\substack{d | n \\ d\leq n-m}}\mu(d)\sum_{j=1}^{(n-m)/d}(-1)^j\binom{n/d}{j} \\
		&=\delta_{n,1}+(-1)^{m+1}\delta_{n,2}+
		\frac{(-1)^{n-m+1}}{n}\sum_{\substack{d | n \\ d <  m}}\mu(d)\sum_{j=1}^{\left\lceil m/d\right\rceil -1}(-1)^{j+n/d}\binom{n/d}{j}.
	\end{align*}
	where $\mu$ is the Möbius function on $\NN$.
\end{lemma}

\begin{observation}\label{observation:BinomialIdentity} 
	Binomial coefficients satisfy the following identity
	$$
	\sum_{i=0}^k(-1)^i\binom{n}{i}=(-1)^k\binom{n-1}{k}.
	$$
	or equivalently
	 $$\sum_{i=1}^k\binom{n}{i}(-1)^i=(-1)^k\binom{n-1}{k}-1.$$
\end{observation}

\begin{lemma}
	When $p$ is prime and $n=p$, the polynomial $U_{\uroot{n}{n} }\left(x\right)$ is divisible by $1+x$ and the quotient has non-negative coefficients. 
\end{lemma}
\begin{proof}
	Observe that roots of unity of order $p$ in $\symm{p}$ are either cycles of length $p$ or the identity permutation.
	Replace
	 $$\sum\limits_{j=1}^{(n-m)/d}(-1)^{j}\binom{n/d}{j}$$ 
	 by
	  $$(-1)^{(n-m)/d}\binom{n/d-1}{(n-m)/d}-1$$ in 
	$$|\unimodal_{n}^{m}| =\delta_{n,1}+
	\frac{(-1)^{n-m}}{n}\sum_{\substack{d | n \\ d\leq n-m}}\mu(d)\sum_{j=1}^{(n-m)/d}(-1)^j\binom{n/d}{j}$$
	to obtain
	$$
	|\unimodal_n^m| = \delta_{n,1}+
	\frac{(-1)^{n-m}}{n}\sum_{ \substack{d | n \\ d \leq n-m} }
	\mu(d)\left((-1)^{(n-m)/d}\binom{n/d-1}{(n-m)/d}-1\right).
	$$
	Note that no unimodal cycles of length $p$ have their maximum at $p$; the only unimodal permutation with the maximum at $p$ is the identity permutation. The case $p=2$ has already been handled in Theorem \ref{theorem-CdesInvolutions}. Therefore, we assume that $p \geq 3$ is an odd prime. Since $p$ is prime, the only divisor $d$ of $p$ with $d < p-m$ for $m \geq 1$ is $d=1$ with $\mu(d)=1$ and obviously $\delta_{p,1}=0$. 
	Thus
	\begin{equation*}
		\begin{split}
			\Delta_p^m 
			&= \frac{1}{p}\left(\binom{p-1}{p-m}-(-1)^{p-m}\right) \\
			&= \frac{1}{p}\left(\binom{p-1}{m-1}+(-1)^{m}\right).
		\end{split}
	\end{equation*}
	Substituting $u_m(\uroot{p}{p})=\unimodal_p^m$ for $1 \leq m \leq p-1$ and $u_p(\uroot{p}{p})=1$ we obtain
	\begin{equation*}
	\begin{split}
		U_{\uroot{p}{p}}(x)
			&= x^p+\frac{1}{p}\sum_{m=1}^{p-1}\left(\binom{p-1}{m-1}
			+(-1)^m\right)x^m\\
			&= x^p+\frac{1}{p}\left(\sum_{m=1}^p\binom{p-1}{m-1}x^m
			+\sum_{m=1}^{p}(-1)^mx^m\right)\\
			&=\frac{1}{p}\left(\sum_{m=1}^p\binom{p-1}{m-1}x^m+px^p
			+\sum_{m=1}^{p}(-1)^m x^m\right)\\
			&=\frac{1}{p}\left(\sum_{m=1}^{p}\binom{p-1}{m-1}x^m
			+\sum_{m=1}^{p-1}\left(x^p-(-1)^mx^m\right) \right)\\
			&=\frac{1}{p} \left( \left(1+x\right)^{p-1} +
			 \left(1+x\right)\sum_{m=1}^{p-1}\sum_{j=m}^{p-1}(-1)^j x^j
			\right) \\
			&=\frac{1+x}{p}\left(\sum_{m=1}^{p-1}\binom{p-2}{m-1}x^m
			+\sum_{m=1}^{p-1}(-1)^m mx^m\right)\\
			&=(1+x)\sum_{m=1}^{p-1}\frac{1}{p}\left(\binom{p-2}{m-1}+(-1)^m m\right)x^m.
		\end{split}
	\end{equation*}
	We have to show that the coefficients in the summation are non-negative integers.
	For $2 \le m\leq p-2$, $\binom{p-2}{m-1}\geq p-2$ and also $\binom{p-2}{0} - 1 = 0$ (for $m=1$) and $\binom{p-2}{p-2}+(p-1)=p$ (for $m=p-1$). Therefore the coefficient is non-negative, and it remains to prove that they are integers. It is sufficient to show that
	$$\binom{p-2}{m-1}+(-1)^{m}{m}\equiv 0\Pmod{p}$$ for every $1\leq m \leq p-1$ or equivalently, 
	$$\binom{p-2}{m}+(-1)^{m+1}(m+1)\equiv0\Pmod{p}$$ for every $0 \leq m \leq p-2$.
	
	We proceed by reverse induction. For $m=p-2$ 
	$$\binom{p-2}{p-2}+(-1)^{p-1}(p-1)=p\equiv 0\Pmod{p}.$$ 
	Assume that for some $m<p-2$,
	$$\binom{p-2}{m+1}\equiv(-1)^{m+1}(m+2)\Pmod{p}$$ holds.
	Recall that
	$$\binom{p-2}{m}=\binom{p-2}{m+1}\frac{m+1}{p-2-m}.$$
	But 
	$$\binom{p-2}{m+1}\equiv (-1)^{m+1}(m+2)\Pmod{p},$$ therefore
	\begin{equation*}
		\begin{split}
			\binom{p-2}{m} 
			& \equiv (-1)^{m+1}(m+2)\frac{m+1}{p-2-m}\Pmod{p} \\
			& \equiv (-1)^{m+1}(m+2)\frac{m+1}{-(m+2)}\Pmod{p} \\
			& \equiv (-1)^m(m+1)\Pmod{p}
		\end{split}
	\end{equation*}
	as desired.
	We have proved that	the coefficients of $\frac{U_{\uroot{p}{p}}(x) }{1+x}$ are non-negative integers.
\end{proof}

Applying Corollary~\ref{corolloary: UnimodalPolynomialCdesCondition}, we immediately obtain the following.
\begin{corollary}
	For every prime $p$, the set $\uroot{p}{p}$ has a cyclic descent extension.
\end{corollary}


\begin{corollary}\label{corollary:involutionUnimodal}
	There exists an involution  $$\phi:\uroot{p}{p}\cap\unimodal(p)\to\uroot{p}{p}\cap\unimodal(p)$$
	such that $m_{\phi(\pi)} = m_\pi \pm 1$ where $m_{\phi(\pi)} = m_\pi \pm 1$ and $m_{\phi(\pi)}$ denote the positions of maxima of $\phi(\pi)$ and $\pi$ respectively.
\end{corollary}

\section{Gluing unimodal permutations}\label{section-Ch5Gluing}
In this section, we study the cycle structure of unimodal permutations. Intuitively, we show the existence of a total ordering $<_\Delta$ on \emph{unimodal full cycles} and define a partial action $\oplus$ on unimodal full cycles, such that all unimodal permutations can be described as $\sigma_1\oplus\dots\oplus \sigma_k$ where $\sigma_1,\dots,\sigma_k$ are unimodal full cycles. We then proceed to show that $\sigma_1,\dots,\sigma_k$ is defined if any only in $\sigma_1,\dots,\sigma_k$ is a unimodal sequence relative to $<_\Delta$.

Most of the definitions, results and basic notation in this section are based on the work of Gannon~\cite{Gannon-Unimodal} on the structure of unimodal permutations. Details and further references can be found there. There are
other approaches to the enumeration of unimodal permutations with given cycle structure and maximum, for example, in the work of Thibon~\cite{Thibon}, who employed methods similar to those of Gessel and Reutenauer~\cite{Gessel-Reutenauer}.

\medskip

We use the following notation. $\unimodal(n)$ denotes the set of unimodal permutations in $\symm{n}$. The set of unimodal {\em full cycles} in $\symm{n}$ will be denoted by $\unimodal_n$. Denote by $\unimodal(*)=\biguplus\limits_{n\in \NN}\unimodal(n)$ the set of all unimodal permutations, and $\unimodal_* = \biguplus\limits_{n\in\NN}\Delta_n$ the set of all unimodal full cycles. The {\em maximum} of a permutation $\pi\in\symm{n}$ is $m_\pi=\pi^{-1}(n)$.

For a linearly ordered set $J=\left\{j_1<\dots<j_n\right\}$,
let $\Omega_J:[n]\to J$ be the unique order preserving bijection, namely $\Omega_J(i)=j_i$.
For a linearly ordered set $X$ and a function $f:X\to X$, when $f$ preserves a finite $J\subseteq X$ and $f|_J$ is a bijection, the {\em permutation type} or {\em shape} of $f|_J$ is the permutation in $\pi\in\symm{n}$ obtained by relabeling $j_i$ with $i$ or equivalently $\pi = \Omega_J^{-1}\circ f|_J \circ \Omega_J$.

We proceed to describe the key results from~\cite{Gannon-Unimodal} that we use in our work.

\begin{lemma}\cite{Gannon-Unimodal}\label{lemma:unimodalSubCycle}
	Let $\pi \in \unimodal(n)$ and $J\subseteq [n]$ such that $\pi$ preserves $J$. Then $\pi|_J$ is unimodal.
	In particular, if $\pi= \sigma_1\dots \sigma_k$ is the cycle decomposition of $\pi$, then $\sigma_i$ is unimodal for every $1 \leq i \leq k$.
\end{lemma}

\begin{definition}\label{defintion:JSum}
Let $\pi_i \in \unimodal(n_i)$ for $i=1,\dots,k$, $n=n_1+\dots+n_k$ and disjoint $J_i\subseteq [n]$ with $|J|_i=n_i$, and $\Omega_{J_i}:[n_i]\to J_i$ be the unique order preserving bijections. Define the permutation $\pi \in S_n$ by
$$\pi(x)=
\left(\pi_{1}\oplus_{J_1}\dots\oplus_{J_{k-1} } \pi_k\right)(x)
=\left(\Omega_{J_i}\circ\pi_i\circ\Omega_{J_i}^{-1}\right)(x)\ \text{for}\ x\in J_i.
$$
Note, that $\pi$ is always well defined, because $\left\{J_1,\dots J_k,\right\}$ is a set partition of $[n]$. 
\end{definition}
\begin{remark}
	An alternative notation will be used.
	When $k=2$, we sometimes write $J$ and $\overline{J}$ instead of $J_1$ and $J_2$ respectively. When we do that, we also write $\Omega$ and $\overline{\Omega}$  instead of $\Omega_{J_1}$ and $\Omega_{J_2}$.
	Often, we write $\Omega_i$ instead $\Omega_{J_i}$ and $\oplus_i$ instead of
	 $\oplus_{J_i}$. When the context is clear, the indices may be omitted.
\end{remark}

Using the $\oplus$ notation, we make the following observation.
Let $\pi\in\unimodal(*)$ and let ${\pi=\sigma_1\circ\dots\circ \sigma_k}$ be the cycle decomposition of $\pi$. Let $J_i$ be the subset of $[n]$ preserved by $\sigma_i$, $n_i$ be the maximum of $J_i$ and $m_i = \sigma_i^{-1}(n_i)$. Let $\pi_i$ be the shape of $\sigma_i$.

Assume that $m_i < m_{i+1}$, since it is always possible to enumerate the cycles in this way.
Then, $\pi=\pi_1\oplus_{J_1}\dots\oplus_{J_{k-1}}\pi_k$ and the sequence $n_1,\dots,n_k$ has to be unimodal.
We obtain the following observation.
\begin{observation}\label{observation:unimodalCycleSumStructure}
Every unimodal permutation $\pi\in\unimodal(*)$ can be uniquely represented as 
$\pi_1\oplus_{J_1}\dots\oplus_{J_{k-1}}\pi_k$ with $\Omega_i(m_i)<\Omega_{i+1}(m_{i+1})$ and the sequence $\Omega_i(n_i)$ is unimodal.
\end{observation}

From now on, unless stated otherwise, when  
$\pi=\pi_1\oplus_1\dots\oplus_{k-1} \pi_{k}$ is unimodal and $m_i$ denote the maxima of $\pi_i$, then $i<j$ implies that $\Omega_i(m_i) < \Omega_j(m_j)$.
In particular, when indices are omitted and the sum $\pi_1\oplus\dots\oplus\pi_k$ is unimodal, we always mean that $\Omega_i(m_i) < \Omega_j(m_j)$ for $i<j$.

The lemma below is essential for building unimodal permutations from smaller ones. A simple observation is that if $\pi = \pi_1 \oplus_J \pi_2$ holds, then the maximum $m$ of $\pi$ is equal to either $\Omega(m_1)$ or $\overline{\Omega}(m_2)$ where $m_i$ denotes  the maximum of $\pi_i$.
\begin{lemma}\cite[Lemma 1]{Gannon-Unimodal}\label{lemma:UnimodalSum}
	Let $\pi_1\in\unimodal(n_1), \pi_2\in\unimodal(n_2)$, $n=n_1+n_2$ and $J\subseteq [n]$ with $|J|=n_1$ such that $\Omega(m_1) < \overline{\Omega}(m_2)$. Let $\pi=\pi_1\oplus\pi_2$. Use $m_1,m_2,m$ to denote maxima of $\pi_1,\pi_2,\pi$ respectively. Assume that $\pi_1\oplus\pi_2\in\unimodal(n)$. Then for every $i_1\in[n_1]$ and $i_2\in[n_2]$ we have:
	\begin{enumerate}
		\item If $i_1 > m_1$ and $i_2 \geq m_2$, then $\Omega(i_1) > m$ and $\overline{\Omega}(i_2) \geq m$.
		\item If $i_1 \leq m_1$ and $i_2 < m_2$ then $\Omega(i_1) \leq m$ and $\overline{\Omega}(i_2) < m$.  
		\item If $i_1 > m_1$ and $i_2 < m_2$ then $\Omega(i_1) > \overline{\Omega}(i_2)$.
		\item If $i_1 \leq m_1$ and $i_2 \geq m_2$ then $\Omega(i_1) < \overline{\Omega}(i_2)$.
	\end{enumerate}
\end{lemma}
We have the following relation on the maxima unimodal sums, which we use extensively in our calculations.
\begin{corollary}\label{corollary:UnimodalSumMax}
	Let $\pi_1,\pi_2,\pi,m_1,m_2,m,\Omega,\overline{\Omega}$ be as in Lemma~\ref{lemma:UnimodalSum}. Then
	\begin{enumerate}
	\item $m=m_1+m_2-1$ if $m=\Omega(m_1)$.
	\item $m=m_1+m_2$ if $m=\overline{\Omega}(m_2)$.
	\end{enumerate}
\end{corollary}
\begin{proof}
	By parts 1 and 2 of \ref{lemma:UnimodalSum}, we have 
	$\big|\{i\in[n_1]:\Omega(i)\leq m\}\big|=m_1$
	and
	$$
	\big|\{j\in[n_2]:\overline{\Omega}(j) \leq m\}\big|=\left\{ 
		\begin{array}{rl}
			m_2 - 1 & \text{when}\ \Omega(m_1) = m \\
			m_2 & \text{when}\ \Omega(m_2) = m.
		\end{array}	
	\right.
	$$
	Next, observe that 
	$$m=\big|\{i\in[n_1]:\Omega(i)\leq m\}\big| + 
	\big|\{j\in[n_2]:\overline{\Omega}(j)\leq m\}\big|.
	$$
	When $m=\Omega(m_1)$, we have $\big|\{j\in[n_2]:\overline{\Omega}(j) \leq m\}\big|=m_2-1$, and therefore $m=m_1+m_2-1$. Otherwise, we have $m=m_1+m_2$, as desrired.
\end{proof}

The following theorem is essential for building unimodal permutations with a given cycle structure.
\begin{theorem}\cite[Theorem 3]{Gannon-Unimodal}\label{theorem-unimodalCycleRelativeOrder}
Let $\pi_1\in\unimodal_{n_1}$ and $\pi_2\in\unimodal_{n_2}$ be unimodal full cycles with maxima $m_1$ and $m_2$ respectively.
There exists a unique $J$ with $\pi_1\oplus_J\pi_2$ unimodal satisfying $\Omega(m_1)<\overline{\Omega}(m_2)$. 
\end{theorem}
\begin{definition}\cite[discussion before Theorem 3]{Gannon-Unimodal}\label{definition:cyclePairElementsOrder}
	Let $\pi_1\in\unimodal_{n_1}$ and $\pi_2\in\unimodal_{n_2}$ be two unimodal cycles. Let $J$ and $\Omega$ be as in Theorem~\ref{theorem-unimodalCycleRelativeOrder}. For every $(i_1,i_2)\in[n_1]\times[n_2]$, write $i_1\prec i_2$ if $\Omega(i_1) < \overline{\Omega}(i_2)$ and $i_1 \succ i_2$ otherwise.
\end{definition}
\begin{remark}
	To avoid ambiguity, we sometimes write $\prec_{ij}$ and $\succ_{ij}$ when working with more than two cycles.
\end{remark}

\begin{remark}
Theorem~\ref{theorem-unimodalCycleRelativeOrder} implies that for unimodal full cycles $\pi_1\in\unimodal_{n_1}$ and $\pi_2\in\unimodal_{n_2}$, there exist unique $J_1,J_2$ such that $\pi_1\oplus_{J_1}\pi_2$ and $\pi_2\oplus_{J_2}\pi_1$ are unimodal.
Observe that when $\pi_1 = \pi_2$ then $J_1 = J_2$, and there exists exactly one unimodal sum $\pi\oplus\pi$. 
The uniqueness of the set together with the convention that $\Omega(m_1)<\overline{\Omega}(m_2)$ allow us to omit the index $J$ in the sum $\pi_1\oplus_J\pi_2$ when it is unimodal. 
\end{remark}

The following result is a corollary of Theorem~\ref{theorem-unimodalCycleRelativeOrder} and determines the structure of unimodal sum $\pi\oplus\pi$ of unimodal cycle $\pi$ with itself. In particular, it determines the maximum of this sum.
\begin{corollary}\cite[Corollary 4]{Gannon-Unimodal}\label{corollary: CycleMultiplicityMaximum}
	Let $\pi \in \unimodal_n$ with the maximum at $m$ and let $J$ be the unique set such that $\pi\oplus_J\pi$ is unimodal.
	Then the following hold.
	\begin{enumerate}
		\item The set $J$ contains exactly one element from from $\left\{2k-1,k\right\}$ for every $1\leq k \leq n$. 
		\item If $n \equiv m \pmod{2}$ then $(\pi\oplus\pi)^{-1}(2n) = 2m$ and $(\pi\oplus\pi)^{-1}(2n)=2m-1$ otherwise. 
	\end{enumerate}
\end{corollary}

\begin{definition}\cite{Gannon-Unimodal}
	We say that $\pi\in\unimodal_{n}$ is {\em acute} if $m\equiv n \Pmod{2}$
	or equivalently $(\pi\oplus\pi)^{-1}(2n-1)<(\pi\oplus\pi)^{-1}(2n)$ (the maximal of subcycles of $\pi\oplus\pi$ run diagonally SW-NE like), and we say that $\pi\in\unimodal_{n}$ is {\em grave} otherwise, namely, the maxima of the subcycles of $\pi\oplus\pi$ run diagonally NW-SE (which implies that $2n$ occurs before $2n-1$.)
	
	For a unimodal permutation $\pi\in\unimodal(n)$ (not necessarily a cycle), let $\mathcal{A}(\pi)$, denote the set of all acute unimodal cycles $\sigma\in\unimodal_*$ for which exist $J\subseteq[n]$ such the shape of $\pi|_J$ that have shape $\sigma$. Similarly, $\mathcal{G}(\pi)$, will denote the grave shapes in $\pi$. Define $\mathcal{A}_>(\pi)$ as the multiset $\sigma\in\mathcal{A}(\pi)$, for which exists $J\subseteq[n]$, that have shape $\sigma$ and the maximum of $\pi|_J$ is larger than the maximum of $\pi$.
	$\mathcal{A}_<(\pi), \mathcal{G}_>(\pi), \mathcal{G}_<(\pi)$ are defined equivalently.  
\end{definition}

\begin{proposition}\cite[Proposition 5]{Gannon-Unimodal}\label{proposition:cycleMultisum}
	Let $\pi_i \in \unimodal_{n_i}$ be unimodal cycles for $i=1,2,3$ with maxima $m_i$ respectively. For $i<j$ let $\pi_{ij}$ be the unimodal sum $\pi_i\oplus\pi_j$ and $\prec_{ij},\succ_{ij}$ the order corresponding to $\prec$ Definition~\ref{definition:cyclePairElementsOrder} for $\pi_i$ and $\pi_j$.  
	\begin{enumerate}
		\item Assume that either $n_1 \prec_{12} n_2$ or $n_2 \succ_{23} n_3$. Then for every $a\in[n_1],b\in[n_2],c\in[n_3]$ we have both
		\begin{enumerate}
			\item $a \prec_{12} b$ and $b \prec_{23} c$ implies $a \prec_{13} c$.
			\item $a \succ_{12} b$ and $b \succ_{23} c$ implies $a \succ_{13} c$.
		\end{enumerate}
		\item There exist $J_i$ such that $\Omega_1(m_1) < \Omega_2(m_2) < \Omega_3(m_3)$ and $\pi_1\oplus_{J_1}\pi_{J_2}\oplus_2\pi_3$ is unimodal if and only if  either $n_1 \prec_{12} n_2$ or $n_2 \succ_{23} n_3$. Moreover, when such $J_i$ exist, they will be unique.
		\item Suppose that $\pi\in\unimodal_{k}$ and $\sigma\in\unimodal_{l}$. and let $I$ and $J$ be the unique $I$ and $J$ such that $\pi\oplus_I\sigma$ and $\sigma\oplus_J\pi$ are unimodal. 
		Then $k \prec_I l$ if and only if  $l \succ_J k$. 
	\end{enumerate}
\end{proposition}
\begin{remark}\label{remark:OPlusIndex2}
	Since $J_1,J_2,J_3$ such that $\Omega_1(m_1)<\Omega_2(m_2)<\Omega_3(m_3)$ and the sum $\pi=\pi_1\oplus_{J_1}\pi_2\oplus_{J_2}\pi_3$ is unimodal, we will write $\pi_1\oplus\pi_2\oplus\pi_3$ to denote this sum. 
\end{remark}

Part 3 of the preceding theorem says that if $\pi\in\Delta_k$ and $\sigma \in \Delta_l$ are unimodal full cycles, then either the maxima of $\pi\oplus_I\sigma$ and $\sigma\oplus_j\pi$ occur both either in $\Omega(\pi)$ or in $\Omega(\sigma)$. Therefore, we can define a total ordering on $\Delta_*$.
\begin{definition}\cite{Gannon-Unimodal}\label{definition:unimodalCycleOrder}
Part 3 of the previous proposition defines total order $\unimodless$ on $\unimodal_{*}$ as follows:
For $\pi_1\in\unimodal_{n_1}$ and $\pi_2\in\unimodal_{n_2}$, we say that $\pi_1\unimodless \pi_2$ if $n_1 \prec n_2$ in the unimodal sum $\pi_1\oplus\pi_2$.
\end{definition}

\begin{observation}\label{observation-UnimodalMaxima}
	Let $\pi_i\in\unimodal_{n_i}$ for $i=1,\dots,k$, then the terms of the sequence $n_1,\dots n_k$ are totally ordered by the relation $\prec$, defined as follows: $n_i \prec n_j$ when $i<j$ and $n_i \prec_{ij} n_j$ or $j<i$ and $n_j \succ_{ji} n_i$. Stated explicitly, $n_i \prec n_j$ if and only if, 
	\begin{enumerate}
		\item  $\pi_i\neq\pi_j$ and $\pi_i\unimodless\pi_j$ for $\pi_i\neq\pi_j$.
		\item  $\pi_i=\pi_j$, $i<j$ and $\pi$ is acute.
		\item  $\pi_i=\pi_j$, $j<i$ and $\pi$ is grave.
	\end{enumerate}
	Note that every two terms with different entities are treated as different objects. In addition, $\prec$ depends on $\pi_1,\dots,\pi_k$.
\end{observation}

The following theorem provides the conditions under which the unimodal sum of two unimodal permutations exists.
\begin{theorem}\cite[Theorem 6]{Gannon-Unimodal}\label{theorem:unimodalSumExistenceConditions}
	Let $\pi_i \in \unimodal(n_i)$ for $i=1,2$. Let $m_i = \pi_i^{-1}(n_i)$ and $\hat{\pi}_i$ be the shape of the cycle containing $m_i$ in $\pi_i$.
	Then 
	\begin{enumerate}
		\item if either $\mathcal{A}_>(\pi_1)\cap\mathcal{A}_>(\pi_2)$ or $\mathcal{G}_<(\pi_1)\cap\mathcal{G}_<(\pi_2)$ are nonempty, then there no unimodal sums of form $\pi_1\oplus\pi_2$ or $\pi_2\oplus\pi_1$.
		\item if instead $\hat{\pi}_1\in\mathcal{G}_<(\pi_2)$ or $\hat{\pi}_2\in\mathcal{A}_>(\pi_1)$, then no unimodal sum of form
		$\pi_1\oplus\pi_2$ with $m_1\prec m_2$ exists.
		\item otherwise, there is exactly one unimodal sum of form $\pi_1\oplus\pi_2$ with $m_1 \prec m_2$.
	\end{enumerate}
\end{theorem} 

\begin{definition}\cite{Gannon-Unimodal}
	When exists, the unique unimodal sum $\pi_1\oplus_J\pi_2$ with $\Omega(m_1)\prec \overline{\Omega}(m_2)$ will be denoted by $\pi_1\boxplus\pi_2$.
\end{definition}
\begin{definition}
	We say that $\pi_1$ and $\pi_2$ in $\unimodal(*)$ are cycle disjoint, if there is no $\sigma\in\unimodal_*$ such that both $\pi_1$ and $\pi_2$ have a subcycle of shape $\sigma$.
\end{definition}

\begin{corollary}\cite{Gannon-Unimodal}\label{corollary:UnimodalSumExistance}
	If both $\pi_1,\pi_2\in\unimodal(*)$ are shape disjoint, then both $\pi_1\boxplus\pi_2$ and $\pi_2\boxplus\pi_1$ exist.
\end{corollary}
\begin{observation}
	We may not omit braces when using $\oplus$ since $\oplus$ is not necessarily associative.
	For example, for $\pi_1,\pi_2,\pi_3\in\unimodal_*$ $\pi_1 \unimodless \pi_2 \unimodless \pi_3$, $\pi_2\oplus(\pi_1\oplus\pi_3)$ is defined, while $\pi_2\oplus\pi_1\oplus\pi_3$ is not. 
	For example, let $\pi_1 \unimodless \pi_2$ be two unimodal cycles, such that $\pi_1$ is grave. 
	Then $(\pi_1\oplus\pi_1)\oplus\pi_2$ = $\pi_1\oplus(\pi_2\oplus\pi_1)$, while $\pi_1\oplus(\pi_1\oplus\pi_2)$ does not exist by Part 2 of Theorem~\ref{theorem:unimodalSumExistenceConditions}.

	One case where braces can be omitted without risking ambiguity or having the sum not defined is the case where all $\pi_i$ are cycles of the same shape since the sum $\underset{k\ times}{\underbrace{\pi\oplus\dots\oplus\pi}}$ exists and must be unique by Theorem~\ref{theorem:unimodalSumExistenceConditions} and the relations $\prec_{ij}$ and $\succ_{ij}$.
\end{observation}

\begin{definition}
	For $k > 0$ and $\pi\in\unimodal_{*}$ write
	$$k\pi=\underset{k\,times}{\underbrace{\pi\oplus\dots\oplus\pi}}$$
	to be the unimodal sum $\pi\oplus\dots\oplus\pi$ where $\pi$ appears exactly $a$ times. We say that a unimodal permutation $\pi\in\unimodal(*)$ is {\em pure} if all the subcycles of $\pi$ have the same shape.
\end{definition}
\begin{remark}
	From now on, we use the notation $\pi=\pi_1\oplus\dots\oplus\pi_k$ where for every $i$, $\pi_i$ is pure and $\pi$ is defined. 
\end{remark}

The following lemma gives the sufficient and necessary conditions for the existence of the unimodal sum $\pi=\pi_1\oplus\dots\oplus\pi_k$ when $\pi_i$ is a full cycle for every $1\leq i \leq k$.

\begin{lemma}\label{lemma-unimodalSumWithoutBracesDefined}
	Let $\pi_1\in \unimodal_{n_1},\dots,\pi_k\in\unimodal_{n_k}$. Then there exists a unique unimodal $\pi$ such 
	$$\pi=\pi_1\oplus\dots\oplus\pi_k$$
	if and only if the sequence $n_1,\dots,n_k$ is unimodal when compared by $\prec$ as in Observation~\ref{observation-UnimodalMaxima}.
\end{lemma}
\begin{proof}
	If a unimodal permutation $\pi=\pi_1\oplus\dots\oplus\pi_k$ is unimodal permutation, then we have $\Omega_i(m_i) < \Omega_{i+1}(m_{i+1})$ and $\pi(\Omega_i(m_i))=\Omega_i(n_i)$. Therefore, $\Omega_1(n_1),\dots,\Omega_k(n_k)$ is a unimodal sequence, which is equivalent to $n_1,\dots,n_k$ being unimodal when compared by $\prec$.

	Uniqueness follows from Theorem~\ref{theorem-unimodalCycleRelativeOrder} and Proposition~\ref{proposition:cycleMultisum}.
	To see that, observe that if $\pi=\pi_1\oplus\dots\oplus\pi_k$ is unimodal, for every $1\leq i < j \leq k$ and $a\in[n_i], b\in[n_j]$, $\Omega_i(a) < \Omega_j(b)$ if and only if $a \prec_{ij} b$.
	But $\prec_{ij}$ does not depend on $\Omega_i$ and $\Omega_j$.
	Therefore, there is at most one way to order the elements of $\biguplus_{i=1}^k[n_i]$.
	(Note again, that for every $i\neq j$ the elements of $[n_i]$ and $[n_j]$ are treated as different entities).
	Therefore, if $\Omega_1,\dots,\Omega_k$, must be unique.

	We prove the existence of the unimodal sum $\pi=\pi_1\oplus\dots\pi_k$ by induction. The case of $k=1$ is trivial, and when $k=2$, there is a unique unimodal $\pi$ such that $\pi = \pi_1\oplus\pi_2$ by Theorem~\ref{theorem-unimodalCycleRelativeOrder}. Now assume that the statement holds for $k$ and show that it holds for $k+1$. Let $n_1,\dots,\pi_k, n_{k+1}$ ordered by $\prec$.
	By definition, $\mathcal{A}_{>}(\pi_{k+1})=\mathcal{G}_<(\pi_{k+1})=\emptyset$. 
	Therefore, by Theorem~\ref{theorem:unimodalSumExistenceConditions}, 
	$\pi\oplus\pi_{k+1}$ exists unless $\pi_{k+1} \in \mathcal{A}_>(\pi)$ 
	which is possible if and only if $\pi_{k+1}$ is acute and there exist 
	$1\leq j<k$ such that $n_j\succ n_i$ for every $i\neq j$ with $1 \leq i \leq k$ and $k \geq l>j$ such that $\pi_l=\pi_{k+1}$.
	But then we have $n_l \prec n_{k+1}$ and $n_l \prec n_j$, and $j < l < k+1$, which contradicting the unimodality of $n_1,\dots,n_{k+1}$.
\end{proof}

\section{Proof of Theorem 5.1.6}\label{section-Ch5FinalProof}
To complete the proof of the theorem, explore the connection between the maxima of unimodal full cycles and unimodal permutations that can be built from them by using the unimodal sum operation described in the previous section.
As noted in Observation~\ref{observation:unimodalCycleSumStructure}, every unimodal permutation can be represented as the unimodal sum of its cycle shapes.

We begin with the simplest case.

\begin{lemma}\label{lemma:purePermutationMaxL}
	Let $\pi\in\unimodal_*$ with the maximum $m$. Then, for every $k$,  the maximum of $k\pi$ is equal to
	\begin{enumerate}
		\item $k(m-1)+1$ when $\pi$ is grave.
		\item $km$ when $\pi$ is acute.
	\end{enumerate}
\end{lemma}
\begin{proof}
	The case for $k=1$ is trivial, and the case of $k=2$ is part 2 of Corollary~\ref{corollary: CycleMultiplicityMaximum}. Denote $m_k$ the maximum of $k\pi$.
	Assume that the statement is true for $k-1$.
	Then $k\pi = \pi\oplus((k-1)\pi)$ with $J=J_1$ where $J_1\dots J_k$ are the unique subsets of $[kn]$ such that $k\pi$ is unimodal and $J_i(m)<J_{i+1}(m)$ for every $1\leq i < k$. 
	
	If $\pi$ is grave, then $n \succ_1i$ holds for every $1<i\leq k$, which implies that maximum of $\pi\oplus((k-1)\pi)$ of is $\Omega(m)$.  By induction assumption $\pi_{k-1}=(k-1)(m-1)+1 $ by Corollary~\ref{corollary:UnimodalSumMax}, the maximum of $k\pi$ is 
	$m_k=m+m_{k-1}+1=m+(k-1)(m-1)-1=k(m-1)+1$ as desired.
	
	If $\pi$ is acute, then we have $n \prec_1i$ for every $1<i\leq k$, which implies that the maximum of $\pi\oplus((k-1)\pi)$ is $\overline{\Omega}(m_{k-1})$, and $m_k=m+m_{k-1}$ by Corollary~\ref{corollary:UnimodalSumMax}. By induction assumption, $m_{k-1}=(k-1)m$, therefore $m_k=m+(k-1)m=km$ as desired.
\end{proof}

We are interested in determining the maxima of permutations with a given number of subcycle shapes.
We start with a simple observation.
\begin{lemma}\label{lemma-orderedCyclePositioning}
	Let $\pi_i\in\unimodal_{n_i}$ for $i=1,\dots,k$ 
	such that $\pi=\pi_1\oplus\dots\oplus\pi_k$ is defined. Let $\sigma,\tau\in\unimodal_*$ such that $\sigma \unimodless \tau$.
	Then for any $1\leq i \leq k$ with $\pi_i=\tau$,  
	\begin{enumerate}
		\item If $\sigma$ is acute, there is at most one $j$ with $i < j$ such that $\pi_j=\sigma$.
		\item If $\sigma$ is grave, there is at most one $j$ with $j<i$ such that $\pi_j = \sigma$.		
	\end{enumerate} 
\end{lemma}

\begin{proof}
	By Lemma~\ref{lemma-unimodalSumWithoutBracesDefined} $\pi=\pi_1\oplus\dots\oplus\pi_k$ is defined if and only if $n_1,\dots,n_k$ ordered by $\prec$ is unimodal. If the sequence $\sigma$ is acute and there exist $j_2 > j_1$ such that $\pi_{j_1} = \pi_{j_2} = \sigma$, then we have $n_i \succ n_{j_1} \prec n_{j_2}$, which contradicts unimodality. Similarly, if $\sigma$ is grave and there exist $j_1 < j_2 < i$ such that $\pi_{j_1} = \pi_{j_2}$ then, we have $n_{j_1} \succ n_{j_2} \prec n_i$, which again contradicts the unimodality of $n_1,\dots,n_k$.
\end{proof}

\begin{definition}
	Let $a_1,\dots,a_l$ positive integers and $l$ distinct unimodal $\sigma_1,\unimodless \sigma_k\in\unimodal_*$ unimodal full cycles.
	Define
	$$\spanbox{a_1\sigma_1,\dots,a_l\sigma_l}=
	\left\{\pi\in\unimodal(*):
	\text{$\pi$ has exactly $a_i$ subcycles of shape $\sigma_i$ and no other}\right\}.$$
\end{definition}

\begin{corollary}\label{corollary:SpanMultiplicity}
Let $\sigma_1 \unimodless \dots \unimodless \sigma_l$ be $l$ different unimodal cycles, and $a_1,\dots,a_l$ postive integers.
Then the number of unimodal permutations in $\unimodal(*)$ with exactly $a_i$ subcycles of shape $\sigma_i$ is $2^{l-1}$.
\end{corollary}
\begin{proof}
	We prove the statement by induction. If $l=1$, there is one such permutation, namely $a_l\sigma_l$. Next, assume that the statement is true for $l$ and prove it for $l+1$. Let $\sigma_1,\dots,\sigma_{l+1}$ be $l+1$ unimodal full cycles, $a_1,\dots,a_{l+1}$ positive integers. By the induction assumption, there are $2^{l-1}$ unimodal permutations with $a_i$ subcycles of shape $\sigma_i$ for $2 \leq i \leq l+1$. Next, observe that if $\pi_1,\dots,\pi_k\in\unimodal_*$ and $\pi=\pi_1\oplus\dots\oplus\pi_k$ is defined, if $\pi_i < \pi_j$ and $\pi_i < \pi_k$, then we can't not have $j < i < k$, as it would imply $n_j \succ n_i \prec n_k$ contradicting Lemma~\ref{lemma-unimodalSumWithoutBracesDefined}. Intuitively, "largest shapes appear in the middle." Therefore, every $\sigma \in \spanbox{a_1\sigma_1,\dots,a_{l+1}\sigma_{l+1}}$ is of form $i\sigma_1 \oplus \pi \oplus j\sigma_2$ where $\pi\in\spanbox{a_2\sigma_2,\dots,a_{l+1}\sigma_{l+1}}$, and $i,j$ are nonnegative integers and $i+j=a_1$. 
	
 	By Lemma~\ref{lemma-orderedCyclePositioning}, we have $j\leq 1$ when $\sigma$ is acute, $i \leq 1$ when $\sigma$ is grave. 
	Therefore, for every $\pi=\pi_1\dots\oplus\dots\pi_m\in\spanbox{a_2\sigma_2,\dots,a_{l+1}\sigma_{l+1}}$ there are at most 2 unimodal permutations in $\spanbox{a_1\sigma_1,\dots,a_{l+1}\sigma_{l+1}}$ such that all the cycles of shape $\sigma_2,\dots,\sigma_{l+1}$ appear in exactly the same order they appear in $\pi$. 
	On the other hand, both unimodal sums are defined. 
	To see that, suppose that $\pi\in\spanbox{a_2\sigma_2,\dots,a_{l+1}\sigma_{l+1}}$ is a unimodal permutation, expressed as a unimodal sum of full cycles $\pi=\pi_1\oplus\dots\oplus\pi_m$. Let $n_i$ denote the maximum of $\pi_i$ and $n_\sigma$ the maximum of $\sigma_1$. Then, by Lemma~\ref{lemma-unimodalSumWithoutBracesDefined} the sequence $n_1,\dots,n_m$ ordered by $\prec$ is unimodal. But then, by definition of $\prec$ the sequences $$\underset{a_1\ times}{\underbrace{n_\sigma,\dots,n_\sigma}},n_1,\dots,n_m$$ and $$\underset{a_1-1\ times}{\underbrace{n_\sigma,\dots,n_\sigma}},n_1,\dots,n_m,n_\sigma$$ are unimodal and therefore both unimodal sums $$\underset{a_1\ times}{\underbrace{\sigma_1\oplus\dots\oplus\sigma_1}}\oplus\pi_1\oplus\dots\oplus\pi_m$$
	and 
	$$\underset{a_1-1\ times}{\underbrace{\sigma_1\oplus\dots\oplus\sigma_1}}\oplus\pi_1\oplus\dots\oplus\pi_m$$ are defined. When $\sigma_1$ is grave, the argument is symmetric.
\end{proof}

\begin{corollary}\label{corollary-boxSpanRepresentation}
	When $\sigma_1\unimodless\dots\unimodless \sigma_l$,
	Every element of $\spanbox{a_1\sigma_1,\dots,a_l\sigma_l}$ has unique representation, of $a_1\sigma_1\boxplus\pi$ or $\pi\boxplus a_l\sigma_1$ where $\pi\in\spanbox{a_2\sigma_2,\dots,a_l\sigma_l}$.
\end{corollary}
\begin{proof}
	By the proof of Corollary~\ref{corollary:SpanMultiplicity}, for every $\sigma\in\boxspan{a_1\sigma_1,\dots,a_l\sigma_l}$, there exist $\pi\in\boxspan{a_2\sigma_2,\dots,a_l\sigma_l}$,  $i,j$ such that $i+j=a_1$ and $i\sigma_1\oplus\pi\oplus j \sigma_1$. 
	If $\sigma$ is acute, then $j=0$ or $j=1$. 
	By definition of $\boxplus$, if $j=0$, then $$\sigma=i\sigma_1\oplus\pi=a_1\sigma_1\oplus\pi=a_1\sigma_1\boxplus\pi,$$
	and if $j=1$ then
	$$i\sigma=\sigma_1\oplus\pi\oplus j\sigma=
	(a_1-1)\sigma_1\oplus\pi\oplus\sigma_1=\pi\boxplus a_1\sigma_1.$$
	Similarily, when $\sigma$ is grave, we have either $i=1$ or $i=0$, 
	$$\sigma=\pi\oplus a_1\sigma_1=\pi\boxplus a_1\sigma_1$$ 
	when $i=0$ and 
	$$\sigma=\sigma_1\oplus\pi\oplus(a_1-1)\sigma_1$$
	when $i=1$.
\end{proof}

In the following proofs, we often rely on the observation stated below.
\begin{observation}\label{observation-ObviousUnimodPolynomial}
Let $A\subseteq\symm{n}$ be a set of unimodal permutations in $\symm{n}$.
Then unimodal maxima polynomial $U_A(x))$, just counts the number of elements in $A$ with given maximum and can be expressed as
$$\sum_{\pi\in A}x^{m_\pi},$$
where $m_\pi$ is the maximum of $\pi$.
\end{observation}

\begin{lemma}\label{lemma:cycleCombinationMaximum}
	Let $a_1,\dots,a_l>0$ and $\sigma_1\unimodless\dots\unimodless \sigma_l$ different full cycles, and let $m_i$ denote the maximum of $a_i\sigma_i$ for $1\leq i \leq l$ and let $m=m_1+\dots+m_l$.  Let $U_{\spanbox{a_1\sigma_1,\dots,a_l\sigma_l}}(x)$ be the unimodal maxima polynomial of $\spanbox{a_1\sigma_1,\dots,a_l\sigma_l}$.
	Then, 
	$$U_{\spanbox{a_1\sigma_1,\dots,a_l\sigma_l}}(x)=x^{m-l+1}(1+x)^{l-1}.$$
	In particular, for every $2<l$ $$1+x|U_{\spanbox{a_1\sigma_1,\dots,a_l\sigma_l}}(x)$$ and the quotient has non-negative coefficients.
\end{lemma}
\begin{proof}
	Let $m_1$,\dots$m_l$ be the maxima of $a_1\sigma_1,\dots,a_l\sigma_l$ respectively. We proceed by induction. When $l=1$, $U_\spanbox{a_1\sigma_1}(x)=x^m_1$ and result is trivial. When $l=2$, by Corollary~\ref{corollary-boxSpanRepresentation}, the elements of $\spanbox{a_1 \sigma_1,a_2 \sigma_2}$ are $a_1\pi_1\boxplus a_2\pi_2$ and $a_2\pi_2\boxplus a_1\pi_1$. By Observation~\ref{corollary:UnimodalSumMax}, and by definition of $\unimodless$, the maximum of $a_1\pi_1\boxplus a_2\pi_2$ is $m_1+m_2$ and the maximum of $a_2\pi_2\boxplus a_1\pi_1$ is $m_1+m_2-1$, and therefore
	$U_{\spanbox{a_1\sigma_1,a_2\sigma_2}}(x)=x^{m_1+m_2-1}(1+x)$ and the statement is true for $l=2$. 
	
	Now, assume that the statement is true for $l-1$. By Corollary~\ref{corollary-boxSpanRepresentation}, for every element of $\sigma\in\spanbox{a_1\sigma_1,\dots,a_l\sigma_l}$ exists a unique $\pi\in\spanbox{a_2\sigma_2,\dots,a_l\sigma_l}$ such that either
	$\sigma=a_1\sigma_1\boxplus \pi$ or $\sigma=\pi\boxplus a_1\sigma_1$.
	
	Denote $m_\pi,m_\sigma$ denote the maxima of $\pi,\sigma$ respectively.	By applying Observation~\ref{corollary:UnimodalSumMax} and the definition of $\unimodless$ again, 
	we have $m_\sigma=m_1+m_\pi$ for $\sigma=a_1\sigma_1\boxplus\pi$ and
	$m_\sigma=m_1+m_\pi-1$ for $\sigma=\pi\boxplus a_1\sigma_1$ is $m_1+m_\pi-1$.
	By Observation~\ref{observation-ObviousUnimodPolynomial},
	\[
	\begin{split}
	U_{\spanbox{a_1\sigma_1,\dots,a_l\sigma_l}}(x) &=
	\sum_{\sigma\in\spanbox{a_1\sigma_1,\dots,a_l\sigma_l}} x^{m_\sigma} 
	\\&=\sum_{\pi\in\spanbox{a_2\sigma_2,\dots,a_l\sigma_l}}x^{m_1+m_\pi-1}+x^{m_1+m_\pi}
	\\ &= x^{m_1-1}(1+x)\sum_{\pi\in\spanbox{a_2\sigma_2,\dots,a_l\sigma_l}}x^{m_\sigma}
	\\&=x^{m_1-1}(1+x)U_{\spanbox{a_2\sigma_2,\dots,a_l\sigma_l}}(x)
	\end{split}
	\]
	
	By the induction hypothesis,  
	$$U_{\spanbox{a_2\sigma_2,\dots,a_l\sigma_l}}(x)=x^{m^\prime-l+2}(1+x)^{l-2}$$ where $m^\prime=m_2+\dots+m_l$.
	Therefore, 
	\[ \begin{split}
	U_{\spanbox{a_1\sigma_1,\dots,a_l\sigma_l}}(x) &= x^{m_1-1}(1+x)x^{m^\prime-l+2}(1+x)^{l-2} \\
	&=x^{m_1+m^\prime-l-1}(1+x)^{l-1}=x^{m-l+1}(1+x)^{l-1}
	\end{split} \]
	as desired.
\end{proof}

\begin{observation}\label{observation:acutePure}
	When $\pi\in\unimodal(n)$ is a pure permutation, with maximum $m$, the maximum of $\pi\oplus\pi$ behaves precisely in the same way as the unique shape $\sigma$ of its subcycles.
	Namely, if $\sigma$ is acute, the maximum of $\pi+\pi$ is $2m$ when $\sigma$ is acute, and $2m-1$ when $\sigma$ is grave, parity of $n-m$. Hence, in the following discussion, we will refer to the pure permutations as {\em acute} or {\em grave}, even when they are not cycles.
\end{observation}

The involution $\phi$ from Corollary~\ref{corollary:involutionUnimodal} plays a crucial role in proving our theorem.
For a permutation $\pi$, let $m_\pi$ denote the minimum of $\pi$. 
Recall that $m_\phi(\pi)=m_\pi\pm 1$.

In the proof of Theorem~\ref{theorem-mainUrootPrimeVersion}, we construct unimodal permutations from $\pi$ and $\phi(\pi)$. The following lemma is crucial for understanding the behavior of maxima of unimodal permutations built this way.
\begin{lemma}\label{lemma:GraveAcuteSumMax}
	Let $\pi$ and $\sigma$ be two pure permutations, such that $m_\sigma=m_\pi+1$ holds.
	Then for every $k$ and $l$ the set of maxima of 
	$k\pi\boxplus l\sigma$ and $l\pi\boxplus k\sigma$ is
	\begin{enumerate}
		\item $\left\{(k+l)m_\pi + 1, (k+l)m_\pi \right\}$ when $\pi$ is acute.
		\item $\left\{(k+l)m_\pi + l-k +1,(k+l)m_\pi + l - k\right\}$ when $\pi$ is grave.
	\end{enumerate}
\end{lemma}
\begin{proof}
	By Observation \ref{corollary:UnimodalSumMax} the minima of $k_\pi\boxplus l\sigma $ and $l\sigma\boxplus k\pi$ are $m_{k\pi}+m_{l\sigma}$ and $m_{k\pi}+m_{l\sigma}-1$. Also, both numbers appear as the minimum of one of the sums, regardless of the order between $\sigma$ and $\pi$.

	If $\pi$ is acute, then $\sigma$ is grave, the maximum of $k\pi$ is $km_\pi$, and the maximum of $l\sigma$ is $l(m_\sigma-1)+1$. By substituting $(m_\sigma - 1)$ with $m_\pi$ and applying Corollary~\ref{corollary:UnimodalSumMax}, shows conclude the set of the maxima of $k\pi\boxplus l\sigma$ and $l\sigma \boxplus k\pi$ is $\left\{(k+l)m_\pi+1(k+l)m_\pi\right\}$.
	
	Similarily, if $\pi$ is grave then $\sigma$ is acute, the the maximum of $k\pi$ is $k(m_\pi-1)+1$ and the maximum of $l\sigma$ is $lm_\sigma$. By substituing $m_\sigma$
	with $m_\pi+1$ and applying
	Corollary~\ref{corollary:UnimodalSumMax}, we show that the set of maxima of $k\pi\boxplus l\sigma$ and $l\sigma \boxplus k\pi$ is $\left\{m_{k\pi}+m_{l\sigma},m_{k\pi}+m_{l\sigma}-1\right\}$, as desired.
\end{proof}

We are ready to complete the proof of Theorem~\ref{theorem-mainUrootPrimeVersion}.

\begin{theorem}
	For a prime $p$, a cyclic descent extension for $\uroot{p}{n}$ exists if and only if $p \mid n$.
\end{theorem}
\begin{proof}
	Let $p$ be a prime number. We prove that if $p|n$, then a cyclic descent extension exists for $\uroot{p}{n}$. The proof of the converse, namely that if a cyclic descent extension exists for $\uroot{p}{n}$, then $p|n$ was given in Lemma~\ref{lemma:easyMainConjecturePart}.
    
    By the given assumption, $p$ divides $n$, and hence $q=\frac{n}{p}$ is an integer. For any unimodal element $\pi \in \uroot{p}{n}$, there exists a unique set of unimodal cycles $\sigma_1,\dots,\sigma_l$ such that $\sigma_i\in\unimodal_{p}\cup\unimodal_1$ for $1\leq i \leq l$, and positive numbers $a_1,\dots a_l$, such that $\pi$ has exactly $a_i$ subcycles of shape $\sigma_i$ and no other subcycle shapes. In other words, $\pi\in\spanbox{a_1\sigma_1,\dots,a_l\sigma_l}$.
    
    Let $S$ be the set of all possible combinations of ${{a_1\sigma_1,\dots,a_l\sigma_l} }$ such that $${\spanbox{a_1\sigma_1,\dots,a_l\sigma_l}\cap\uroot{n}{p}\neq \emptyset}.$$ Then
    $$\uroot{p}{n}\cap \unimodal(n) = \biguplus_{{a_1\sigma_1,\dots,a_l\sigma_l}\in S}\spanbox{a_1\sigma_1,\dots,a_l\sigma_l}.$$
	Obviously, 
	$$U_{\uroot{p}{n}}(x)=U_{\uroot{p}{n}\cap\unimodal(n)}(x).$$
	Therefore, 
	$$U_{\uroot{p}{n}}(x)=
	\sum_{\{a_1\sigma_1,\dots,a_l\sigma_l\} \in S}U_{\spanbox{a_1\sigma_1,\dots,a_l\sigma_l}}(x).$$
	By Lemma~\ref{lemma:cycleCombinationMaximum}, when $l\geq 2$. $U_{\spanbox{a_1\sigma_1,\dots,a_l\sigma_l}}(x)$ is divisible by $1+x$, with the coefficients of the quotient being non-negative.

	o finish the proof, we need to handle the case when $l=1$, as $U_{\spanbox{q\sigma}}=x^{m_{q\sigma}}$ is not divisible by $1+x$. To address this issue, we use the involution $\phi$ from Corollary~\ref{corollary:involutionUnimodal}.
    
    By using $\phi$, we can partition $\uroot{p}{p}\cap\unimodal(p)$ into pairs ${\pi,\phi(\pi)}$. Without loss of generality, we assume that $m_{\phi(\pi)}=m_\pi+1$. If this is not the case, we can replace $\pi$ with $\phi(\pi)$.
    
    Let $S_\pi$ denote all possible pairs ${a_\pi\pi,a_{\phi(\pi)}\phi(\pi)}$ with $a_\pi+a_{\phi(\pi)}=q$, where $q=\frac{n}{p}$. If $\pi=e$ or $\phi(\pi)=e$, we replace $e$ with $p(1)$ when definitions require it. Observe that $e=p(1)$ can be treated as a regular acute full cycle for calculation purposes by Observation~\ref{observation:acutePure}.
    
    Note that the parities of $m_\pi$ and $m_{\phi(\pi)}$ are different, meaning that $\pi$ is acute if and only if $\phi(\pi)$ is grave.
    
    We proceed to calculate
    $$\sum_{\{a_\pi\pi,a_{\phi(\pi)}\phi(\pi)\}\in S_\pi}U_{\spanbox{a_\pi\pi,a_{\phi(\pi)}\phi(\pi)}}(x).$$

	If $\pi$ is acute and $\phi(\pi)$ is grave, then the maxima of $\spanbox{a_\pi\pi,a_{\phi(\pi)}\phi(\pi)}$ for  every $1\leq a_\pi \leq q-1$  are $qm_\pi+1$ and $qm_\pi$, with maximum of $q\pi$ being $qm_\pi$ and the maximum of $q\phi(\pi)$ being $qm_\pi+1$, by Lemma~\ref{lemma:GraveAcuteSumMax}.
	Hence,
	$$\sum_{\{a_\pi\pi,a_{\phi(\pi)}\phi(\pi)\}\in S_\pi}U_{\spanbox{a_\pi\pi,a_{\phi(\pi)}\phi(\pi)}}(x)=q(x^{qm_\pi}+x^{qm_\pi+1})$$ which is divisible by $1+x$ with quotient $qx^{qm_\pi}$.
	If $\pi$ is grave then the maxima of $q\pi$ is $q(m_\pi-1)+1$ and the maximum of $q\phi(\pi)$ is $q(m_\pi+1)$.
	For every $1<a_\pi<q-1$, by Lemma~\ref{lemma:GraveAcuteSumMax} the maxima of 
	$\spanbox{a_\pi\pi,a_{\phi(\pi)}\phi(\pi)}$ are $(a_\pi+a_{\phi(\pi)})m_\pi+q-2a_\pi+1$ and $(a_\pi+a_{\phi(\pi)})m_\pi+q-2a_\pi$.
	Thus, we have
	\begin{enumerate} 
		\item $U_{\spanbox{q\pi}}=x^{qm_\pi-q+1}$,
	 	\item $U_{\spanbox{q\phi(\pi)}}(x)=x^{q(m_\pi+1)}$
	 	\item $U_{\spanbox{a_\pi\pi,a_{\phi(\pi)}\phi(\pi)}}=x^{qm_\pi+q-2a_\pi+1}
	 			+x^{qm_\pi+q-2a_\pi}$ for $1\leq a_\pi \leq q-1$.
	\end{enumerate}
 	Therefore,
 	\[
 	\begin{split}
 		&\sum_{\spanbox{a_\pi\pi,a_{\phi(\pi)}\phi(\pi)}\in S_\pi}U_{\spanbox{a_\pi\pi,a_{\phi(\pi)}\phi(\pi)}}(x)=\\
 		&x^{qm_\pi-q+1}+\sum_{a_\pi=1}^{q-1}(x^{qm_\pi+q-2a_\pi+1} +x^{qm_\pi+q-2a_\pi})+x^{q(m_\pi+1)}=\\
 		&x^{qm_\pi-q+1}\sum_{i=1}^{2q-1}x^{i}		
	\end{split}
 	\]
 	which is also divisible by $1+x $, with the quotient having negative coefficients.
 	Hence, the terms of $$U_{\uroot{p}{n}}(x)$$ coming from $S_\pi$ sum up to polynomial
 	divisible by $1+x$ with non-negative coefficients in the quotient, and other terms
 	are form $U_{\spanbox{a_1\sigma_1,\dots,a_l\sigma_l}}(x)$ with $l\geq 2$ and are also divisible by $1+x$ with non-negative coefficients in the quotient.
 	Thus $U_{\uroot{p}{n}}(x)$ is divisible by $1+x$, and the quotient has non-negative coefficients. The result follows.
\end{proof}
\cleardoublepage

\end{document}